\documentclass[a4paper,11pt,french,leqno]{smfart}
 
  \usepackage[francais]{babel}
\usepackage{hyperref}
\usepackage{amsmath}
\usepackage{amssymb}
\usepackage{amsthm}
\usepackage{amsopn}
\usepackage{amscd}
\usepackage{fullpage}
\usepackage[utf8]{inputenc}
\usepackage{graphics, xcolor}
\usepackage{textcomp} 
\usepackage[all]{xy}
\usepackage{lscape}
\usepackage{url}
\usepackage{verbatim}
\usepackage{mathrsfs}

\newcommand{\scrB}{\mathscr{B}}
\newcommand{\scrC}{\mathscr{C}}

\newcommand{\scrE}{\mathscr{E}}

\newcommand{\scrH}{\mathscr{H}}

\newcommand{\scrP}{\mathscr{P}}

\newcommand{\frI}{\mathfrak{I}}

\newcommand{\frS}{\mathfrak{S}}

\newcommand{\fra}{\mathfrak{a}}

\newcommand{\frg}{\mathfrak{g}}
\newcommand{\frh}{\mathfrak{h}}

\newcommand{\frl}{\mathfrak{l}}

\newcommand{\frqqq}{\mathfrak{q}}

\newcommand{\fru}{\mathfrak{u}}

\newcommand{\bbC}{\mathbb{C}}

\newcommand{\bbN}{\mathbb{N}}

\newcommand{\bbQ}{\mathbb{Q}}
\newcommand{\bbR}{\mathbb{R}}

\newcommand{\bbZ}{\mathbb{Z}}

\newcommand{\caA}{\mathcal{A}}
\newcommand{\caB}{\mathcal{B}}
\newcommand{\caC}{\mathcal{C}}

\newcommand{\caG}{\mathcal{G}}
\newcommand{\caH}{\mathcal{H}}

\newcommand{\caL}{\mathcal{L}}
\newcommand{\caM}{\mathcal{M}}

\newcommand{\caQ}{\mathcal{Q}}
\newcommand{\caR}{\mathcal{R}}
\newcommand{\caS}{\mathcal{S}}
\newcommand{\caT}{\mathcal{T}}
\newcommand{\caU}{\mathcal{U}}

\newcommand{\caX}{\mathcal{X}}

\def\U{\mathbf{U}}

\newcommand{\coker}{\mathrm{coker} }

\newcommand{\GL}{\mathbf{GL}}

\newcommand{\Hom}{\mathrm{Hom}}

\newcommand{\Ind}{\mathrm{Ind}}

\newcommand{\tr}{\mathrm{Tr}\, }
\newcommand{\Id}{\mathrm{Id}}

\renewcommand{\Ind}{\mathrm{Ind}}

\newcommand{\SL}{\mathbf{SL}}

\newcommand{\SO}{\mathbf{SO}}
\newcommand{\Sp}{\mathbf{Sp}}
\newcommand{\Or}{\mathbf{O}}

\newcommand{\Tr}{\mathrm{Tr}}

\newcommand{\Trans}{\mathrm{Trans}}
\newcommand{\St}{\mathbf{St}}
\newcommand{\Speh}{\mathbf{Speh}}
\newcommand{\sgn}{\mathbf{sgn}}
\newcommand{\Triv}{\mathbf{Triv}}

\newcommand{\Std}{\mathbf{Std}}
\newcommand{\Res}{\mathrm{Res}}

\theoremstyle{plain}
\newtheorem{thm}{Théorème}[section]
\newtheorem{lemme}[thm]{Lemme}
\newtheorem{cor}[thm]{Corollaire}
\newtheorem{prop}[thm]{Proposition}

\theoremstyle{definition}
 \newtheorem{defi}[thm]{Définition}

\newtheorem{rmq}[thm]{Remarque}

\def \dem {\noindent \underline{\sl Démonstration}. }

\begin{document}

\numberwithin{equation}{subsection}

\title{Paquets d'Arthur des groupes classiques et unitaires}

\author{Nicol\'as Arancibia}
\address{Institut Mathématique de Jussieu} 
\email{nicolas.arancibia@imj-prg.fr}
 
  \author{Colette Moeglin}
 \address{CNRS, Institut Mathématique de Jussieu } 
 \email{colette.moeglin@imj-prg.fr}

  \author{David Renard  }
 \address{Centre de Mathématiques
  Laurent Schwartz,  Ecole Polytechnique} 
\email{david.renard@polytechnique.edu}

\begin{abstract} 
Soit $G=\mathbf G(\bbR)$ le groupe des  points réels d'un groupe algébrique connexe réductif quasi-déployé défini sur $\bbR$. 
Supposons de plus que $G$ soit un groupe  classique  (symplectique, spécial orthogonal ou unitaire).
Nous montrons que les  paquets de représentations   irréductibles unitaires et   cohomologiques  définies par
 Adams et  Johnson en 1987 co\"incident
avec ceux  definis plus récemment par  J. Arthur dans son travail sur la  classification 
du  spectre automorphe discret des  
groupes classiques  (C.-P. Mok pour les  groupes unitaires). Pour cela, nous calculons  le  transfert endoscopique 
des  distributions stables 
sur  $G$ supportées par ces paquets vers le groupe  $\GL_N$ tordu  en termes de modules standard  et nous montrons 
 qu'il est égal à la  trace  tordue prescrite par  Arthur. 
\\

{\bf  Abstract.}\, ---\, 
Let $G=\mathbf G(\bbR)$ be the group of real points of a quasi-split connected reductive  algebraic group
defined over $\bbR$. Assume furthermore  that $G$ is a classical group (symplectic, special orthogonal or unitary).
We show that the packets of  irreducible unitary  cohomological representations defined by Adams and Johnson in 1987 coincide
with  the ones defined recently by J. Arthur in his work on the classification of the discrete automorphic spectrum of 
classical groups (C.-P. Mok for unitary groups). For this, we compute the endoscopic transfer of the stable distributions
on $G$ supported by these packets to twisted $\GL_N$ in terms of standard modules and show that it coincides
with the twisted trace prescribed by Arthur. 
\end{abstract}

\date{\today}

\thanks{Le troisième auteur a bénéficié d'une aide de  l'agence nationale de la recherche 
ANR-13-BS01-0012 FERPLAY}

\maketitle

\section{Introduction}

En \cite{Art13}, Arthur a donn\'e une description des repr\'esentations automorphes de carr\'e 
int\'egrable des groupes classiques  quasi-d\'eploy\'es, le  cas des 
groupes unitaires ayant été traité ensuite dans \cite{mok}.   Dans les deux cas cette  description 
 ram\`ene la situation à celle  des groupes g\'en\'eraux   lin\'eaires, les groupes ${\mathbf{GL}}_N$ 
 sur un corps de nombres.  Pour ces groupes ${\mathbf{GL}}_N$, c'est le  th\'eor\`eme de multiplicit\'e 
 un fort qui est utilis\'e en plus  de la description des repr\'esentations automorphes de carr\'e int\'egrable: 
  la situation en les places non ramifi\'ees pour une repr\'esentation automorphe induite 
  de repr\'esentations de carr\'e int\'egrable, d\'etermine la situation en toutes les places. En \cite{Art13} et ses
   g\'en\'eralisations, un th\'eor\`eme du m\^eme ordre est d\'emontr\'e. Notons $N$ la dimension de la repr\'esentation
    naturelle du $L$-groupe du groupe classique ou unitaire consid\'er\'e, groupe que l'on note $\mathbf G$; 
    ce groupe  est suppos\'e  quasi-d\'eploy\'e et d\'efini sur un corps de nombres que l'on ne nomme pas, car il n'appara\^{\i}t pas 
     dans l'article.  Le th\'eor\`eme de multiplicit\'e un fort est presque vrai pour $\mathbf G$,
      c'est-\`a-dire que la situation aux places non ramifi\'ees d'une repr\'esentation automorphe
       de carr\'e int\'egrable de $\mathbf G$ d\'etermine uniquement une repr\'esentation automorphe de 
       ${\mathbf{GL}}_N$,  obtenue comme induite de repr\'esentation de carr\'e int\'egrable. 
 Ceci est \'etabli via  la stabilisation de la formule des traces tordue et la fonctorialit\'e de Langlands en les places non ramifi\'ees.
   En les autres places, 
 la repr\'esentation de ${\mathbf{GL}}_N$ d\'etermine \`a son tour une repr\'esentation virtuelle 
 de longueur finie, qu'il vaut mieux  voir   comme une somme (avec a priori des multiplicit\'es) de 
 repr\'esentations irr\'eductibles.    Le r\'esultat complet d\'emontr\'e en \cite{Art13}
 et \cite{mok}, nous dit quand un produit (restreint) en toutes les places d'une de ces repr\'esentations irr\'eductibles
 est effectivement une repr\'esentation automorphe de carr\'e int\'egrable de $\mathbf G$ et avec quelle multiplicit\'e 
 elle intervient.  Les  multiplicit\'es se décomposent en deux types de facteurs : des  multiplicit\'es de nature globale,
   calcul\'ees explicitement en \cite{Art13} et \cite{mok},  et des multiplicit\'es locales devant elles  \^etre
calcul\'ees par des propri\'et\'es de transfert endoscopique. Ce qui n'apparaît pas 
dans les travaux sus-cités est  la description explicite  de ces repr\'esentations irr\'eductibles 
et de ces multiplicit\'es locales aux places ramifiées.  Le cas des  places finies a été réglé dans \cite{MoeM1}, 
 et en particulier en ces places les multiplicit\'es locales sont un.
 Dans cet article, nous nous intéressons aux places archimédiennes.
 
Notons donc maintenant  $\mathbf G$ un groupe classique
quasi-d\'eploy\'e défini sur $\bbR$, et $G=\mathbf G(\bbR)$ le groupe de ses points réels. Supposons pour simplifier que 
$\mathbf G$ soit symplectique ou spécial orthogonal, nous ferons quelques remarques sur le cas des groupes 
unitaires à la fin de cette introduction. Soit 
$\Std_G: {}^LG\rightarrow \GL_N(\bbC)$ la représentation standard de son $L$-groupe. 
 Soit $\psi_G : W_\bbR\times\SL_2(\bbR) \rightarrow {}^LG$ un paramètre d'Arthur (voir section \ref{ArtP}), et notons 
 $\psi=\Std_G\circ \psi_G$, que l'on peut voir comme un paramètre d'Arthur pour le groupe $G_N:=\GL_N(\bbR)$. 
 Comme nous l'avons dit au début de cette introduction, le spectre automorphe discret 
est bien compris pour $\GL_N$,
 et en particulier, on sait associer à $\psi$ une représentation
irréductible unitaire de $G_N$, composante locale d'une représentation automorphe de carré intégrable, et que nous notons
$\Pi_\psi$ (Section \ref{PaqArtGN}). Le paramètre $\psi$ provenant d'un paramètre pour un groupe classique, la représentation
$\Pi_\psi$ est auto-duale, c'est-à-dire qu'elle est stable sous l'action de l'automorphisme extérieur $\theta_N$ de $\GL_N$
 (Section \ref{espacetordu}). 
   Posons $G_N^+=G_N \rtimes \langle \theta_N \rangle$ et  $\widetilde G_N:=G_N^+\setminus G_N$. 
 Ce dernier est un espace tordu pour $G_N$, et la représentation $\Pi_\psi$  admet des  extensions comme représentations de 
 l'espace tordu  $\widetilde G_N$, deux d'entre-elles étant égales à un multiple scalaire près.
  Lorsqu'on fixe une donnée de Whittaker (Section \ref{NormExt}) pour le groupe  $G_N$, 
 on peut distinguer grâce à la théorie des fonctionnelles de Whittaker l'une de ces  extensions, que l'on note 
 $\Pi_\psi^+$. On note $\Tr_{\theta_N}(\Pi_\psi)$
la trace de  $\Pi_\psi^+$, c'est une distribution sur  $\widetilde G_N$.

 D'autre part, Arthur définit dans \cite{Art13} un paquet $\Pi_{\psi_G}$ attaché au paramètre $\psi_G$, c'est-à-dire 
 un ensemble fini de représentations unitaires irréductibles de $G$, ayant toute même caractère infinitésimal
 (facilement déterminé par $\psi_G$). Ce paquet est caractérisé par le fait que certaines distributions sur $G$, qui sont des combinaisons
 linéaires de caractères-distribution des membres du paquet doivent vérifier  
 des identités de transfert endoscopique. Il y a deux types de telles identités. Les premières sont  les identités endoscopiques
 ordinaires qui relient les combinaisons linéaires mentionnées ci-dessus à d'autres,  qui vivent sur un groupe endoscopique
 $H$ (qui est un produit de deux groupes classiques) de $G$, via la théorie de Langlands-Shelstad. 
  Les secondes sont les identités endoscopiques tordues, où le groupe $G$ et la représentation standard
 $\Std_G$ de son $L$-groupe sont des données endoscopiques elliptiques pour le groupe tordu $(G_N,\theta_N)$. 
 Ce sont ces identités qui sont cruciales dans ce travail. Le point de départ est que le paquet d'Arthur $\Pi_{\psi_G}$ est le support 
 d'une distribution stablement invariante sur $G$, que nous notons $\Theta_{\Pi_{\psi_G}}^{st}$, et qui est donc une certaine
 combinaison linéaire bien déterminée  de caractères-distribution des représentations dans $\Pi_{\psi_G}$.
 La théorie de l'endoscopie tordue (\cite{KS},\cite{Shel}, \cite{mezo}) transfère de telles distributions stablement invariantes sur $G$ vers 
 des distributions sur $\widetilde G_N$, invariantes sous l'action par conjugaison  de $G_N$,  par une application 
 $\Trans_G^{\widetilde G_N}$ (transfert spectral). Cette application est ici totalement déterminée, et pas seulement à une constante multiplicative
 près, et ceci grâce à la donnée de Whittaker qui a été fixée sur $G_N$.
  L'identité endoscopique tordue, fondamentale pour nous, s'écrit alors
 \begin{equation}\label{IdEndTor}
 \Trans_G^{\widetilde G_N}(\Theta_{\Pi_{\psi_G}}^{st})=\Tr_{\theta_N}(\Pi_\psi).
 \end{equation}
 Le problème, comme nous l'avons mentionné, c'est que dans ces travaux d'Arthur, les membres des
 paquets  $\Pi_{\psi_G}$ ne sont pas identifiés dans une classification connue, par exemple dans la classification de Langlands.
 Il en est de même de la forme exacte des combinaisons linéaires des caractères-distributions formées
 à partir de ce paquet et qui entrent dans les identités endoscopiques. Or, pour un paramètre d'Arthur $\psi_G$
 comme ci-dessus, d'autres constructions de paquets ont été proposés.
 Dans \cite{ABV}, ceci est fait 
 dans un cadre totalement g\'en\'eral (c'est-\`a-dire plus g\'en\'eral que les groupes classiques). Un
  paquet $\Pi_{\psi_G}^{ABV}$ y est défini,  celui-ci est le support d'une distribution stablement invariante
 $\Theta_{\Pi_{\psi_G}^{ABV}}^{st}$
  et  il est conjectur\'e que cette construction co\"{\i}ncide dans les cas \'etudi\'es avec celles de \cite{Art13} .
 D'autre part, il est établi dans \cite{ABV} que les paquets qui y sont définis vérifient les identités endoscopiques ordinaires. 
 Ainsi, pour montrer la conjecture, il suffit d'établir que ces derniers satisfont aussi à l'identité endoscopique tordue (\ref{IdEndTor}), 
 c'est-à-dire que l'on peut remplacer  le terme de gauche de cette identité par $ \Trans_G^{\widetilde G_N}(\Theta_{\Pi_{\psi_G}^{ABV}}^{st})$.
 Les identités endoscopiques ordinaires ne jouent donc aucun rôle dans cet article et nous n'en dirons rien de plus.
 Notons au passage que les constructions de \cite{ABV} utilisent des invariants géométriques sophistiqués des representations
 (les \og cycles caractéristiques\fg), et que ceux-ci sont incalculables en pratique, ce qui fait que les membres 
 des  paquets  $\Pi_{\psi_G}^{ABV}$ ne sont pas identifiés eux non plus dans une classification connue.
 
 Pour un certain type de paramètres d'Arthur $\psi_G$, auxquels nous allons nous référer dans cette article sous la dénomination
 peut-être abusive de \og paramètres d'Adams-Johnson\fg, 
 une autre construction avait été proposé antérieurement à \cite{ABV}
 par Adams et Johnson (\cite{AdJo}). La propriété fondamentale de ces paramètres est que le caractère infinitésimal qui 
 leur est associé est entier et régulier.
  Là encore, un paquet $\Pi_{\psi_G}^{AJ}$ est défini, et celui-ci est 
 le support d'une distribution stablement invariante $\Theta_{\Pi_{\psi_G}^{AJ}}^{st}$.
 Cette fois, les membres du paquet sont bien identifiés, ce sont des représentations cohomologiques unitaires, c'est-à-dire des modules
 $A_\frqqq(\lambda)$ de Vogan-Zuckerman (\cite{VZ}, voir  \cite{KV} Chapter 5).
  Il a \'et\'e v\'erifi\'e par les auteurs de \cite{ABV}  que leur  construction coïncide avec celle   de \cite{AdJo}. 
 Si les paramètres d'Adams-Johnson ne sont qu'un type particulier  de paramètre  d'Arthur,  ils jouent néanmoins
 un rôle important dans tous les problèmes liés aux relations entre formes automorphes et cohomologie des variétés.

 Nous pouvons maintenant énoncer le résultat principal de cet article :
 \begin{thm}\label{MainThm} Soit $\psi_G$ un paramètre d'Adams-Johnson du groupe classique $G$. Alors  
 \begin{equation}\label{MR}
 \Trans_G^{\widetilde G_N}(\Theta_{\Pi_{\psi_G}^{AJ}}^{st})=\Tr_{\theta_N}(\Pi_\psi).
 \end{equation}
 En conséquence, le paquet d'Arthur $\Pi_{\psi_G}$ est égal au paquet d'Adams-Johnson $\Pi_{\psi_G}^{AJ}$.
 \end{thm}
 Donnons maintenant quelques détails sur la façon d'établir ce résultat. 
 L'idée est de passer par le transfert spectral tempéré établi par Mezo (\cite{mezo}) et précisé dans l'appendice.
 Soit $\phi_G:\, W_\bbR\rightarrow {}^LG$ un paramètre de Langlands. Pour les groupes réels, Langlands  (\cite{Langl1})
 a défini le paquet correspondant $\Pi_{\phi_G}$. Lorsque le paquet est tempéré, la somme des caractères-distribution des membres 
 du paquet est une distribution stablement invariante sur $G$. Ce n'est plus vrai si le paquet n'est pas tempéré, mais dans ce cas, on définit
 le pseudo-paquet $\widetilde \Pi_{\phi_G}$ dont les éléments sont les représentations standard, {\sl i.e.} des induites paraboliques
 de représentations essentiellement tempérées dont les éléments de $\Pi_{\phi_G}$ sont les sous-représentations de Langlands
 (dans cet article, nous utilisons la version de la classification de Langlands en termes de  sous-représentations des représentations standard,
  et non de quotients comme il est le plus souvent l'usage).
 Alors la somme des caractères-distribution des membres 
 du pseudo-paquet $\widetilde \Pi_{\phi_G}$ est une distribution stablement invariante notée 
 $ \Theta_ {\widetilde \Pi_ {\phi_G}}$. Soit $\phi=\Std_G\circ \phi_G$,  paramètre de Langlands
 pour $G_N$, et soit $\widetilde \Pi_\phi$ la représentation standard $\theta_N$-stable de $G_N$ 
qui lui est associée.  Le résultat de Mezo est alors que 
\begin{equation}
 \Trans_G^{\widetilde G_N}( \Theta_ {\widetilde \Pi_ {\phi_G}})=\Tr_{\theta_N}(\widetilde \Pi_\phi).
 \end{equation}
 Mezo démontre cette égalité à une constante multiplicative près, et c'est cette ambiguïté qui est levée dans l'appendice.
 Dans leur article, Adams et Johnson ont aussi écrit explicitement $\Theta_{\Pi_{\psi_G}^{AJ}}^{st}$
 comme combinaison linéaire de $  \Theta_ {\widetilde \Pi_ {\phi_G}}$, ceci vient des résolutions
de Johnson pour les $A_\frqqq(\lambda)$ de \cite{Joh1}. Disons que l'on a 
\begin{equation}\label{AJID} \Theta_{\Pi_{\psi_G}^{AJ}}^{st}   =\sum_{\phi_G}   a_{\phi_G}  \, \widetilde \Theta_{\phi_G}^{st}. \end{equation}
 Les  coefficients $a_{\phi_G}$ sont en fait des signes et la somme est bien entendu à support fini.
 Grâce au résultat de Mezo, le membre de gauche s'écrit donc  (toujours en posant $\phi=\Std_G\circ \phi_G$)
 \begin{equation}\label{AJID2}\Trans_G^{\widetilde G_N}( \Theta_{\Pi_{\psi_G}^{AJ}}^{st} )  =\sum_{\phi_G}   a_{\phi_G}  \, 
 \Tr_{\theta_N}(\widetilde \Pi_{\phi}), 
 \end{equation}
et l'on est ramené à démontrer que 
 \begin{equation}\label{AJID3}\ \sum_{\phi_G}   a_{\phi_G}  \, 
 \Tr_{\theta_N}(\widetilde \Pi_{\phi})=\Tr_{\theta_N}(\Pi_\psi). 
 \end{equation}
Notre approche est de le faire par récurrence sur la longueur $r$ d'une décomposition 
$\psi=\oplus_{i=1,\ldots,r}\psi_i $
du paramètres $\psi$ en \og paramètres élémentaires\fg. Le cas où $r=1$, c'est-à-dire que $\psi$ est élémentaire, 
se sépare en deux sous-cas disjoints. Le premier est celui où $\psi_G$ est le paramètre d'Arthur d'un caractère 
quadratique de $G$. On montre par un argument global que le paquet d'Arthur  $\Pi_{\psi_G}$ est un singleton (Section \ref{dim1}).
L'identité (\ref{MR}) est alors établie par Arthur dans \cite{Art13}, et (\ref{AJID3}) est alors une conséquence qui nous servira
dans l'étape de récurrence. Le second cas élémentaire  est celui où $\psi$ est le paramètre 
d'une représentation de Speh $\theta_N$-stable de $G_N$, $N$ pair. 
Il est traité dans la section \ref{ParSpehE}, à partir de la formule pour la trace tordue
d'une telle représentation établie dans la section \ref{caractor}. Pour obtenir cette formule, nous partons du fait qu'une représentation de 
Speh est un cas particulier de $A_\frqqq(\lambda)$ de Vogan-Zuckerman, et qu'elle
admet donc une résolution de Johnson par des modules standard, c'est-à-dire que l'on a un complexe exact de la forme
\begin{equation}\label{introresSpeh}
0\rightarrow \Speh\rightarrow X_0\rightarrow X_1\rightarrow \cdots \rightarrow X_{\ell_{max}}\rightarrow 0
\end{equation}
où les $X_i$ sont des sommes directes de représentations standard $X(w)$, indexées par $w\in \frS_{N/2}$, et 
$X_i$ est la somme des $X(w)$ pour $w\in \frS_{N/2}$ de longueur $i$ (la longueur usuelle dans le groupe symétrique), et 
bien sûr $\Speh=\Pi_\psi$ est une représentation de Speh.
 Une telle résolution implique immédiatement une 
formule pour le caractère-distribution  de cette représentation comme somme alternée des caractère des $X_i$. Mais ici, ce n'est pas le 
caractère-distribution de $\Speh$ que nous voulons calculer, mais sa trace tordue.  Nous montrons, en suivant les argument de 
Johnson que le complexe   (\ref{introresSpeh}) de représentations de $G_N$ peut être muni d'une structure
de complexe de représentations du groupe $G_N^+$. Les représentations standard $\theta_N$-stables
du complexe, les seules qui vont contribuer à la trace tordue, sont indexées par les involutions $w$ du groupe
$ \frS_{N/2}$. Elles sont donc étendues au groupe $G_N^+$, cette extension étant l'une des deux possibles $X(w)^\pm$, mais 
qui n'est pas nécessairement celle, disons $X(w)^+$, déterminée par la donnée de Whittaker.
Nous calculons alors par un argument combinatoire le signe par lequel l'extension de $X(w)$, $w$ involution de  
$\frS_{N/2}$, diffère de $X(w)^+$. Ceci nous donne la contribution de ce  $X(w)$ à la trace tordue. Nous introduisons
une fonction longueur sur l'ensemble  $\frI_{N/2}$ des involutions de  $\frS_{N/2}$, appelée $\theta$-longueur
et notée $\ell_\theta$, et au final, la trace tordue de la Speh s'écrit
\begin{equation}\label{introTTSpeh}
\Tr_{\theta_N}(\Speh) =\sum_{w\in \frI_{N/2}} (-1)^{\ell_\theta(w)} \Tr_{\theta_N} (X(w)).
\end{equation}
Il reste donc à comparer les termes extrêmes de 
 \begin{equation}\label{AJID4}\ \sum_{\phi_G}   a_{\phi_G}  \, 
 \Tr_{\theta_N}(\widetilde \Pi_{\phi})=\Tr_{\theta_N}(\Pi_\psi)=\Tr_{\theta_N}(\Speh) =\sum_{w\in \frI_{N/2}} (-1)^{\ell_\theta(w)} \Tr_{\theta_N} (X(w)). 
 \end{equation}
Ceci se fait en exhibant une bijection entre l'ensemble des $\phi_G$ intervenant dans le terme de gauche avec $\frI_{N/2}$
et en vérifiant que si $\phi_G$ correspond à $w$, alors $X(w)=\widetilde \Pi_{\phi}$ et 
$(-1)^{\ell_\theta(w)}=a_{\phi_G}$ (rappelons que d'après Adams-Johnson
$a_{\phi_G}$ est bien un signe). Cette comparaison passe par les représentations des groupe unitaires $\U(b,c)$ avec $b+c=N/2$.
Les représentations de $\U(b,c)$ sont reliées à celles de $G$ par le foncteur d'induction cohomologique de Zuckerman, 
et à celles de $G_N$ par changement de base vers $\GL_{N/2}(\bbC)$ suivi d'une induction cohomologique de 
$\GL_{N/2}(\bbC)$ vers $G_N$.

Une fois démontré le résultat principal pour les paramètres élémentaires, on passe à la démonstration du cas général 
par récurrence sur la longueur de la décomposition de $\psi$ en paramètres élémentaires. Dans une telle décomposition, 
on a au plus un seul paramètre élémentaire du premier type décrit ci-dessus, et donc dans l'étape de récurrence, on 
suppose qu'on ajoute un paramètre élémentaire $\psi_1$ du second type (paramètre d'une Speh) à un paramètre $\psi'$. 
 L'idée est la suivante : de même que le transfert endoscopique ordinaire
commute à l'induction parabolique, le transfert endoscopique (spectral) tordu de $G$ vers $\widetilde G_N$ commute à un 
certain foncteur d'induction, un peu délicat à décrire exactement malheureusement, mais dont voici l'ingrédient principal.
Supposons que $N=N_1+N'$, avec $N_1=2n_1$, et soit $\Speh$ une représentation de Speh $\theta_{N_1}$-stable de $G_{N_1}$.
Soit $M$ le sous-groupe de Levi standard de $G_N$ isomorphe à $G_{N_1}\times G_{N'}$. Posons $\theta_M=\theta_{N_1}\times \theta_{N'}$, 
que l'on voit comme un automorphisme d'ordre 2 de $M$, puis $M^+=M\rtimes \langle \theta_M  \rangle$ et 
$\widetilde M=M^+\setminus M$. On considère la  catégorie des représentations $\pi_M^+$ de $\widetilde M$ dont la restriction $\pi_M$
à $M$ est de longueur finie,  et telle que tous  les sous-quotients irréductibles de $\pi_M$ sont des produits tensoriels 
de la forme $\pi_1\otimes \pi'$, avec $\pi_1$ une représentation irréductible de $G_{N_1}$, sous-représentation de 
Langlands de l'une des représentations $X(w)$, $w\in \frS_{n_1}$, apparaissant dans la résolution de Johnson de $\Speh$, et $\pi'$
une représentation irréductible de $G_{N'}$ dont on suppose seulement 
que le caractère infinitésimal $\lambda'$ est fixé, et vérifie une certaine condition ({\sl cf.} hypothèse  (\ref{hyplemmegen}), 
le paramètre de la Speh doit être 
grand devant  $\lambda'$, et cette condition  va être vérifiée lorsqu'on part d'un paramètre d'Adams-Johnson).
 On montre tout d'abord (lemme \ref{irredIndpiM}) que si $\pi_1$ et $\pi'$ vérifient ces hypothèse, 
alors $\Ind_P^{G_N}(\pi_1\otimes \pi')$ est irréductible, $P$ étant ici le sous-groupe parabolique standard de facteur de Levi $M$.
Grâce aux propriétés fines  des opérateurs d'entrelacement (Section \ref{opeent}),
 on construit un foncteur d'induction
$\Ind_{\widetilde M}^{\widetilde G_N}$ de cette catégorie de représentations de $\widetilde M$ vers les représentations de  
$\widetilde G_N$, ayant la propriété suivante : 
la donnée de Whittaker pour $G_N$ induit une donnée de Whittaker pour  $M$. Supposons que $\pi_M$ soit 
une représentation irréductible $\theta_M$-stable de $M$ dans la catégorie décrite ci-dessus, et soit 
 $\pi_M^+$ extension à $\widetilde M$ distinguée par la donnée de Whittaker. Alors $\Ind_{\widetilde M}^{\widetilde G_N}(\pi_M^+)$
est l'extension de  $\Ind_{M}^{G_N}(\pi_M)$ distinguée par la donnée de Whittaker. 
Ces résultats sont les plus techniques, mais aussi les plus prometteurs de l'article, car on espère les réutiliser dans un contexte
plus général que celui des paramètres d'Adams-Johnson, et ils faisaient partie de la thèse du premier auteur.

Pour le changement de base des  groupes unitaires, les résultats analogues avaient été obtenus par Johnson dans \cite{Joh2}, en s'appuyant
sur le transfert  tempéré établi par Clozel  \cite{CloBC} dans ce cas (et dont le résultat de Mezo est une généralisation).
Nous reprenons rapidement notre démarche pour montrer qu'elle s'adapte sans difficulté  aux groupes unitaires
(c'est même plus simple techniquement).

\medskip 

Décrivons brièvement le contenu de l'article. 
Les premières sections sont purement formelles, elles servent à rappeler les notations et résultats de la littérature essentiels
à la suite de l'article. La section \ref{LetA} concerne les paramètres de Langlands et d'Arthur, et la classification de Langlands.
On y rappelle la définition des distributions stables attachées aux paramètres de Langlands (lorsque le paramètre n'est pas tempéré,
on introduit les \og pseudo-paquets\fg, constitués de représentations induites à partir de tempérées).
La section \ref{GLNnot} présente la classification du dual admissible des groupes $\GL_N(\bbR)$.
On y rappelle aussi la définition de la représentation attachée à un paramètre d'Arthur (dans le cas de $\GL_N$, les paquets
d'Arthur sont des singletons), le cas basique étant celui des représentations de Speh.
La section \ref{BBUpq} introduit la classification de Beilinson-Bernstein pour les groupes unitaires.
Comme nous l'avons expliqué, ceci va nous permettre de faire le lien entre  les termes de gauche et de droite de (\ref{AJID3}).
La section   \ref{torduGLN}   concerne les repr\'esentations de l'espace tordu $\widetilde G_N$ (ou de manière presque équivalente, du 
groupe $G_N^+$) et de  leur normalisation \`a l'aide des mod\`eles de Whittaker. Ceci permet de définir précisément
le terme de droite  de (\ref{MR}) (autrement seulement défini à une constante multiplicative près).
La section \ref{caractor} est consacrée au calcul de la trace tordue des représentations de Speh auto-duales et la 
 section \ref{IndTor} au foncteur d'induction parabolique \og tordu \fg.
La section \ref{ParAJ} 
rappelle les résultats d'Adams et Johnson dans un contexte général, en rappelant en particulier la définition des paquets
$\Pi_{\psi_G}^{AJ}$, de la combinaison linéaire stable des caractères-distributions  $\Theta_{\Pi_{\psi_G}^{AJ}}^{st}$
et son expression en terme de pseudo-paquets obtenue à partir des résolutions de Johnson. 
La section \ref{GCET} introduit les groupes classiques orthogonaux et symplectiques considérés, leur représentation standard, 
 l'application de transfert spectral $\Trans_G^{\widetilde G_N}$ et détaille la forme des paramètres d'Adams-Johnson pour ces groupes.
La section \ref{principal} est consacrée à la démonstration du résultat principal.
La section \ref{GrUnit} adapte l'énoncé du  résultat principal et sa démonstration aux groupes unitaires, en parcourant rapidement
les modifications superficielles à effectuer. Enfin l'appendice précise le résultat de Mezo en levant l'ambiguïté de celui-ci.
Il présente un intérêt indépendant du reste de l'article et le résultat est obtenu par des méthodes globales.

\medskip 

Le résultat principal de cet article est utilisé comme hypothèse dans de nombreux travaux. Citons en particulier
\cite{BMM1}, \cite{BMM2}, \cite{ChRe}, \cite{ChLa}, \cite{Kott},  \cite{Morel}, \cite{Tai}.
\footnote{Dans  \cite{BMM1}, \cite{BMM2}, \cite{ChLa}, l'hypothèse en question permet simplement de simplifier
la démonstration de certains résultats.}

\medskip 

Dans certains cas particuliers, le premier auteur avait d\'ej\`a obtenu dans sa thèse le r\'esultat principal  (sous 
l'hypoth\`ese que les représentations de Speh  vivent dans des $\GL_{2n}(\bbR)$ avec  $n\leq 4$). 
Une partie  des résultats de cet article sont issus -tels quels ou sous une forme ayant évoluée avec le temps- de cette thèse 
du premier auteur \cite{Ara},  en  particulier dans le chapitre  7.

\tableofcontents

\section{Paramètres de Langlands et d'Arthur}\label{LetA}
\subsection{Paramètres de Langlands et pseudo-paquets}
\label{LangL}
Nous renvoyons le lecteur à \cite{Bor} pour plus de détails sur les objets introduits dans cette section. 

Soient $\mathbf{G}$ un groupe algébrique réductif connexe défini sur $\bbR$, $\widehat G$ 
son dual de Langlands, et ${}^LG= \widehat G \rtimes W_\bbR $ son $L$-groupe, où bien sûr 
$W_\bbR$ est le groupe de Weil de $\bbR$.

Le groupe $\widehat G$ agit par conjugaison sur l'ensemble des paramètres de Langlands $\phi:\, W_\bbR\rightarrow {}^LG$,  et l'on note 
$\Phi(G)$ l'ensemble de ces classes de conjugaison. On note $\Phi_{\mathrm{temp}}(G)$  l'ensemble
 des classes de conjugaison de paramètres de Langlands d'image bornée.

Le théorème de classification de Langlands donne l'existence d'une partition de l'ensemble $\Pi(G)$ des classes d'équivalence  
 de représentations irréductibles du groupe $G:=\mathbf G(\bbR)$  
\begin{equation}\label{LanglClass} \Pi(G)=\coprod_{\phi \in \Phi(G)}  \Pi_\phi \end{equation}
en  $L$-paquets (ou paquets de Langlands) $\Pi_\phi$. Le point essentiel est bien entendu
les propriétés de cette partition. Donnons-en quelques unes.
Tous les éléments d'un paquet ont même caractère
infinitésimal. Pour les représentations tempérées, on a 
\begin{equation}\label{LanglClassTemp}\Pi_{\mathrm{temp}}(G)=\coprod_{\phi \in \Phi_{\mathrm{temp}}(G)}  \Pi_\phi. \end{equation}
 Lorsque $\pi$ est une représentation de longueur finie
 de $G$, notons $\Theta_\pi$ son caractère : c'est une distribution invariante sur $G$. 
 Pour tout  $\phi \in \Phi_{\mathrm{temp}}(G)$, notons 
\begin{equation}\label{StabTemp} \Theta_{\Pi_\phi}= \sum_{\pi \in \Pi_{\phi}} \Theta_\pi. \end{equation}
C'est une  distribution {\sl stablement invariante} sur $G$. 
Pour ce qui concerne cette notion, nous renvoyons par exemple à \cite{Shelstad79} et \cite{Boua}.
 Ceci n'est plus vrai pour un paquet non tempéré. Soit $\phi \in \Phi(G)$, non nécessairement tempéré.  
 Toute représentation $\pi\in \Pi_\phi$
est obtenue dans la classification de Langlands comme l'unique sous-représentation irréductible
d'une représentation standard $I(\pi)$, induite parabolique d'une représentation 
essentiellement tempérée.

\begin{defi}\label{pseudopaquets} Soit $\phi\in \Phi(G)$ un paramètre de Langlands.
Appelons \og pseudo-paquet\fg, 
et notons $\widetilde \Pi_\phi$, l'ensemble
des $I(\pi)$ pour $\pi \in \Pi_\phi$. Posons :
\begin{equation}\label{StablePsP}  \Theta_ {\widetilde \Pi_ {\phi}}= \sum_{\pi \in \Pi_{\phi}} \Theta_{I(\pi)}. 
\end{equation}
\end{defi}

Il est bien connu (\cite{Shelstad79}, \cite{AdJo},  Lemma 4.3) 
 que  $ \Theta_ {\widetilde \Pi_ {\phi_G}}$ est une  distribution stablement invariante sur $G$. 
Remarquons que lorsque $\phi$ est tempéré, on a  $ \Theta_ { \Pi_ {\phi}}= \Theta_ {\widetilde \Pi_ {\phi}}$.

\subsection{Paramètres d'Arthur}\label{ArtP}
Les notations sont les m\^emes que dans la section précédente.

\begin{defi}\label{ArtPar}
Un paramètre d'Arthur pour  $G$ est un morphisme de groupes continu
\begin{equation*}
\psi: \, W_\bbR \times \SL_2(\bbC) \longrightarrow {}^ L G
\end{equation*}
tel que 
\begin{itemize}
\item[(i)] la restriction de $\psi$ à  $W_\bbR$ est un paramètre de Langlands tempéré,
\item[(ii)]  la restriction de  $\psi$ à $\SL_2(\bbC)$ est algébrique.
\end{itemize}

\medskip

Le groupe $\widehat G$ agit par conjugaison sur l'ensemble des paramètres d'Arthur, et l'on note 
$\Psi(G)$ l'ensemble de ces classes de conjugaison. On identifie  $\Phi_{\mathrm{temp}}(G)$ à l'ensemble
 des paramètres d'Arthur de restriction triviale à $\SL_2(\bbC)$.
\end{defi}

A tout paramètre d'Arthur  $\psi$, on associe un paramètre de Langlands 
\begin{equation}\label{ArtLang}
\phi_\psi: \, W_\bbR \longrightarrow {}^ L G, \qquad w \mapsto \psi(w,
  \left(  \begin{matrix}  \vert w \vert^{\frac{1}{2}} & 0\\ 
0 &   \vert w \vert^{-\frac{1}{2}}   \end{matrix} \right)) , 
\end{equation}
où $w\mapsto \vert w \vert$ est le morphisme de groupe de $W_\bbR$ dans $\bbR^\times_+$ 
défini de la manière suivante. Rappelons que $W_\bbR$ peut-être vu comme le groupe engendré
par $\bbC^\times$ et un élément $j$ tel que $j^2=-1$ ({\sl cf.} \cite{ABV}, Definition 5.2). On pose alors 
 $\vert j\vert =1$ et $\vert z\vert =z\bar z$ si $z\in \bbC^\times$. 

Dans \cite{Art84}, \cite{Art89}, J. Arthur conjecture  (pour des groupes définis
sur un corps local $F$) l'existence de paquets $\Pi_\psi$  attachés aux 
paramètres $\psi \in \Psi(G)$,   devant posséder certaines propriétés, énoncées dans l'introduction. 

Pour $F$ archimédien,  Adams et Johnson ont proposé  pour une certaine classe de paramètres $\psi$ 
(les \og paramètres d'Adams-Johnson\fg, voir section \ref{ParAJ}) une définition 
de paquets  $\Pi^{AJ}_\psi$ possédant aussi un certain nombre des propriétés requises, en particulier
l'existence d'une distribution stable $\Theta_{\Pi_\psi^{AJ}}^{st}$ explicite, et la compatibilité à l'endoscopie standard.
Comme il a été expliqué dans l'introduction, J.Arthur pour les groupes classiques, puis C-P. Mok 
pour les groupes unitaires ont donné une définition des paquets  ayant toutes les propriétés voulues, 
en les caractérisant par les identités endoscopiques standard (transfert spectral vers leurs groupes endoscopiques 
standard) et par les identités endoscopiques tordues venant du fait que ces groupes peuvent être vu comme
faisant partie d'une donnée endoscopique tordue d'un groupe général linéaire. Notre but dans cet article est de montrer
que les paquets définis par Adams-Johnson et ceux définis par Arthur et Mok coïncident (pour un même paramètre, bien entendu)
 dans les cas où ils sont tous deux définis. Pour cela, il suffit donc de montrer que les premiers
 satisfont les  identités endoscopiques tordues avec les groupes généraux linéaires.

\section{ Le groupe $G_N=\GL_N(\bbR)$ et ses représentations}\label{GLNnot}
Dans cette section, nous introduisons des notations concernant le groupe général linéaire et ses représentations.
 Soit $N$ un entier positif. Notons 
  $\mathbf{G}_N$ le groupe général linéaire $\GL_N$ défini sur le corps des réels
  (avec la convention que $\mathbf{G}_0$ est le groupe trivial), 
    $G_N=\GL_N(\bbR)$  le groupe de ses points réels, et   $\widehat G_N$ le groupe dual de Langlands 
    (isomorphe à  $\GL_N(\bbC)$). 

 Soit   $ N_1,\ldots ,N_r \in \bbN^\times $   tels que $\displaystyle \sum_{i=1}^r N_i=N$. 
Le sous-groupe $M=M_{N_1,\ldots ,N_r}$ des matrices diagonales par blocs de taille respective
$N_1,\ldots, N_r$, isomorphe à 
$G_{N_1} \times G_{N_2}\times ...\times G_{N_r}$
est un sous-groupe de Levi standard de $G_N$, et le sous-groupe parabolique 
 $P=P_{N_1,\ldots, N_r}$  contenant $M$ et le sous-groupe de  Borel des  matrices triangulaires supérieures 
 est un sous-groupe parabolique standard de radical unipotent  $N=N_{N_1,\ldots, N_r}$.
Pour tout  $1\leq i\leq r$, soit   $\pi_i$ une représentation de $G_{N_i}$ de longueur finie. 
Nous notons alors  $\pi_1\times \pi_2\times ...\times \pi_r$ la représentation obtenue par induction parabolique 
(normalisée) à partir de la représentation $\pi_1\otimes \pi_2\otimes ...\otimes \pi_r$ de $M$ 
relativement au sous-groupe parabolique $P$.

 Notons $\nu=\nu_N$ le caractère $g\mapsto \vert \det g \vert$ de $G_N$. Pour toute représentation 
$\pi$ de $G_N$ et tout $s \in \bbC$,  notons de manière abrégé  $\pi\nu^s$ le produit tensoriel $\pi\otimes \nu^s$ .
 
 \medskip 
 Soit $\frh_{d,N}$  la sous-algèbre de Lie  des matrices diagonales de $\caM_N(\bbC)$, 
 identifiée naturellement à $\bbC^N$, et de même pour son dual $\frh_{d,N}^*$.
Via l'isomorphisme d'Harish-Chandra, un caractère infinitésimal pour $G_N$ est donné
par un élément de $\frh_{d,N}^*$ et donc  par un élément $\underline \lambda=
(\lambda_1,\ldots, \lambda_N) \in \bbC^N$, où plutôt
 par une orbite de tels éléments sous l'action du groupe de Weyl, ici identifié au groupe $\frS_N$. 
 Un tel caractère infinitésimal est entier si  les $\lambda_i-\lambda_j$ sont entiers, et régulier
 si les $\lambda_i$ sont distincts.

\medskip 

Nous allons maintenant rappeler  la classification de Langlands de 
$\Pi(\GL):=\coprod_{N\in \bbN}\Pi(G_N)$  en termes de  représentations irréductibles de 
$G_1$ et $G_2$. Par souci  d'alléger un peu la terminologie, nous appelons \og séries discrètes unitaires \fg\, 
 ce que l'on devrait appeler \og séries discrètes modulo le centre \fg \,  et \og séries discrètes \fg \, les
  représentations   obtenues  par produit tensoriel  d'une série discrète unitaire et d'un caractère.  

Le groupe   $G_1=\GL_1(\bbR)\simeq \bbR^\times $ est abélien, et ses représentations 
irréductibles sont de la forme 
\begin{equation}\label{ParG1} \gamma(s,\epsilon): \bbR^\times \rightarrow \bbC^\times, \quad x\mapsto 
\vert x \vert^s \sgn(x)^\epsilon, \quad 
(s \in  \bbC), \; (\epsilon \in \{0,1\}). \end{equation}
Un tel caractère est unitaire si et seulement si $s\in i\bbR$. En général, posons
\begin{equation}\label{exposantG1}
e(\gamma_(s,\epsilon))=\Re e(s).
\end{equation}
Considérons deux caractères $\gamma(s_i,\epsilon_i)$, $i=1,2$ comme ci-dessus .
La série principale $\gamma(s_1,\epsilon_1)\times \gamma(s_2,\epsilon_2)$ de $G_2$ 
est  réductible si et seulement si  $s_1-s_2=n\in \bbZ^\times$ et  
$\epsilon_1+\epsilon_2\equiv n+1 \mod 2$ et dans ce cas l'un des deux sous-quotient irréductible
de celle-ci est une  série discrète  que l'on note $\delta(s_1,s_2)$. Les équivalences entre ces séries discrètes font 
l'on peut imposer    $n=s_2-s_1\in \bbN^\times$. 
A équivalence près, les séries discrètes de $G_2$ sont donc les 
\begin{equation}\label{ParG2}
\delta(s_1,s_2), \qquad s_1, s_2 \in \bbC,\;  s_2-s_1 \in \bbN^\times .\end{equation}
\bigskip 
Une telle série discrète est unitaire si et seulement si $s_1+s_2\in i\bbR$.  En général, posons
\begin{equation}\label{exposantG2}
e(\delta(s_1,s_2))=\frac{\Re e(s_1+s_2)}{2}.
\end{equation}
Si   $\tau$ est l'un des caractères $\gamma(\epsilon,s)$ ou bien l'une des séries discrètes $\delta(s_1,s_2)$, on a
\[ \tau= \nu^{e(\tau)}  \tau^u, \]
où $\tau^u$ est de même type que $\tau$, mais unitaire.

\begin{defi}\label{OrdStd} Supposons que pour tout $i=1, \ldots,l$, on se donne une série discrète  $\tau_i$
qui est  soit l'un des séries discrètes (\ref{ParG1}) de $\GL_1(\bbR)$ ou bien l'une des séries discrètes  (\ref{ParG2}) de $\GL_2(\bbR)$
On dit que  $\underline \tau=(\tau_1,\ldots,\tau_l)$ est écrit dans {\sl un ordre standard} si 
$ e(\tau_1) \leq  \ldots \leq e(\tau_l)$ 
et dans  {\sl un ordre standard inverse } si $ e(\tau_1) \geq  \ldots \geq e(\tau_l)$. 
\end{defi}

Le théorème de classification de Langlands s'énonce alors ainsi : 

\begin{thm}\label{Lgclass} Soit   $\underline \tau=(\tau_1,\ldots,\tau_l)$ comme dans la définition ci-dessus, 
   écrit dans  un ordre standard, alors  :
\begin{itemize}
 \item[(i)] La  représentation 
$ X(\underline \tau)= \tau_1 \times \ldots  \times \tau_l $
possède une unique sous-représentation irréductible  $\bar X(\underline \tau)$, apparaissant avec multiplicité un  dans
la suite de  Jordan-Hölder de   $X(\underline \tau)$. 
Cette représentation $\bar X(\underline \tau)$, est aussi l'unique quotient 
irréductible de la représentation 
\[ \widetilde X(\underline \tau):= \tau_l\times \tau_{l-1} \times \ldots \times \tau_2 \times \tau_1. \]


 \item[(ii)]  L'application $ \underline \tau \mapsto \bar X(\underline \tau) $
réalise une bijection entre  l'ensemble des   ensembles  avec multiplicités de séries discrètes et $\Pi(\GL)$. 
\end{itemize}

Les représentations $X(\underline \tau)$, $\widetilde X(\underline \tau)$ sont appelées représentations standard.
\end{thm}

\begin{rmq} \label{defstdgln} Les représentations  $X(\underline \tau)$ et 
$\bar X(\underline \tau)$ ne dépendent pas du choix d'un ordre standard sur les $\tau_i$ (parmi tous les ordres standard possibles).
Il en est de même de  $\widetilde X(\underline \tau)$ et de l'ordre standard inverse choisi pour l'écrire.
C'est un cas particulier d'un théorème de B. Speh rappelé plus loin (Théorème \ref{SpehIrrStd}).
Ceci nous permet d'adopter les notations suivantes : si  $\underline \tau=(\tau_1,\ldots,\tau_l)$, on note 
 \[ \displaystyle \times_{i=1,\ldots,r}\;  \tau_i =\tau_1 \times \ldots \times \tau_l \]
et 
 \[  \times_{i=1,\ldots,r}^\rightarrow  \;  \tau_i , \qquad \times_{i=1,\ldots,r}^\leftarrow  \;  \tau_i\]
le produit des $\tau_i$ obtenu en les permutant pour les mettre dans  un ordre standard 
et un ordre standard inverse respectivement. 
\end{rmq}

\begin{rmq}
Le caractère infinitésimal de $\bar X(\underline \tau)$ est donné par $(\lambda_1,\ldots ,\lambda_N)$
si et seulement si l'ensemble avec multiplicité $\{ \lambda_1,\ldots ,\lambda_N \}$ 
 est constitué des $ s_i$  pour les $i$ tels que  $ \tau_i=\gamma(s_i,\epsilon_i)$ et
 des  $ s_{i,1},s_{i,2} $  pour les $i$ tels que  $\tau_i=\delta(s_{i,1},s_{i,2})$.   
\end{rmq}

\subsection{Classification de Langlands pour $\GL_N(\bbR)$ avec $L$-groupe }

Pour $\mathbf{G}=\GL_N$, les paquets de Langlands sont des singletons, et 
(\ref{LanglClass}) constitue donc une classification de $\Pi(G_N)$. Ce fait nous autorise à noter de la m\^eme
manière le paquet associé à un paramètre $\phi$, et l'unique représentation qu'il contient : $\Pi_\phi$. Il en est de même du pseudo-paquet
$\widetilde \Pi_\phi$ attaché à $\phi$.

Comme le $L$-groupe est ici un produit direct, ${}^LG_N=\widehat G_N \times W_\bbR $, on peut voir simplement
un paramètre de Langlands  $\phi$ comme une représentation  de dimension $N$ de $W_\bbR$.
Cette représentation est complètement réductible. 
Les représentations irréductibles de $W_\bbR$ vont paramétrer les séries discrètes. 
Elles sont donc de dimension $1$  ou $2$. Rappelons tout ceci brièvement.

Les  représentations irréductibles de  $W_\bbC\simeq \bbC^\times$ sont de dimension $1$ puisque  $\bbC^\times$ est
 abélien. Elles sont paramétrées par les couples     
$(s_1, s_2)\in \bbC\times \bbC$ avec  $s_1-s_2\in \bbZ $ de la manière suivante : 
\begin{equation}\label{Chisn} \chi_{s_1,s_2}(z)=\vert z \vert^{s_1+s_2} \left( \frac{z}{\vert z \vert}\right)^{s_1-s_2} =
z^{s_1}\bar z^{s_2}
\end{equation}
La représentation  $\chi_{s_1,s_2}$ est unitaire si  $s_1+s_2 \in i\bbR$.

Les  représentations irréductibles de  $W_\bbR$ sont facilement obtenues à partir des  $\chi_{s,n}$ par 
la théorie de Mackey. Posons 
\begin{equation}\label{Vs1s2} V(s_1,s_2)=\Ind_{\bbC^\times}^{W_\bbR} (\chi_{s_1,s_2}).\end{equation}
Si  $s_1-s_2\neq 0$, $V(s_1,s_2)$ est une représentation irréductible de  $W_\bbR$.
De plus  $V(s_1,s_2)\sim V(s_1',s'_2)$ si et seulement si  $\{s_1,s_2\}=\{s'_1,s'_2\}$.

Si  $s_1-s_2=0$, $V(s_1,s_2)$ est  réductible. Notons  respectivement  $\Triv$ and $\sgn$
les caractères   de  $W_\bbR$  obtenus par relèvement des caractères de   $W_\bbR/(W_\bbR)_0\simeq \bbZ/2\bbZ$.
Notons aussi  $\vert\; \vert^s$ le caractère de  $W_\bbR$ donné explicitement sur les générateurs par 
$ z\mapsto \vert z \vert^s, \, (z\in \bbC^\times), \, j\mapsto 1$.  
Alors  $V(s,s)$ se décompose comme  
$ V(s,s)\simeq \Triv \otimes \vert\; \vert^{2s} \oplus \sgn \otimes \vert\; \vert^{2s} $.  
Posons 
 \begin{equation}W(s,\epsilon)= \begin{cases} \Triv \otimes \vert\; \vert^s \text{ si } \epsilon=0\\ 
       \sgn \otimes \vert\; \vert^s \text{ si } \epsilon=1
\end{cases}.        \end{equation}

\begin{prop} \label{IrrWR} Les représentations irréductibles de  $W_\bbR$ sont à équivalence près : 
\begin{itemize} 
\item[(i)] les  représentations  $W(s,\epsilon)$, $s\in \bbC$, $\epsilon \in \{0,1\}$, de dimension 1, 

\item[(ii)] les représentations $V(s_1,s_2)$, $s_1,s_2\in \bbC, \, s_2-s_1\in \bbN^\times $, de dimension 2.
\end{itemize}
Elles sont unitaires si et seulement si $s\in i\bbR$ (cas  $(i)$) et $s_1+s_2\in i\bbR$ (cas  $(ii)$).
\end{prop}
\begin{prop}
Les séries discrètes de $G_1$ et $G_2$ sont respectivement en bijection avec les représentations irréductibles 
de $W_\bbR$ de dimension $1$ et $2$,  cette correspondance étant 
\begin{align*}
\gamma_{s,\epsilon} &\leftrightarrow  W(s,\epsilon), \qquad  s\in \bbC, \epsilon \in \{0,1\},\\
\delta(s_1,s_2)&\leftrightarrow V(s_1,s_2), \qquad s_1,s_2 \in \bbC, s_2-s_1\in \bbN^\times,\; 
\end{align*}
où $\gamma_{s,\epsilon}$ est défini en (\ref{ParG1}) et $\delta(s_1,s_2)$  est défini en (\ref{ParG2}).
  \end{prop}

La bijection entre $\Pi(G_N)$ et $\Phi(G_N)$ s'en déduit alors de la manière suivante. 
Soit $\phi\in \Phi(G_N)$, considéré  comme une représentation de dimension $N$ de  $W_\bbR$.
Cette représentation se décompose en une somme de représentations irréductibles de $W_\bbR$, 
 écrivons ceci $\phi=\oplus_{i=1}^r \phi_i$ avec chaque $\phi_i$ équivalente à l'une des représentations
 de $W_\bbR$ dans $\GL_1(\bbC)$ ou $\GL_2(\bbC)$ définies ci-dessus. Notons $\delta_i$ la série discrète
  correspondant à   $\phi_i$. Le multi-ensemble de séries discrètes $\{ \delta_i \}_{i=1,\ldots,r}$
  paramètre  une classe d'équivalence de représentations irréductibles $\Pi_\phi$ de 
  $G_N$ d'après le théorème \ref{Lgclass}.
 Ceci  définit  la bijection voulue : 
 \begin{equation}
 \Phi(G_N)\longrightarrow \Pi(G_N), \qquad  \phi \mapsto \Pi_\phi.
 \end{equation}

\subsection{Paramètres et paquets d'Arthur pour $G_N$}\label{PaqArtGN}

Un paramètre d'Arthur pour $G_N$ est  un morphisme continu :
$ \psi: \; W_\bbR \times \SL_2(\bbC)\longrightarrow {}^LG_N=\widehat G_N \times W_\bbR $
vérifiant les propriétés énoncées dans la définition \ref{ArtPar}.
Comme pour les paramètres de Langlands, le fait que le $L$-groupe de $G_N$  soit un produit direct 
nous autorise à  considérer $\psi$ comme un morphisme  de $W_\bbR\times \SL_2(\bbC)$ dans 
$\widehat G_N=\GL_N(\bbC)$, c'est-à-dire comme une représentation de dimension $N$ de $W_\bbR\times \SL_2(\bbC)$.
Là encore, comme pour les paramètres de Langlands, une telle représentation est complètement réductible. Elle 
s'écrit donc comme une somme directe 
\begin{equation}
\label{opsii} \psi=\oplus_{i=1,\ldots,r} \psi_i, \qquad \psi_i :\; W_\bbR\times \SL_2(\bbC) \rightarrow \GL_{N_i}(\bbC)
  \end{equation}
avec $\psi_i$ irréductible et $\sum_{i=1}^r N_i=N$. Les  représentations irréductibles
de  $W_\bbR\times \SL_2(\bbC) $ sont des produits tensoriels de représentations irréductibles de $W_\bbR$
avec des représentations irréductibles de  $\SL_2(\bbC)$. Les représentations irréductibles
de  $W_\bbR$ ont été décrites ci-dessus. Celles qui apparaissent ici ont en plus la propriété d'être à image
bornée,  ce sont donc les représentations $W(s,\epsilon)$, $s\in i\bbR$, $\epsilon\in \{0,1\}$,  de dimension 1, et les 
représentations  $V(s_1,s_2)$, $s_1+s_2\in i\bbR$, $s_2-s_1 \in \bbN^\times$,  de dimension $2$.
 Les classes de représentations irréductibles de
 $\SL_2(\bbC)$ sont déterminées par leur dimension, et l'on note $R_d$ un choix de représentation irréductible 
 de dimension $d$ de $\SL_2(\bbC)$ (ou sa classe d'équivalence). Les représentations irréductibles de 
  $W_\bbR\times \SL_2(\bbC)$ qui nous intéressent sont donc à équivalence près 
  \begin{equation}
  W(s,\epsilon)\otimes R_n, \; V(s_1,s_2)\otimes R_n, \qquad s_2-s_1 \in \bbN^\times, \; \epsilon\in \{0,1\}, \; s, s_1+s_2\in i\bbR. 
  \end{equation}

Comme les paquets de Langlands, les paquets d'Arthur pour $G_N$ sont des singletons, et l'on a 
donc $\Pi_\psi=\Pi_{\phi_\psi}$, ceci désignant à la fois le paquet et l'unique représentation qu'il contient. 
 Nous allons maintenant décrire  la représentation $\Pi_\psi$ associée à un paramètre d'Arthur 
 $\psi$, en commençant par les $\psi$ irréductibles.
 
 \medskip
  
 Notons   $\mathbf{Triv}$   et $\sgn$ respectivement les caractères
triviaux et $\sgn$ de $\bbR^\times =\GL_1(\bbR)$, (bien que les m\^emes notations désignent aussi des caractères 
de $W_\bbR$, il n'y a pas de risque de confusion).  Si $\varepsilon \in \{ \Triv, \sgn\}$, et si  $n\in \bbN^\times$,   
notons    $\varepsilon_n$  le caractère de $G_n$ obtenu en composant  $\varepsilon$ et le déterminant
  $\det : \, G_n \rightarrow \bbR^\times $.
On a alors 
\begin{align}
&\text{Si } \psi=W(s,0)\otimes R_n, \quad \Pi_\psi=\Triv_n \; \nu^s=\nu^s\\  
&\text{Si } \psi=W(s,1)\otimes R_n, \quad \Pi_\psi=\sgn_n \; \nu^s  
\end{align} 
Remarquons que ces représentations sont des caractères unitaires de $G_n$, car on a pris $s\in i\bbR$.

\begin{defi} [Représentation de Speh]\label{defiSpeh} 
 Soit $\delta$ une série discrète unitaire  de $G_2$.  Considérons le module standard 
\begin{equation}\label{StdSpeh}
I(\delta,n)=\delta \nu^{-\frac{n-1}{2}} \times \delta\nu^{-\frac{n-3}{2}} 
\times \cdots \times \delta\nu^{\frac{n-1}{2}},
\end{equation}
et notons  $\Speh(\delta,n)$ son unique sous-module irréductible.   
Si $\delta=\delta(s_1,s_2)$,  $s_2-s_1 \in \bbN^\times$ et $s_1+s_2\in i\bbR$, on peut réécrire ceci comme 
   \begin{equation}\label{StdSpeh2}
 I\left(\delta\left(s_1,s_2\right),n\right)
 = \times_{i=1}^n \delta\left(s_1-\frac{n+1}{2}+i,s_2 -\frac{ n+1}{2}+i\right) 
\end{equation}
   
 On a alors,  
 \begin{equation}\label{PsiSpeh}
\text{si } \psi=V(s_1,s_2)\otimes R_n, \quad \Pi_\psi=\Speh\left( \delta\left(s_1,s_2\right),n \right).
  \end{equation} 
Les représentations $\Speh(\delta,n)$  sont unitaires.
\end{defi}

Nous venons donc de déterminer $\Pi_{\psi}$ lorsque $\psi$ est irréductible. Pour le cas général,
nous laissons au lecteur le soin de montrer que l'on a le résultat suivant. 
\begin{prop}
 Si $\psi=\oplus_{i=1,\ldots,r} \psi_i$
est une décomposition en irréductibles, alors
\begin{equation}\label{refreq}
 \Pi_\psi= \times_{i} \; \Pi_{\psi_i}. 
\end{equation}
\end{prop}

Par définition, $\Pi_\psi=\Pi_{\phi_\psi}$ et l'exercice consiste à  montrer que le paramètre de 
Langlands de la représentation irréductible $\times_{i} \; \Pi_{\psi_i}$ est bien  $\phi_\psi$.

\begin{rmq} 
Un résultat de Vogan \cite{VogGL} (voir aussi  \cite{Tadic} et  \cite{Baruch}) affirme que cette représentation 
est unitaire et irréductible, en particulier, 
elle ne dépend pas de l'ordre dans lequel on prend le produit.
\end{rmq}

\section{Paramètres de Beilinson-Bernstein pour $U(p,q)$}\label{BBUpq}

Dans cette section, $\mathbf{G}$ est le groupe unitaire $\U(p,q)$. Posons $p+q=N$.
Fixons une involution de Cartan $\tau$ de $G$, et soit $K=G^\tau$ le sous-groupe 
compact maximal de $G$ correspondant. On a $K\simeq U(p)\times U(q)$, et 
$K_\bbC\simeq \GL_p(\bbC)\times \GL_q(\bbC)$. Soit $\scrB$ la variété des drapeaux de $\frg$. 
Donnons un paramétrage combinatoire de $K_\bbC\backslash \scrB$. Pour cela, notons $\frI_N$ l'ensemble des 
involutions du groupe symétrique $\frS_N$ et introduisons
l'ensemble $\frI_N^{p,q,\pm}$ dont les éléments sont des couples $(\eta,f_\eta)$, où 
$\eta \in \frI_N$  et 
$f_\eta $ est une application de l'ensemble des points fixes de $\eta$ à valeurs dans $\{\pm 1\}$
qui vérifie, si l'on note $m$ le nombre de  $2$-cycles dans la décomposition en cycle de $\eta$ : 
\[ m+ \vert  f_\eta^{-1}(\{1\})  \vert =p, \qquad  m+ \vert  f_\eta^{-1}(\{-1\})  \vert =q. \]

Adoptons une  notation symbolique assez commode  pour les éléments de $\frI_N^{p,q,\pm}$ en écrivant
par exemple :
\[  \bar \eta= (+-1+23--312) \]
pour désigner l'élément dont les points fixes $(i)$ de l'involution sous-jacente $\eta$  sont repérés
un signe $\pm$ en position $i$, et bien entendu $f_\eta((i))=\pm$. Les 2-cycles $(ij)$ sont donnés
par  les positions $i$ et  $j$ où apparaissent les mêmes nombres, ici $(3,10)$, $(5,11)$ et $(6,9)$.
L'élément donné ici est dans $\frI_{11}^{5,6,\pm}$.
Par exemple, pour $(p,q)=(2,1)$, on a 
\[ \frI_3^{2,1,\pm} = \{ (++-), (+-+), (-++), (11+), (1+1), (+11)  \}.  \]

Le théorème 2.2.8 de \cite{Yam} donne une paramétrisation explicite de $K_\bbC\backslash \scrB$
par $\frI_N^{p,q,\pm}$. Notons 
\[  \bar \eta \in \frI_N^{p,q,\pm} \mapsto \caQ_{\bar \eta}\in K_\bbC\backslash \scrB \]
cette bijection.  L'ordre naturel sur $K_\bbC\backslash \scrB$ induit par transport de structure 
un ordre sur  $\frI_N^{p,q,\pm} $ qui est décrit en partie dans la section 2.4 de \cite{Yam}.
Seules certaines arêtes du diagramme de Hasse sont données dans \cite{Yam}, pour avoir 
la description complète de l'ordre, il faut rajouter les arêtes obtenues par la \og condition d'échange \fg \, 
 (3) de \cite{McGTr}.
On peut aussi définir facilement une fonction longueur $\ell_\frI$ sur $\frI_N^{p,q,\pm} $
qui va coïncider avec la dimension des orbites via la bijection.  Pour cela, posons
 pour tout $\bar \eta =( \eta,f_\eta) \in \frI_N^{p,q,\pm}$, 
 \[  \ell_\frI (\bar \eta)= \frac{1}{2}(p(p-1)+q(q-1))+ \sum_{\substack{2\text{-cycles } (ij) \text{ de } \eta,\\  i<j}}  
 \Big( (j-i)-\vert \{  2\text{-cycles } (kl) \text{ de } \eta, \;   \;   k<i<l<j \} \vert \Big).   \]

 Nous pouvons maintenant décrire plus précisément l'ordre sur  $\frI_N^{p,q,\pm} $ : il est gradué
par la fonction $\ell_\frI$, et il suffit de donner la liste des couples  d'éléments  
$(\bar \eta, \bar \eta')  \in \frI_N^{p,q,\pm}\times \frI_N^{p,q,\pm} $
tels que  $\ell_\frI(\bar \eta)=\ell_\frI(\bar \eta')+1$   et $\bar \eta \geq \bar \eta'$, c'est-à-dire les ar\^etes
du diagramme de Hasse. On a une ar\^ete lorsque  l'écriture symbolique de $\bar \eta'$ est obtenue à partir 
de celle pour $\bar \eta$ soit en remplaçant 2 symboles consécutifs de la forme $(aa)$, $a\in \bbN^\times$ 
(correspondant à un  
 un 2-cycle $(i,i+1)$) par $(+-)$ ou $(-+)$, soit en permutant deux symboles
consécutifs qui ne sont pas tout deux des signes (par exemple  $(1+1)$ et $(+11)$). 
({\sl cf.} \cite{Yam}), et il faut ensuite ajouter toutes les ar\^etes obtenues par la condition d'échange
(3) de \cite{McGTr}.

\bigskip 

Fixons maintenant une représentation $F$ de dimension finie de $G=\U(p,q)$. Notons $\Pi_F(G)$
l'ensemble des classes d'équivalence de représentations irréductibles  de $G$ ayant même caractère infinitésimal
que $F$. Les sous-groupes de Cartan de $G=\U(p,q)$ étant connexes,
la paramétrisation de Beilinson-Bernstein  { \cite{BeBe}, voir aussi \cite{Mil},\cite{VIC3}  Cor. 2.2 et \cite{Vgreen} Chapter 6)
   nous donne une bijection
\[\frI_N^{p,q,\pm} \simeq K_\bbC\backslash \scrB \simeq \Pi_F(G).   \]

\begin{rmq}\label{LpacUN}
Dans cette description, la partition de $\Pi_F(G)$ en $L$-paquets est particulièrement   simple : 
deux paramètres $\bar \eta, \bar \eta' \in \frI_N^{p,q,\pm}$
correspondent à des représentations dans le même $L$-paquet si et seulement si 
les involutions sous-jacentes $\eta$ et $\eta'$ dans $\frI_N$ sont égales. Ceci nous donne une 
paramétrisation de $\Phi(G)$ par $\frI_N$. On voit  dans ce cas que la longueur ne dépend 
pas des éléments du paquet. Nous verrons plus loin que la longueur 
$\ell_\frI$ est reliée à la longueur de Vogan, notée dans cet article $\ell_V$ (\cite{VIC3} et \cite{VIC4}). La propriété mentionnée ci-dessus 
pour les paquets de Langlands   est en fait une   propriété générale de la longueur de Vogan.
\end{rmq}

\section{L'espace tordu $\widetilde G_N$, le groupe non connexe $G_N^+$ et leurs représentations}\label{torduGLN}

\subsection{L'espace tordu $\widetilde G_N$ et ses représentations\label{espacetordu}}  
  
  On  note $\tau : \, g\mapsto {}^tg^{-1}$   l'involution de Cartan de $G_N$. 
 Soit $J_N\in G_N$ la matrice antidiagonale  {\tiny $ \begin{pmatrix}
 .&.&   &&& 1  \\   
.&. &&&-1&  \\
.&&& 1 && \\
&&.&&&\\
&.&&&&&\\
(-1)^{N-1}&.&.&&.&.\end{pmatrix}$ }
  
  On note $\theta_N :\, \mathbf{G}_N\rightarrow  \mathbf{G}_N$ l'automophisme défini par 
  $g\mapsto J_N\, ({}^tg^{-1})J_N^{-1}$.
  C'est un automorphisme intérieur à l'involution de Cartan $\tau$, mais qui a l'avantage de préserver 
  l'épinglage standard.

Nous pouvons maintenant définir le produit semi-direct $\mathbf{G}_N^+=\mathbf{G}_N \rtimes \langle \theta_N\rangle$.
 C'est un groupe algébrique réductif non connexe. L'automorphisme $\theta_N$ étant d'ordre 2,
 ce groupe compte deux composantes connexes, et l'on note $\widetilde{\mathbf{G}}_N=\mathbf{G}_N \rtimes  \theta_N$ celle qui
  ne contient pas l'élément neutre.  On note  $\widetilde G_N$ l'ensemble des points réels de $\widetilde{\mathbf{G}}_N$.
On obtient ainsi un espace tordu au sens de Labesse \cite{LaWa}. Pour le lecteur intéressé, remarquons que la construction ci-dessus est 
particulièrement simple puisque l'espace tordu en question  est l'ensemble des points réels
d'une composante connexe  d'un groupe algébrique, comme c'est le cas d'ailleurs 
de tous les espaces tordus pour lesquels on sait faire des choses intéressantes (par exemple, la stabilisation de  la formule des traces).

\medskip 

On appelle représentation de l'espace tordu $\widetilde{G}_N$ (ou pour faire bref représentation tordue) 
 tout triplet  $(\tilde \pi, \pi, V)$ où $V$ est un $\bbC$-espace vectoriel muni d'une représentation $\pi$ de $G_N$ et où $\tilde \pi$
est une  application $\tilde \pi: \, \widetilde G_N \rightarrow \GL(V)$
 tels que 
\begin{equation}\label{repestor}
  \tilde \pi (g_1\tilde x g_2)=\pi(g_1) \tilde \pi(\tilde x) \pi(g_2), \qquad (g_1,g_2\in G_N, \; 
\tilde x  \in \widetilde G_N).  \end{equation}

\medskip 

Si $(\pi,V)$ est une représentation de $G_N$, on note $\pi^{\theta_N}$ la représentation définie par 
\begin{equation} \pi^{\theta_N}(g)=\pi(\theta_N(g)), \qquad (g\in G_N). \end{equation}

\medskip 
Si  $(\pi,V)$ est une représentation de $G_N$ telle que $\pi^{\theta_N}$ est équivalente à $\pi$, on dit que 
$\pi$ est $\theta_N$-stable. Si $(\pi,V)$ est une représentation unitaire irréductible, ou bien
un module d'Harish-Chandra irréductible,  ceci est équivalent au fait que $(\pi,V)$ soit équivalente à 
sa contragrédiente. On dit alors aussi que $\pi$ est auto-duale.

\medskip 
 Si $(\tilde \pi, \pi, V)$ est une représentation de l'espace tordu 
$\widetilde G_N$, il résulte de (\ref{repestor})
que $ \tilde \pi(\theta_N)$ entrelace $\pi$ et $\pi^{\theta_N}$. Si l'on est dans une situation où le lemme
de Schur s'applique, par exemple si $\pi$ est une représentation unitaire irréductible de $G_N$, 
alors $\tilde \pi(\theta_N)$ est en fait déterminé par $\pi$ à un scalaire non nul près.
Réciproquement, si $A$ est un opérateur d'entrelacement inversible entre $\pi$ et $\pi^{\theta_N}$, on peut,
 définir $\tilde \pi : \,  \, \widetilde G_N \rightarrow \GL(V)$ par 
$\tilde \pi (g\theta_N)=\pi(g) A$, $g\in G_N$, et ceci fait de $(\tilde \pi, \pi, V)$ une représentation 
de l'espace tordu $\widetilde G_N$.

\medskip 

La définition  de représentation d'espace tordu s'adapte facilement aux modules de Harish-Chandra, nous 
laissons au lecteur le soin de faire ce travail. Soit $(\tilde \pi, \pi, V)$  une représentation de l'espace tordu 
$\widetilde G_N$, telle que le caractère de $\pi$ soit bien défini comme distribution
sur $G_N$ par la théorie d'Harish-Chandra (par exemple, $\pi$ est de longueur finie et admissible, à valeurs
dans un espace topologique localement convexe raisonnable). On peut alors définir de façon analogue le caractère tordu
de $\tilde \pi$. Tout d'abord, remarquons que $\widetilde G_N$ est muni d'une mesure invariante par les deux actions
de $G_N$ : on transporte une  mesure de Haar sur $G_N$ vers $\widetilde G_N$ par le difféomorphisme 
$g \mapsto g\theta_N=g \rtimes \theta_N$. Ensuite, pour toute fonction test $f \in \widetilde \caH_N
=\scrC_c^\infty(\widetilde G_N)$, l'opérateur 
\[  \tilde \pi(f)=\int_{\widetilde G_N}  f(\tilde y) \; \tilde \pi(\tilde y) \; dy    \]
est un opérateur à traces, et 
\begin{equation}\label{DefTr}  f \in   \widetilde \caH_N \mapsto   \Tr(   \tilde \pi(f))    \end{equation}
est une distribution sur $\widetilde G_N$, invariante par l'action adjointe de $G_N$. 

\medskip

 \subsection{Normalisation des extensions à $G_N^ +$}\label{NormExt}
 
Soit $(\tilde \pi, \pi, V)$  une représentation de l'espace tordu 
$\widetilde G_N$, telle que $\tilde \pi(\theta_N)^2=\Id_V$. Alors on peut définir une représentation $\pi^+$  du groupe
$G_N^+$ dont la restriction à $G_N$ est $\pi$ en posant $\pi^+(g\theta_N)=\pi(g)\tilde\pi(\theta_N)$.
Par exemple, si $(\pi,V)$ est une représentation à laquelle s'applique le lemme de Schur (une représentation 
irréductible unitaire, ou un module d'Harish-Chandra irréductible), et que $\pi^{\theta_N}$ est équivalente
à $\pi$, alors ce lemme entraîne que l'on peut choisir l'opérateur d'entrelacement $A$
réalisant cette équivalence de sorte que $A^2=\Id$, et ainsi permettre de définir deux extensions 
(les seules possibles) à $G_N^+$, disons $\pi^+$ et $\pi^-$, en posant respectivement 
$\pi^+(\theta_N)=A$ et $\pi^+(\theta_N)=-A$. Réciproquement, une représentation du groupe $G_N^+$
définit naturellement une représentation de l'espace tordu $\widetilde G_N$. Ainsi, l'on voit que 
l'on dispose de trois notions très proches : 
\begin{align}
 &\text{Les représentations $(\pi,V)$ de $G_N$ telles que $\pi^{\theta_N}$ est équivalente à $\pi$.}\label{piV}  \\
& \text{ Les représentations $(\tilde \pi, \pi,V)$ de l'espace tordu $\widetilde G_N$.}\label{tilde} \\
& \text{ Les représentations $(\pi^+,V)$ du groupe $G_N^+$}\label{plus}.
\end{align}
 Il sera commode d'énoncer certains résultats  en terme de représentations du groupe $G_N^+$. 
 Par exemple, si $(\pi^+,V)$ est une telle représentation avec $A=\pi^+(\theta_N) \in \GL(V)$,
 et si $(\tilde \pi, \pi,V)$ est la  représentation de  l'espace tordu associé, on pose 
\begin{equation} \label{Trplus} \Tr(\pi^+(f))=  \Tr(   \tilde \pi(f))  \qquad (f\in \widetilde \caH_N).
\end{equation}
Ceci définit une distribution sur $\widetilde G_N$, invariante pour l'action par conjugaison de $G_N$.
On peut exprimer  $\tilde \pi(f)$  en fonction de $\pi:=\pi^+_{\vert G_N}$ et $A$ par : 
\begin{equation} \tilde \pi(f)= \int_{\widetilde G_N}  f(\tilde y) \; \tilde \pi(\tilde y) \; dy =  
\int_{G_N}  f(g\theta) \;  \pi(g) \, A  \; dg .\end{equation}

Remarquons que si la représentation  $(\pi^+,V)$ de $G_N^+$ est donnée par une représentation
$(\pi,V)$ de $G_N$  et un opérateur d'entrelacement inversible $A$ entre $\pi$ et $\pi^{\theta_N}$
vérifiant $A^2=\Id_V$, et que $(\pi^-,V)$ est l'extension de $\pi$ à $G^+_N$ obtenue en remplaçant $A$
par $-A$, alors 
\begin{equation}\label{Trpm}
 \Tr(\pi^+(f))=- \Tr(\pi^-(f)), \qquad (f\in \widetilde \caH_N).
\end{equation}

Soit $(\pi,V)$ une représentation irréductible de $G_N$ à laquelle s'applique le lemme de Schur et supposons que 
 $ \pi^{\theta_N}$ est équivalente à $\pi$. 
 Nous allons expliquer comment  imposer un choix entre les deux extensions $\pi^\pm$ de $\pi$ à $G_N^+$, grâce aux fonctionnelles de Whittaker.

   On fixe une {\sl donnée de Whittaker} $(B_d,\chi)$ de $G_N$ de la manière suivante. 
On considère  sur $\bbR$  le caractère additif $\psi:\, x \mapsto \exp 2i\pi x$ 
 et l'on définit sur le radical unipotent $N_d$
du sous-groupe de Borel $B_d$ des matrices triangulaires supérieures dans $G_N$ le caractère
\[ \chi: (n_{ij}) \mapsto \psi(n_{12}+n_{23}+\cdots ++n_{N-1,N}). \]

Soit $(\pi,V)$ une représentation de $G_N$ dans un espace vectoriel topologique localement convexe raisonnable
(Hilbert, Banach, Fréchet, limite inductive de Fréchet, ...). Soit $V_\infty$ le sous-espace des vecteurs
$\scrC^\infty$ de $V$ et $V_\infty^*$ son dual topologique. Une fonctionnelle de Whittaker sur $(\pi,V)$
est un élément $\Omega \in V_\infty^*$ telle que 
\begin{equation}
\Omega(\pi(n) v)=\chi(n) \, \Omega(v), \qquad (v\in V_\infty), (n\in N_d).
\end{equation}

\begin{rmq}\label{chitheta}
Le sous-groupe unipotent $N_d$ est stable par $\theta_N$ et $\chi(\theta_N(n))=\chi(n)$ pour 
tout $n \in N_d$.
\end{rmq}

Supposons tout d'abord $(\pi,V)$ tempérée (en particulier unitaire) et irréductible. Alors, un résultat de Shalika
\cite{Shal} affirme que  $(\pi,V)$ admet une fonctionnelle de Whittaker non nulle, unique à un scalaire non nul près.

Soit  $P=MN$  un sous-groupe parabolique standard  de $G_N$ de facteur de Levi $M$ et de radical unipotent $N$.
La donnée de Whittaker $(B_d,\chi)$ sur $G_N$ définit par restriction une donnée de Whittaker 
$(B_M,\chi_M)=(B_d\cap M, \chi_{\vert M\cap N_d})$ pour $M$. 
Il est clair 
que le résultat de Shalika s'étend immédiatement aux représentations tempérées de $M$, et même à celles qui sont
essentiellement tempérées, c'est-à-dire produit tensoriel d'une représentation tempérée avec une représentation 
de dimension un. Supposons plus généralement  que $(\sigma, V)$ soit une  représentation de longueur finie de $M$
admettant une fonctionnelle  de Whittaker non nulle  $\Omega_M$ (pour la donnée de Whittaker $(B_M,\chi_M)$).
Rappelons comment définir une fonctionnelle de Whittaker non nulle sur $\Ind_P^{G_N} (\sigma)$.
Pour ceci, nous suivons la discussion dans \cite{Art13}, p. 111, auquel nous renvoyons
pour les notations (usuelles) non introduites ici (voir aussi \cite{Art89b}). Soit $\lambda \in \fra_{M,\bbC}^*$ et 
$\sigma_\lambda$ la torsion de $\sigma$ par le caractère  de $M$ défini par $\lambda$.
On réalise l'espace de la représentation induite  $\Ind_P^{G_N} (\sigma)$ comme un espace de Hilbert $\scrH_P(\sigma)$
de fonctions sur le sous-groupe compact maximal  $K=G_N^\tau$, cet espace restant le même lorsque $\sigma$ est remplacée
par $\sigma_\lambda$ (c'est l'action de $G_N$  qui change).
Si $h\in \scrH_P(\sigma)$, il faut poser 
\begin{equation}\label{fixe}
 h_{\sigma,\lambda}(x)= \sigma(M_P(x))\cdot h(K_P(x))e^{(\lambda+\rho_P)(H_P(x))}, \qquad (x\in G_N)
  \end{equation}
pour obtenir un vecteur dans l'espace usuel de $\Ind_P^{G_N} (\sigma_\lambda)$.
Soient $\bar w_l$ et $\bar w_l^M$ les éléments les plus longs dans les groupes de Weyl  de 
$G_N$ et $M$ respectivement,  posons $\bar w_M=\bar w_l \bar w_l^M$, $M'=\bar w_M\cdot M$ et soit $P'=M'N'$ 
le sous-groupe parabolique standard de $G_N$ de facteur de Levi standard $M'$. Fixons un représentant
$w_M$ de $\bar w_M$ dans $G_N$ (nous le ferons explicitement plus tard en (\ref{wM})).
Pour tout $h\in \scrH_P(\sigma)$,  l'intégrale de Whittaker 
\begin{equation}\label{IntWh} \mathrm{Wh}(h,\sigma_\lambda)=\int_{N'} \Omega_M(h_{\sigma,\lambda}(w_M^{-1}n')) \; 
\chi(n')^{-1}\; dn', \qquad 
(h\in \scrH_P(\sigma)) \end{equation}
converge absolument lorsque $\Re e(\lambda)$ se trouve dans un certain cône et admet un prolongement analytique
comme fonction de $\lambda \in \fra_{M,\bbC}^*$ à $\fra_{M,\bbC}^*$  tout entier 
(\cite{Shahi} Prop. 3.2, voir aussi \cite{Sha10}, lemma 3.6.8 et corollary 3.6.11) .
La fonctionnelle
\begin{equation}\label{WhMG}
\Omega : h \in \scrH_P(\sigma) \mapsto  \mathrm{Wh}(h,\sigma)
\end{equation}
est alors une fonctionnelle de Whittaker non nulle pour  $\Ind_P^{G_N} (\sigma)$.
Remarquons que ces définitions dépendent du choix de $w_M$.

Lorsque la fonctionnelle de Whittaker $\Omega_M$ est unique à une constante multiplicative près, 
il est expliqué  dans \cite{Sha10} que  l'induite $\Ind_P^{G_N}(\sigma)$ admet à un scalaire près une unique   
fonctionnelle de Whittaker, qui est donc celle notée   $\Omega$ définie ci-dessus. 
En particulier, les représentations standard
 admettent une fonctionnelle de Whittaker non nulle, unique à un scalaire non nul près.

\bigskip

Soit $(\pi,V)$  une représentation irréductible $\theta_N$-stable de $G_N$.
On fixe une extension  $\pi^+$ de $\pi$ à $G_N^+$ de la manière suivante.
Soit $(\rho,W)$ la représentation standard d'unique sous-représentation irréductible
$(\pi,V)$. Comme nous l'avons remarqué plus haut, la représentation standard $(\rho,W)$ admet 
 une unique droite de fonctionnelles de Whittaker. Fixons l'une d'elle, $\Omega$, non nulle, mais remarquons
 que ce choix  n'a aucune incidence sur ce qui suit. Il découle des théorèmes de classification que 
  $\rho$ est aussi $\theta_N$-stable. Soit $\caA$ un opérateur d'entrelacement non nul  entre $\rho$ et $\rho^{\theta_N}$. 
  En tenant compte de  la remarque \ref{chitheta}, on calcule pour tout $v \in V_\infty$ et tout $n \in N_d$,
\[   \Omega \circ \caA (\rho(n)v)=\Omega (\rho^{\theta_N}(n)\caA(v))= \Omega (\rho(\theta_N(n))\caA(v))=\chi(n) \; 
 \Omega \circ  \caA (v) ,   \]
ce qui montre que $ \Omega \circ \caA $ est aussi une fonctionnelle de Whittaker non nulle pour $(\rho,V)$.
Il existe donc $c \in \bbC^\times$ tel que $\Omega \circ \caA=c\,  \Omega$. Posons 
$\rho^+(\theta_N)=c^{-1}\, \caA$. On a alors $\Omega=\Omega\circ \rho^+(\theta_N)$ et  $\rho^+(\theta_N)^2=\Id_W$, 
ce qui fait que l'on définit ainsi une extension $\rho^+$ de $\rho$ à $G_N^+$.

 L'opérateur d'entrelacement  $\caA$   fixé comme ci-dessus  préserve l'unique  sous-représentation irréductible $\pi$, 
 et définit ainsi une extension   de $\pi$ à $G_N^+$.
 Celle-ci ne dépend pas des choix faits pour la construire (autres que celui de la donnée de Whittaker). 
 On note encore $\pi^+=(\tilde \pi,\pi,V)$ la représentation de
l'espace tordu $\widetilde G_N$ donnée par $\pi^+$.

\begin{defi} \label{ExtCanO} Soit $(\pi,V)$  une représentation irréductible $\theta_N$-stable de $G_N$.
L'extension  $\pi^+$ de $\pi$ à $G_N^+$ ou à $\widetilde G_N$ construite ci-dessus sera appelée {\sl extension canonique}, ou encore
{\sl extension déterminée par la donnée de Whittaker}. \end{defi}

\begin{rmq} Ceci n'est pas la façon dont procède Arthur dans \cite{Art13}, p. 63-64, pour définir les extensions
canoniques, mais les deux définitions
sont en fait équivalentes. Pour vérifier cette assertion, remarquons que celle-ci est tautologique pour les 
représentations tempérées.  Grâce aux  propriétés d'analyticité de la normalisation d'Arthur (\cite{Art13}, p. 64) et à 
 celles des fonctionnelles de Whittaker obtenues  par prolongement analytique de l'intégrale (\ref{IntWh}), 
 il suffit de vérifier l'assertion dans le domaine de convergence, domaine dans lequel un calcul direct permet de conclure.   
\end{rmq}

\begin{defi}\label{TrT}
Soit $\pi$ une représentation auto-duale de $G_N$ et soit $\pi^+$ son extension à $G^+_N$ ou à $\widetilde G_N$ déterminée 
par la donnée de Whittaker. On pose alors 
\[ \Tr_{\theta_N}(\pi(f))=\Tr(\pi^+(f)) \qquad (f \in \widetilde \caH_N). \]
\end{defi}

\section{Trace  tordue des représentations de Speh autoduales}  \label{caractor}

\subsection{Résolution de la représentation triviale de $\GL_n(\bbC)$}\label{cartorC}
On adopte dans cette section les notations sur les représentations des groupes complexes de \cite{BV}, appliquées au groupe $\GL_n(\bbC)$.
On note  $\rho$ la demi-somme des racines de $\GL_n(\bbC)$ et  pour tout $w\in \frS_n$ (identifié au groupe de Weyl de
$\GL_n(\bbC)$), on considère la
 série principale $X(\rho,-w\rho)$.  Si $w=1$, cette série principale est tempérée et si $w=w_0$, l'élément de plus grande longueur  
de $\frS_N$, cette série principale a pour unique sous-module irréductible la représentation triviale $\Triv_n$. 
 En général, pour tout $w\in \frS_n$, $X(\rho,-w\rho)$ a pour unique sous-module irréductible une représentation que l'on note 
 $\bar X(\rho,-w\rho)$.

On note $\overline{\theta}$ l'automorphisme de $\GL_n(\bbC)$ conjugué par $J_n$ (cf.section \ref{espacetordu}) 
de l'automorphisme $g\mapsto \, ^t\overline{g}^{-1}$. On remarque que l'image par $\overline{\theta}$ de
 $X(\rho,-w\rho)$ est la série principale $X(\rho, -w^{-1}\rho)$.
On munit le groupe de Coxeter $\frS_n$ de son  ordre de Bruhat  $\leq_B$ et de sa longueur usuelle $\ell_\frS$.
Lorsqu'il n'y a pas de risque de confusion avec une autre longueur, ce qui est le cas dans ce paragraphe, nous la notons simplement $\ell$.
 Pour tout $w,w'\in \frS_n$ tel que 
$\ell(w)=\ell(w')+1$ et $w'\leq_B w$, il existe un morphisme surjectif de $X(\rho,-w\rho)$ 
sur $X(\rho,-w'\rho)$, unique à un scalaire près. L'existence d'un tel morphisme surjectif est un exercice facile dans $\GL(2,\bbC)$. 
L'unicité à un scalaire près  vient de ce que $X(\rho,-w'\rho)$ a un unique sous-module irréductible qui intervient avec multiplicité 
un dans $X(\rho,-w\rho)$, c'est une propriété des polynômes de Kazhdan-Lusztig. Notons $f_{w,w'}$ un tel morphisme et posons 
$f_{w,w'}=0$ si la condition sur les longueurs est vérifiée mais pas celle sur l'ordre.

Pour $i= 0,\ldots ,\ell(w_0)$, on note $X_i:=\oplus_{w| \ell(w)=i}X(\rho,-w\rho)$. Se donner pour tout $i\in 1,\ldots ,\ell(w_0)$ un 
morphisme $\phi_i$  de $X_i$ dans $X_{i-1}$, revient bien évidemment à se donner une famille de morphismes 
$(f_{w,w'})_{\{(w,w')\vert \ell(w)=i, \, \ell(w')=i-1\}}$
 et réciproquement.

Le résultat suivant  est un cas particulier de ceux de  \cite{Joh1}.
\begin{prop}\label{ResCTriv}
(i) Il existe des choix de $f_{w,w'}$ comme ci-dessus de sorte que la suite d'applications:
\begin{equation}\label{EqResTriv}0\rightarrow \Triv_n \rightarrow X(\rho,-w_0\rho)\rightarrow \cdots\rightarrow X_i
\stackrel{\phi_i}{\longrightarrow} X_{i-1}\rightarrow \cdots
 \rightarrow X(\rho,-\rho)\rightarrow 0
\end{equation}
soit un complexe exact.

(ii) On peut fixer la famille des $f_{w,w'}$ de sorte que pour tout $w \in \frS_n$, il existe un morphisme $\bar \theta$-invariant
$A^{\bar{\theta}}_w$ de $X(\rho,-w\rho)$ dans $X(\rho, -w^{-1}\rho)$ vérifiant, quels que soient
 $w,w'$  avec $\ell(w')=\ell(w)-1$,   $A^{\bar{\theta}}_w \circ A^{\bar \theta}_{w^{-1}}=1$ et 
$A^{\bar{\theta}}_{w'}\circ f_{w^{-1},w'}\circ A^{\bar{\theta}}_{w}= f_{w,w^{'-1}}$.
\end{prop}
\dem Le (i) est dû à Johnson dans un cadre beaucoup plus général. Donnons
 le  point clé de la démonstration qui permet d'obtenir aussi (ii). Johnson construit inductivement
 (ici en faisant décroître l'indice) les applications $f_{w,w'}$.
  Supposons que les applications sont construites jusqu'à $\phi_{i+1}:X_{i+1}\rightarrow X_i$, alors 
  Johnson montre que pour tout $w'\in \frS_n$ de longueur $i-1$, 
 $\bar X(\rho,-w'\rho)$ intervient avec multiplicité un dans le conoyau $Y_i=X_i/\phi_{i+1}(X_{i+1})$
 et que $\dim \Hom_{\GL_n(\bbC)}(Y_i,X(\rho, -w' \rho))=1$. Soit $f_{i,w'}$ un 
  morphisme non nul dans $\Hom_{\GL_n(\bbC)}(Y_i,X(\rho, -w' \rho))$. En sommant sur les $w'$, 
  on obtient un morphisme de $X_i$ dans $X_{i-1}$ qui convient. On a donc une totale liberté sur les choix de $f_{i,w'}$. 
  Avec cela on obtient aussi facilement  (ii). \qed

\subsection{Résolution de Johnson  de  $\Speh\left(\delta \left(-\frac{p}{2},\frac{p}{2}\right), n\right)$}
\label{SecSpehRes}
Soit $p,n\in \bbN^ \times$  tels que $p > n-1$
  et soit $\Speh\left(\delta \left(-\frac{p}{2}, \frac{p}{2}\right), n\right)$ 
la représentation de Speh de $G_N=G_{2n}=\GL_{2n}(\bbR)$ de la définition \ref{defiSpeh}. 
Il se trouve que cette représentation est un $A_\frqqq(\lambda)$ de Vogan-Zuckerman (voir par exemple \cite{KV}, p. 330
pour une définition de ces représentations).
 En effet, 
il existe une sous-algèbre parabolique complexe  $\tau$-stable $\frqqq=\frl\oplus \fru$  de $\caM_{2n}(\bbC)$
telle que le sous-groupe de Levi $L$ associé soit la copie de $\GL_n(\bbC)$ naturellement contenue dans  $\GL_{2n}(\bbR)$. 
Le caractère  $\lambda$ de $L=\GL_n(\bbC)$ que l'on induit cohomologiquement est 
$\mu_{p,n} : \, g\mapsto (\det g)^{p-n}$.
Le caractère infinitésimal de $\Speh\left(\delta \left(-\frac{p}{2},\frac{p}{2}\right), n\right)$ est donné par 
{\footnotesize  \begin{equation}\label{CISp} 
 \left(\frac{-p-(n-1)}{2}, \frac{-p-(n-3)}{2}, \ldots ,  \frac{-p+(n-1)}{2}, \frac{p-(n-1)}{2},\frac{p-(n-3)}{2},
 \ldots ,  \frac{p+(n-1)}{2} \right)
 \end{equation}}
et la condition $p>n-1$ assure que ce caractère infinitésimal est bien entier et régulier.

\medskip

On obtient une résolution de Johnson de $\Speh\left(\delta \left(-\frac{p}{2}, \frac{p}{2}\right), n\right)$ de la manière suivante : 
on tensorise par $\mu_{p,n}$ la résolution de la représentation triviale de $\GL_n(\bbC)$ de la section \ref{cartorC}.
 Avec les notations de cette section,
 les modules standards  apparaissant dans cette résolution sont donc les $X(\rho, -w\rho)\otimes \mu_{p,n}$, $w\in \frS_n$.
Comme dans \cite{Joh1}, on induit ensuite cohomologiquement de $\GL_n(\bbC)$ à $\GL_{2n}(\bbR)$ 
ces modules, pour obtenir des modules standard que nous notons $X(ww_0)=X(ww_0,n,p)$ 
qui fournissent les termes d'une  résolution  de $\Speh\left(\delta \left(-\frac{p}{2},\frac{p}{2}\right), n\right)$. 
Remarquons que nous avons inversé l'ordre de l'indexation sur le groupe $\frS_n$ en multipliant par l'élément le plus long $w_0$.
En utilisant la comparaison entre la classification des représentations irréductible \og à la Langlands\fg\, et la classification \og à la 
Vogan-Zuckerman\fg  ({\sl cf.} \cite{Vgreen} Chapter 6 ou \cite{KV},  Chapter 11) on calcule que le  
 module standard  $X(w)=X(w,n,p)$ de $\GL_{2n}(\bbR)$  indexé par $w$ dans cette résolution peut-être décrit d'une autre manière, comme
induite parabolique ordinaire d'une représentation essentiellement de carré intégrable modulo le centre, à savoir 
\begin{equation}\label{Xss}
X(w)= X(w,n,p)=\times_{i=1}^{n}   \delta\left(  \frac{-p-(n-1)}{2}+(i-1),\frac{p-(n-1)}{2}+(w(i)-1)    \right) .
\end{equation}
 La résolution de Johnson s'écrit alors 
\begin{equation} \label{ResSpeh0} 0\rightarrow  \Speh\left(\delta \left(-\frac{p}{2}, \frac{p}{2}\right), n\right)
\rightarrow X_0 \rightarrow X_1 \rightarrow \cdots
 \rightarrow X_{\frac{n(n-1)}{2}} \rightarrow 0\end{equation}
où $X_0=X(1)$ et plus généralement $X_i$ est la somme des $X(s)$, $s\in \frS_n$ avec $\ell(s)=\ell_\frS(s)=i$.
Remarquons que le module standard $X(1)$ avait été noté  $I\left(\delta \left(-\frac{p}{2}, \frac{p}{2}\right), n\right)$ en 
(\ref{StdSpeh2}).
La résolution de Johnson induit une identité de représentations virtuelles dans le groupe de Grothendieck :
\begin{equation}\label{ResSpeH}
[\Speh\left(\delta \left(-\frac{p}{2}, \frac{p}{2}\right), n\right)]=\sum_{s\in \frS_n} (-1)^{\ell_{\frS}(s)} [X(s)].
\end{equation}
ce qui donne en prenant le caractère une  identité de distributions 
 \begin{equation}\label{ThetaSpeH}
  \Theta_{\Speh\left(\delta \left(-\frac{p}{2}, \frac{p}{2}\right), n\right)} = \sum_{s\in \frS_N}
  (-1)^{\ell_{\frS}(s)} \Theta_{X(s)}. \end{equation} 
  Dans cet article, ce ne sont pas les identités (\ref{ResSpeH}) et  (\ref{ThetaSpeH}), qui nous sont utiles, mais leurs analogues
  qui donnent la trace tordue de  $\Speh\left(\delta \left(-\frac{p}{2}, \frac{p}{2}\right), n\right)$ et que nous établissons dans la section suivante.
  
  \subsection{Trace tordue  de  $\Speh\left(\delta \left(-\frac{p}{2}, \frac{p}{2}\right), n\right)$}
\label{carspeh}
Nous voudrions maintenant obtenir une formule analogue à (\ref{ThetaSpeH}), mais pour la trace  tordue 
$\Tr_{\theta_N}\left( \Speh\left(\delta \left(-\frac{p}{2}, \frac{p}{2}\right) \right),n \right) $. 
Pour cela, introduisons la définition suivante :
\begin{defi}\label{thetalentgth} Notons $\frI_n$ l'ensemble des involutions de $\frS_n$.
Soit $w\in\frI_n$. On définit la $\theta$-longueur de $w$ en posant:
$$
\ell_\theta(w)=\frac{\ell(w)}{2}+\frac{|i\in [1,n]; w(i)>i| }{2}.
$$
\end{defi}
C'est le nombre d'orbites pour l'action de $\theta$ dans l'ensemble $\mathcal{E}_w$ des couples  formés de deux caractères 
de la forme $ z\mapsto z^i/\overline{z}^{w(i)}; i\in [1,n] $, ces couples étant  indexés par les  $(i,j)\in [1,n]^2$ 
vérifiant  $i<j;$ et $w(i)>w(j)$.  Evidemment $|\mathcal{E}_w|=\ell(w)$.
Cette $\theta$-longueur a aussi été introduite et étudiée dans  \cite{Incitti}, via la remarque suivante.

\begin{rmq}
Nous avons déjà rencontré cette fonction longueur dans la section \ref{BBUpq}. En effet, 
\begin{equation}\label{complinv}  \ell_\frI(\eta)= \ell_\theta(\eta) +\frac{1}{2}(p(p-1)+q(q-1)),   \end{equation}
où ici, la longueur du terme de gauche est celle définie sur $\frI_n$ dans la section \ref{BBUpq}. 
 La vérification de cette égalité est un exercice combinatoire un peu fastidieux, pour lequel la lecture de \cite{Incitti} nous a été 
d'une aide  précieuse. On peut s'en servir pour démontrer le résultat 
de la section \ref{ParSpehE}.   Mais comme nous donnerons aussi un autre argument qui permet
d'éviter d'utiliser cette égalité (et même qui la démontre), nous le laissons  au lecteur.
\end{rmq}

On rappelle que chaque module standard $X(w)$, pour $w$ une involution, est muni d'une action de $\theta=\theta_N$.
 Cette action est celle qui laisse invariante toute fonctionnelle de Whittaker.
\begin{thm}\label{THMTTSpeh}
 La trace tordue  ({\sl cf.} Définition \ref{TrT}) 
de $\Speh\left(\delta \left(-\frac{p}{2},\frac{p}{2}\right),n \right)$ est donnée par 
\begin{equation}\label{FormTTSpeh}
\Tr_{\theta_N}\left( \Speh\left(\delta \left(-\frac{p}{2},\frac{p}{2}\right),n \right) \right)
=\sum_{w\in \frI_n} (-1)^{\ell_\theta(w)} \Tr_{\theta_N}\left( X(w)\right).\end{equation}
\end{thm}
\dem  On note ici simplement $\ell=\ell_\frS$ la longueur usuelle dans le groupe symétrique. Le complexe  (\ref{ResSpeh0}) qui résout 
$\Speh\left(\delta \left(-\frac{p}{2}, \frac{p}{2}\right),n\right)$, est muni d'une action de $\theta=\theta_N$. Il est très vraisemblable que cette action 
 est l'image fonctorielle de celle de  $\bar\theta$ de la proposition \ref{ResCTriv},
  mais plutôt que de vérifier cela, il est plus simple de reprendre l'argument donné dans cette proposition  pour ajuster les choix des $f_{w,w'}$.
   Les morphismes dans la résolution
(\ref{ResSpeh0}) sont obtenus comme somme de morphismes  $f_{w,w'}: \, X(w)\rightarrow X(w')$ où ici $\ell_{w'}=\ell(w)+1$.
Rappelons que nous avons renversé l'indexation entre (\ref{EqResTriv}) et (\ref{ResSpeh0}) en multipliant par $w_0$, mais 
à ce changement d'indice près, les $f_{w,w'}$ sont les images par le foncteur d'induction cohomologique des $f_{w,w'}$ de (\ref{EqResTriv}).
   On obtient  pour tout $w\in \frS_n$ un morphisme  $A_\theta(w)$ de $X(w)$ dans $X(w^{-1})$.
   Si $w$ n'est pas une involution, $X(w)\oplus X(w^{-1})$ est $\theta$-stable, mais la contribution à la trace tordue est nulle. 
  Pour tout  $w\in \frI_n$, on note alors $\Tr_{A_\theta}(X(w))$  la trace tordue de $X(w)$  pour l'action de $\theta$ donnée par 
  l'opérateur   $A_\theta(w)$. On impose comme cela est loisible, quitte à  multiplier tous les choix de $A_\theta$ 
  par $-1$, que l'action de $A_\theta(1)$,  pour $w=1$ 
   est l'action de $\theta$ normalisée par le modèle de Whittaker. 
   On obtient alors immédiatement que l'injection de la représentation  de $\Speh$ dans $X(1)$ est $\theta$-invariante et que l'on a:
  \begin{equation}\label{TTSp1}
\Tr_{\theta_{N}}\left( \Speh\left(\delta \left(-\frac{p}{2}, \frac{p}{2}\right),n \right) \right)
=\sum_{w\in \frI_n} (-1)^{\ell(w)} \Tr_{A_\theta}\left( X(w)\right).\end{equation}
 Pour démontrer le théorème, il suffit de voir que $A_\theta$ sur $X(w)$ (pour $w$ une involution) diffère de l'action de $\theta$ normalisée 
 par les modèles de   Whittaker  par le signe $(-1)^{\ell_\theta(w)-\ell(w)}$.  C'est ce que nous allons voir
 maintenant, en commençant par un lemme technique.
\begin{lemme} \label{lemtech} 
Soit $w\in\frI_n$ une involution qui n'est pas l'identité. Alors l'un des deux cas suivant a lieu (non exclusivement);

(i) soit il existe une involution, $w'\in \frI_n$ telle que $\ell(w')=\ell(w)-1$ et $w'<w$ et $\ell_\theta(w')=\ell_\theta(w)-1$;

(ii) soit il existe  $w'\in \frI_n$ et $s\in \frS_n$ tels que $\ell(w')= \ell(s)-1=\ell(w)-2$, $w<s<w'$  et $\ell_\theta(w')=\ell_\theta(w)-1$.
\end{lemme}
\dem On peut très bien démontrer ce lemme en exhibant $w',s$ pour tout $w$. Mais un référé nous en a donné une  preuve  élégante.

 Soit $\mathcal{X}$ un sous-ensemble de $\{1,\ldots,n-1\}$ tel que si $i\in \mathcal{X}$ alors $i+1\notin {\mathcal X}$. 
  A ${\mathcal X}$ on associe une involution qui est le produit des involutions échangeant $i$ et $i+1$ pour tout $i\in {\mathcal X}$. 
  On note  $\sigma_{\mathcal X}$ cette involution. Soit $w$ une involution, alors il existe un sous-ensemble ${\mathcal X}$ tel que 
  $w$ soit conjugué dans $\frS_n$ de $\sigma_{\mathcal X}$. Et on vérifie, par exemple par récurrence sur $n$ que l'on peut écrire 
  $w=\tau^{-1}\sigma_{\mathcal X}\tau$ où $\tau\in \frS_n$ avec les longueurs qui s'ajoutent, c'est-à-dire
\begin{equation}\label{longueursplus}
\ell(w)= 2\ell(\tau)+\ell(\sigma_{\mathcal X}).
\end{equation}
Les $\theta$-longueur sont très faciles à calculer: soit $y$ une involution. On note $\mathrm{Inv}(y)$ l'ensemble des couples $(i,j)$
 tels que $1\leq i<j\leq n$ tels que $y(i)>y(j)$ et $\mathrm{Exc}(y)$ l'ensemble des couples $(i,y(i))$ tel que $i<y(i)$. 
 Alors $\ell_{\theta}(y)=\frac {|\mathrm{Inv}(y)|+|\mathrm{Exc}(y)|}{2}.$ En revenant à (\ref{longueursplus}), 
 on trouve que $\ell_\theta(w)= \ell(\tau)+|{\mathcal X}|$.

Supposons que $\tau$ soit l'identité, on fixe alors $i\in {\mathcal X}$ car par hypothèse ${\mathcal X}$ n'est pas vide.
 On note ${\mathcal X}'$ l'ensemble ${\mathcal X}$ privé de $i$ et $w'=\sigma_{{\mathcal X}'}$ répond à (i).

Supposons maintenant que $\tau$ ne soit pas l'identité. On fixe une symétrie élémentaire, $\tau_0$ tel que $\tau':=\tau_0\tau<\tau$ 
et on pose $w'= (\tau')^{-1}\sigma_{\mathcal X}\tau'$ et $s=\tau^{-1}\sigma_{\mathcal X}\tau'$. Ces deux éléments de $\frS_n$ répondent à (ii). \qed

En multipliant par  l'élément le plus long $w_0$, on obtient 
\begin{cor} \label{corlemtech} 
Soit $w\in\frI_n$ une involution qui n'est pas l'élément le plus long $w_0$. Alors l'un des deux cas suivant a lieu (non exclusivement);

(i) soit il existe une involution, $w'\in \frI_n$ telle que $\ell(w')=\ell(w)+1$ et $w'>w$ et $\ell_\theta(w')=\ell_\theta(w)+1$;

(ii) soit il existe  $w'\in \frI_n$ et $s\in \frS_n$ tels que $\ell(w')= \ell(s)+1=\ell(w)+2$, $w'>s>w$  et $\ell_\theta(w')=\ell_\theta(w)+1$.
\end{cor}

\begin{rmq}\label{rmqtech} Plaçons nous dans le cas de (ii) en fixant $w'$. Alors $s$ et $s^{-1}$ sont les seuls éléments $s'$
 vérifiant $w<s'<w'$ et en particulier l'exactitude du complexe (\ref{ResSpeh0}) entraîne que 
 $f_{s,w'}\circ f_{w,s}+f_{s^{-1},w'}\circ f_{w,s^{-1}}=0$.
\end{rmq}
\medskip

Pour tout $w\in \frI_n$, notons $\theta(w)$ l'action de $\theta$ sur $X(w)$ normalisé par les modèles de Whittaker. 
Nous avons vu ci-dessus qu'il suffisait de montrer le résultat suivant pour terminer la démonstration du théorème \ref{THMTTSpeh}.
\begin{lemme}  Pour tout  $w\in \frI_n$, on a $A_\theta(w)=(-1)^{\ell_\theta(w)-\ell(w)}\theta(w)$.
\end{lemme}

\dem
Soit $\epsilon_0$ le signe tel que $A_\theta(w_0)=\epsilon_0 \theta(w_0)$, où $w_0$ est l'élément de plus grande longueur  de $\frS_n$.
 On montre que $A_\theta(w)=\epsilon_0(-1)^{\ell_\theta(w_0)-\ell(w_0)}(-1)^{\ell_\theta(w)-\ell(w)}\theta(w)$ 
  par récurrence  descendante sur la longueur de $w$.
 Ceci est donc   vrai par construction pour l'élément $w=w_0$. On fixe $w\in \frI_n$ qui n'est pas $w_0$
  et on lui applique le corollaire \ref{corlemtech}. Plaçons nous dans le cas (i) de ce lemme. On fixe $w'$ comme dans (i)
   et on remarque que $f_{w,w'}$ permet de remonter une fonctionnelle de Whittaker  $\Omega_{w'}$ sur $X(w')$ en une fonctionnelle de
    Whittaker $\Omega_{w}=\Omega_{w'}\circ f_{w,w'}$ sur $X(w')$.
   Remarquons que les morphismes $f_{w',w}$ sont surjectifs. C'est une propriété héritée
   des morphismes de la résolution (\ref{EqResTriv})  pour lesquels cela avait été noté.
      Notons $\sigma(w)$ le signe que l'on cherche à calculer; par hypothèse de récurrence, on a 
    $\sigma(w'):=\epsilon_0(-1)^{\ell_\theta(w_0)-\ell(w_0)}(-1)^{\ell_\theta(w')-\ell(w')}$. On a donc
$$\Omega_{w}\circ A_\theta(w)=\Omega_{w'}\circ f_{w,w'}\circ A_\theta(w)=\Omega_{w'}\circ A_\theta(w')\circ f_{w,w'}
=\sigma(w')\Omega_{w'}\circ f_{w,w'}=\sigma(w')\Omega_{w}.
$$
On en déduit  $\sigma(w)=\sigma(w')$. Or  $\ell_\theta(w)=\ell_\theta(w')-1$ et $\ell(w)=\ell(w')-1$ et on obtient l'assertion cherchée.

Plaçons nous dans le cas (ii) en fixant $w'$ et $s$ et en reprenant  les  notations $\sigma(w),\sigma(w')$. 
On a $\ell(w')=\ell(w)+2$ et $\ell_\theta(w')=\ell_\theta(w)+1$. On cherche donc à montrer que $\sigma(w)=-\sigma(w')$. 
On fixe encore une fonctionnelle de Whittaker $\Omega_{w'}$ sur  $X(w')$  et on note $\Omega_{w}$ la fonctionnelle de Whittaker telle que 
$\Omega_{w}=\Omega_{w'}\circ f_{s,w'}\circ f_{w,s}$.
On a alors :
$$\Omega_{w}\circ A_\theta(w)=\Omega_{w'}\circ f_{s,w'}\circ f_{w,s}\circ A_\theta(w)
=\Omega_{w'}\circ A_\theta(w')\circ f_{s^{-1},w'}\circ f_{w,s^{-1}}.
$$
Ici on utilise la remarque de \ref{rmqtech} pour remplacer $f_{s^{-1},w'}\circ f_{w,s^{-1}}$ par $-f_{s,w'}\circ f_{w,s}$, et l'on obtient 
$$\Omega_{w}\circ A_\theta(w)
=-\Omega_{w'}\circ A_\theta(w')\circ f_{s,w'}\circ  f_{w,s}=-\sigma(w')\Omega_w.    $$

C'est bien le changement de signe cherché. On calcule $\epsilon_0$ sur l'identité  $1$ de $\frS_n$ :
 par construction $\sigma(1)=1$. Ainsi $\epsilon_0=(-1)^{\ell_\theta(1)-\ell(1)}(-1)^{\ell_\theta(w_0)-\ell(w_0)}=(-1)^{\ell_\theta(w_0)-\ell(w_0)}$
  et ceci termine la démonstration du lemme, et 
 donc du théorème. \qed

\begin{rmq}\label{caslim}
Nous aurons besoin plus loin de pousser les égalités (\ref{ThetaSpeH}) et (\ref{FormTTSpeh}) au cas où $p=n-1$. Le caractère infinitésimal
de $\Speh\left(\delta \left(-\frac{p}{2}, \frac{p}{2}\right), p+1\right)$ n'est plus régulier ($0$ appara\^it
avec une multiplicité $2$) dans le caractère infinitésimal. Par continuation cohérente ({\sl cf.} \cite{Vgreen}, Chapter 7), 
  (\ref{ThetaSpeH}) et (\ref{FormTTSpeh})  restent
 valides, avec la convention
\[ \delta(0,0)= \Triv\times \sgn.  \]
\end{rmq}

\section{Induction parabolique tordue}  \label{IndTor}

Le but de cette section est de définir un foncteur d'induction parabolique \og  tordu \fg \, 
d'un sous-espace de Levi tordu $\widetilde M$ vers l'espace tordu $\widetilde G_N$.
Nous commençons par fixer des notations pour tout le reste de cette section.

\subsection{Notations}\label{NotII}

Soient $a$ un entier au moins égal à $1$, $N_1=2a$  et considérons la représentation    
$\Speh(\delta\left(-\frac{p}{2},\frac{p}{2}\right),a)$  de  $\GL_{N_1}(\bbR)$.
Soit $N'$ un autre entier au moins égal à $1$, et soit $\lambda'=(\lambda'_1, \ldots, \lambda'_{N'})\in \bbC^{N'}$.
On fait dans la suite l'hypothèse suivante :
\begin{equation}\label{hyplemmegen}
\frac{p-(a-1)}{2} >\vert  \lambda'_i  \vert, \qquad  (1\leq  i \leq N').  
\end{equation}

 Soit $M$ le sous-groupe de Levi standard de $\GL_N(\bbR)$ isomorphe à  $\GL_{N_1}(\bbR)\times  \GL_{N'}(\bbR)$
et soit $P=MN$ le sous-groupe parabolique standard de facteur de Levi $M$.

En (\ref{Xss}), nous avons défini des modules standard $X(w)=X(w,a,p)$, $w\in \frS_a$, qui entrent dans la résolution de Johnson de 
$\Speh(\delta\left(-\frac{p}{2},\frac{p}{2}\right),a)$. Ces modules standard sont en position de Langlands et admettent un unique
sous-module irréductible que nous notons $\bar X(w)$.

\begin{defi}\label{srcE}
Soit  $\scrE$ l'ensemble des représentations irréductibles $\pi_M=\pi_1 \otimes \pi'$ de $M$ telles que 
\begin{itemize}
\item[(i)] $\pi_1$ est une représentation irréductible de $\GL_{N_1}(\bbR)$ parmi les $\bar X(w), w\in \frS_{a}$ (ce sont 
les représentations irréductibles de  même caractère infinitésimal que $\Speh(\delta\left(-\frac{p}{2},\frac{p}{2}\right),a)$ 
dont le paramètre dans le  $\caG$-ordre de Bruhat  ({\sl cf.}  \cite{VIC3}, \cite{VIC4})   est inférieur ou égal à celui de 
$\Speh(\delta\left(-\frac{p}{2},\frac{p}{2}\right),a)$).
\item[(ii)] $\pi'$ est une représentation irréductible de $\GL_{N'}(\bbR)$ de caractère infinitésimal $\lambda'$.
\end{itemize}
\end{defi}

\begin{rmq}  \label{heredite}
 Les propriétés du $\caG$-ordre de Bruhat impliquent immédiatement  que 
 l'ensemble $\scrE$ possède la propriété suivante : si $\displaystyle \pi_M=\pi_1 \otimes  \pi'$ 
 est dans  $\scrE$, et si $X_{\pi_M}$ est le module standard dont $\pi_M$ est l'unique sous-module irréductible, 
 alors tous les sous-quotients irréductibles de $X_{\pi_M}$ sont dans $\scrE$.
  \end{rmq}
  
\subsection{Un résultat d'irréductibilité d'induite}\label{new11}

Le résultat d'irréductibilité suivant est crucial pour la suite.

\begin{lemme} \label{irredIndpiM} 
Pour toute représentation $\pi_M \in \scrE$, $\Ind_P^{\GL_N(\bbR)} (\pi_M)$ est irréductible.
\end{lemme}

\dem
On réalise l'induite comme quotient d'un module standard dont les exposants sont dans la chambre de Weyl positive
 et comme sous-module d'un module standard dont les exposants sont dans la chambre de Weyl négative, 
 ces deux modules standard ayant même sous-quotient de Langlands. Cela force l'induite à être ce sous-quotient de 
 Langlands et a fortiori à être irréductible. Donnons une idée de la démonstration, les détails se trouvent dans \cite{Ara}, Lemme 3.1.2 ou 
 dans la version préliminaire de cet article \cite{AMR} Lemme 10.6.

La représentation $\pi_1=\bar X(w)$ est donc l'unique sous-représentation 
irréductible du module standard  $X(w)$. Notons  $I'$ le module standard  de $\GL_{N'}(\bbR)$ dont $\pi'$ est un sous-module irréductible 
et $I$ le module standard de  $\GL_{N}(\bbR)$ dont les exposants sont dans la chambre de Weyl négative et dont le caractère 
est celui de l'induite $X(w)\times I'$  (les représentations $X(w)\times I'$ et $I$ sont des induites à partir de la même 
représentation du sous-groupe de Levi $M$, mais pour $I$, on induit  d'un sous-groupe parabolique qui n'est pas $P$, mais 
celui qui met les exposants de $X(w)\times I'$ dans un ordre standard). 
 Par des opérations élémentaires dans $\GL_3(\bbR)$ et $\GL_4(\bbR)$, on fait des échanges de facteurs pour ramener 
 cette induite à $I$ et on vérifie à chaque fois que ces échanges sont des isomorphismes. On procède ensuite de même avec le parabolique
 qui met les exposants de $X(w)\times I'$ dans un ordre standard inverse.
Ainsi la démonstration repose {\sl in fine}  sur le  résultat suivant  
 dû  à B. Speh (\cite{Speh}).
  
\begin{thm}\label{SpehIrrStd}  
Pour $i=1,2$,  $\eta_i$ désigne soit $\Triv$ soit  $\sgn$, c'est-à-dire une représentation de $G_1=\GL_1(\bbR)$,
auquel cas on pose $n_i=1$,  soit une série discrète $\delta_i=\delta\left(-\frac{p_i}{2},\frac{p_i}{2}\right)$ 
de $G_2=\GL_2(\bbR)$, auquel cas on pose 
$n_i=2$.   
Alors la représentation  
\[ \eta_1\nu^{t_1} \times \eta_2 \nu^{t_2} ,  \quad t_i \in \bbC, \quad i=1,2,\]
 de $G_{n_1+n_2}=\GL_{n_1+n_2}(\bbR)$ est réductible  si et seulement si $t_1-t_2\in \bbR$ et si l'on est dans un des cas suivant : 

\medskip 

\begin{itemize}
\item[(i)]  $n_1=n_2=1$,  $\eta_1 =\eta_2$, $t_1-t_2  \in 2\bbZ+1$.   
\item[(ii)]  $n_1=n_2=1$,  $\eta_1 \neq \eta_2$, $ t_1-t_2  \in 2\bbZ^\times$.   
\item[(iii)]  $n_i=2, n_j=1$, $\{i,j\}=\{1,2\}$,   $t_1-t_2\in \bbZ$, $-\frac{p_i}{2}+\vert t_1-t_2 \vert \in \bbN^\times$.
\item[(iv)]  $n_1=n_2=2$, $t_1-t_2\in \bbZ$, $-\frac{\vert p_1-p_2\vert }{2}+\vert t_1-t_2 \vert \in \bbN^\times$.
\end{itemize}

\end{thm}

En particulier, si $p'$ et $p''$ sont des demi-entiers, 
 une induite de la forme $\delta(-p'',p')\times \tau$, où $\tau$ est une représentation irréductible de 
$\GL_b(\bbR)$ avec $b=1$ ou $2$, est irréductible si le caractère infinitésimal de $\tau$ est de la forme $(m)$ ou $(m',m)$ avec 
$m$ et $m'$ des demi-entiers vérifiant:
\begin{equation}\label{critirr}
p'\geq m \geq m'\geq -p'' .\end{equation}

C'est ce critère qui sert dans les échanges de facteurs mentionnés ci-dessus, et plus précisément, on l'applique avec 
$p'\geq p/2-(a-1)/2$ et $-p''\leq -p/2+(a-1)/2$, 
l'hypothèse (\ref{hyplemmegen}) garantissant que l'on a bien (\ref{critirr}). \qed

Notons la conséquence suivante du lemme \ref{irredIndpiM}.
  \begin{cor} \label{passtdstd}
Soient   $X(w)$, $w\in \frS_n$ une des représentations standard de $\GL_{N_1}(\bbR)$ définie en (\ref{Xss}) et $X'$ un module standard de
$\GL_{N'}$ de caractère infinitésimal $\lambda'$. Considérons 
\[ X(w)\times X' .\]   
Cette représentation n'est pas une représentation standard
 {\sl a priori} car  les facteurs ne sont  pas forcément dans un ordre standard.
Les arguments de la démonstration du   lemme  \ref{irredIndpiM} montrent que l'on peut permuter les facteurs dans l'induite  pour l'amener
dans un ordre standard. Elle admet donc un unique sous-module irréductible.  
   \end{cor}

\subsection{Morphismes entre représentations standards}

Nous allons  introduire les représentations   standard  de $G_{N'}$ de caractère infinitésimal 
$\lambda'$ et maximales
pour le $\caG$-ordre de Bruhat ({\sl cf.} \cite{VIC3}, \cite{VIC4}).    Ce sont des séries principales.
On suppose que $\lambda'=(\lambda'_1,\ldots,\lambda'_{N'})$ est donné dans un ordre standard, 
c'est-à-dire que $(\Re e(\lambda'_i))_i$ est une suite croissante.

\begin{defi}  
Soit   
$\underline \varepsilon'=(\varepsilon_1',\ldots ,\varepsilon_{N'}')\in \{\Triv,\sgn\}^{N'} $. 
Posons alors 
\begin{equation}\label{epsilonprimem} I(\underline \varepsilon',\lambda')=\varepsilon_1' \nu^{\lambda'_1} \times \varepsilon_2'\nu^{\lambda'_2} 
\times \cdots \times \varepsilon_{N'}'\nu^{\lambda'_{N'}}.\end{equation}
C'est une représentation  standard de $G_{N'}$. Posons aussi
\begin{equation}\label{epsilonprime}
I_M^{\underline{\varepsilon'}}=X(1)\otimes I(\underline \varepsilon' ,\lambda')
\end{equation}
 \end{defi}

\begin{rmq}
Soit $\pi$ un sous-quotient irréductible de $I_M^{\underline{\varepsilon'}}$ et soit $X_\pi$ le module standard
dont $\pi$ est l'unique sous-module irréductible. Il découle de la remarque \ref{heredite} que 
 tous les sous-quotients irréductibles de 
$X_\pi$ sont dans $\scrE$. 
\end{rmq}

\bigskip 

  \begin{lemme} \label{morsur}
   Soit $\pi_M= \pi_1\otimes \pi'  \in \scrE$. 
   Soit $X_M$ le module standard de $M$ dont $\pi_M$ est l'unique sous-module irréductible.
  C'est un produit tensoriel des modules standard $X_1$ de   $G_{N_1}$ et $X'$ de  $G_{N'}$ dont  $\pi_1$ et $\pi'$ sont 
   les uniques sous-module irréductibles.
   
   Alors il existe     un morphisme surjectif (de représentations de  $G_{N_1}$) 
  \[ h_1 :\; X(1) \longrightarrow  X_1, \]
 et il existe  $\underline \varepsilon'=(\varepsilon_1',\ldots ,\varepsilon_{N'}')\in \{\Triv,\sgn\}^{N'} $ 
 et un morphisme surjectif (de représentations de $G_{N'}$)
   \[ h' :\; I(\underline\varepsilon',\lambda') \longrightarrow  X'.\]
      \end{lemme}

\dem Le premier point découle immédiatement  par une récurrence sur la longueur du fait que si $w,w'$ sont dans $\frS_a$ avec 
$w<_Bw'$ et $\ell(w)=\ell (w')-1$, alors le morphisme $f_{w,w'}$ entrant dans la résolution (\ref{ResSpeh0}) est surjectif.
Pour le second,  on raisonne par récurrence sur $N'$. Il n'y a rien a montrer si $N'=1$, et si $N'=2$ ou $N'=3$
l'assertion du lemme est facile à vérifier car la structure des représentations standard
est totalement connue dans ces petits rangs. On suppose donc $N'\geq 4$ et le résultat établi
pour $N''<N'$. On écrit 
\begin{equation}\label{wX}  X' = \sigma \times  \rho \end{equation}
où $\sigma$ est une série discrète de $G_d$, $d=1$ ou $2$ et $\rho$ une représentation standard de $G_{N'-d}$. 
 Rappelons que la condition   pour que $X' $ soit une représentation standard admettant une unique 
sous-représentation irréductible est que ce soit une induite de séries discrètes écrites dans un ordre standard.
Or, avec les notation de la définition \ref{OrdStd},  $e(\varepsilon\nu^\lambda)=\Re e(\lambda)$ et 
$e(\delta(\alpha_1,\alpha_2))=\Re e\left(\frac{\alpha_1+\alpha_2}{2}\right)$. Si $d=1$ ceci implique que  
 la représentation $\sigma$  est de la forme $\varepsilon'_1 \nu^{\lambda'_1}$. Si $d=2$, alors nécessairement
 dans l'écriture de $X'$ en produit de séries discrètes dans l'ordre standard inverse,
il apparaît une série discrète de la forme $\delta(\lambda'_1,\lambda'_s)$ (avec en particulier $\Re e(\lambda'_s)>\Re e(\lambda'_1)$)
 et à gauche de celle-ci dans le produit  
  uniquement des séries discrètes $\delta(\lambda'_j, \lambda'_k )$ pour  $j<k  \in \{2,\ldots,s-1\}$ avec 
  $\Re e\left(\frac{\lambda'_j+\lambda'_k}{2}\right)\leq \Re e\left( \frac{\lambda'_1+\lambda'_s}{2}\right)$. 
  Or, d'après le théorème \ref{SpehIrrStd}, on a alors 
 \[ \delta(\lambda'_j, \lambda'_k )\times \delta(\lambda'_1,\lambda'_s)=\delta(\lambda'_1,\lambda'_s) \times 
 \delta(\lambda'_j, \lambda'_k )\]
car les deux membres sont des représentations irréductibles. On peut donc supposer que l'on a 
$\sigma=\delta(\lambda'_1,\lambda'_s)$ dans (\ref{wX}).

Dans le  cas $d=1$, on conclut en appliquant l'hypothèse de récurrence à $\rho$ : 
il existe des caractères $\varepsilon'_i \in \{\Triv,\sgn\}$, $i=2,\ldots,N'$, tels que l'on ait un morphisme surjectif
\[ \varepsilon'_2\nu^{\lambda_2} \times \cdots \times \varepsilon'_{N'} \nu^{\lambda_{N'}}  \twoheadrightarrow \rho \]
et l'on en déduit par fonctorialité de l'induction parabolique un morphisme surjectif
\[  \varepsilon'_1\nu^{\lambda_1} \times \cdots \times \varepsilon'_{N'} \nu^{\lambda_{N'}}  \twoheadrightarrow  
  \sigma \times  \rho = X'.  \]

 Si $\sigma= \delta(\lambda'_1, \lambda'_s )$ pour un $s\in \{2,\ldots,N'\}$, on applique encore 
 l'hypothèse de récurrence à $\rho$ : 
 il existe des caractères $\varepsilon_i \in \{\Triv,\sgn\}$, $i\in\{2,\ldots,N'\}\setminus \{s \}$,
  tels que l'on ait un morphisme surjectif
\[ \varepsilon'_2\nu^{\lambda'_2} \times \cdots  \times \varepsilon'_{s-1} \nu^{\lambda'_{s-1}} \times
 \varepsilon'_{s+1} \nu^{\lambda'_{s+1}} \times   \cdots \times\varepsilon'_{N'} \nu^{\lambda_{N'}}  \twoheadrightarrow \rho .\]

 Par fonctorialité de l'induction parabolique, on a  un morphisme surjectif
\[ \delta(\lambda'_1,\lambda'_s)\times \varepsilon'_2\nu^{\lambda'_2} \times \cdots  \times
\varepsilon'_{s-1} \nu^{\lambda'_{s-1}} \times \varepsilon'_{s+1} \nu^{\lambda'_{s+1}}
  \times  \cdots \times\varepsilon'_{N'} \nu^{\lambda'_{N'}}  \twoheadrightarrow \sigma\times \rho=X'. \]

Si $s=2$, comme il existe $\varepsilon'_1,\varepsilon'_2\in \{\Triv,\sgn\}$ et  un morphisme surjectif
\[ \varepsilon'_1\nu^{\lambda'_1}\times\varepsilon'_2\nu^{\lambda'_2}  \twoheadrightarrow \delta(\lambda'_1,\lambda'_2),\]
on en déduit immédiatement l'assertion. Supposons donc $s\geq 3$.
Si $\Re e(\lambda'_s)=\Re e(\lambda'_2)$, alors $\Re e(\lambda'_i)=\Re e(\lambda'_2)$ pour tout $i\in \{2,\ldots, s\}$ 
et l'on est ramené au cas précédent grâce au théorème \ref{SpehIrrStd}.
 Si $\Re e(\lambda'_2)<\Re e(\lambda'_s)$,  
il existe $\varepsilon'_1\in \{\Triv,\sgn\}$ et  un morphisme surjectif
\[ \varepsilon'_1\nu^{\lambda'_1}\times \delta(\lambda'_2,\lambda'_s)   \twoheadrightarrow \delta(\lambda'_1,\lambda'_s)\times
\varepsilon'_2\nu^{\lambda'_2} \]
(c'est une situation dans $\GL_3(\bbR)$).
On a donc un morphisme surjectif 
\[ \varepsilon'_1\nu^{\lambda'_1}\times \delta(\lambda'_2,\lambda'_s) \times \varepsilon'_3\nu^{\lambda'_3}
 \times \cdots  \times \varepsilon'_{s-1} \nu^{\lambda'_{s-1}} \times \varepsilon'_{s+1} \nu^{\lambda'_{s+1}} \times 
  \cdots \times\varepsilon'_{N'} \nu^{\lambda'_{N'}}   \twoheadrightarrow X'.\]
On applique encore une fois l'hypothèse de récurrence à 
$\delta(\lambda'_2,\lambda'_s) \times \varepsilon'_3\nu^{\lambda'_3} \times \cdots  \times \varepsilon'_{s-1} \nu^{\lambda'_{s-1}}
 \times \epsilon'_{s+1} \nu^{\lambda'_{s+1}} \times 
  \cdots \times\varepsilon'_{N'} \nu^{\lambda'_{N'}}   $
pour conclure.\qed

\subsection{Opérateurs d'entrelacement} \label{opeent}  Définissons l'automorphisme involutif  $\theta_M$ par la formule 
\begin{equation}\label{m}
m= \begin{pmatrix} g_1&  \\ 
                 &  g'        \end{pmatrix}
  \mapsto              \theta_M(m)= \begin{pmatrix} \theta_{N_1}(g_1)& \\ 
                             &  \theta_{N'}(g')                             \end{pmatrix}.
\end{equation} 
Posons $J_M=\begin{pmatrix} J_{N_1}&  \\ 
                 &  J_{N'}
                      \end{pmatrix}$. On a alors $\theta_M(m)= J_M ({}^tm^{-1}) J_M^{-1}$.
  Posons aussi 
  \begin{equation}\label{wM}
  w_M=J_NJ_M^{-1}.
\end{equation}
   Soient $P'$ le sous-groupe parabolique standard de facteur de Levi
  $w_M \cdot M\cdot w_M^{-1}$ et $N'$ son radical unipotent. On définit, pour tout $\underline s=(s_1,s')\in \bbC^2$, le 
  caractère $\chi_{\underline s}$ de $M$ par la formule 
  \[ \chi_{\underline s}(m)= \vert \det (g_1) \vert^{s_1}\times  \vert \det (g') \vert^{s'} \]
  lorsque $m \in M$ est comme en (\ref{m}).

Soit $\pi_M$ une représentation de $M$ (disons de longueur finie).  On définit l'opérateur d'entrelacement 
\begin{equation}
M(\underline s,w_M,\pi_M) : \; \Ind_P^{G_N} (\pi_M\otimes \chi_{\underline s}) \longrightarrow 
\Ind_{P'}^{G_N} (w_M \cdot ( \pi_M\otimes \chi_{\underline s}) )
\end{equation}  
 comme fonction méromorphe de $\underline s$ par prolongement analytique  d'une formule intégrale 
 (\cite{KnSt1},  \cite{KnSt2},\cite{Sch},  voir aussi  \cite{Art89b}, \S1).
  Plus précisément, on considère une autre réalisation de $\Ind_P^{G_N} (\pi_M\otimes \chi_{\underline s})$, où 
 l'espace de la représentation est un espace de fonctions $\scrH_P( \pi_M)$ sur le compact maximal $G_N^\tau$ de $G_N$,
  cet espace  ayant l'avantage de ne pas dépendre de $\underline s$ (c'est l'action de $G_N$ qui en dépend). 
 Si $f\in \scrH_P( \pi_M)$, soit $f_{\underline s}$ la fonction de $\Ind_P^{G_N} (\pi_M\otimes \chi_{\underline s})$
 qui lui correspond par l'isomorphisme (donné en  (\ref{fixe}) avec des notations légèrement différentes)
  entre $\scrH_P( \pi_M)$ et  $\Ind_P^{G_N} (\pi_M\otimes \chi_{\underline s})$.
  On pose alors 
  \begin{equation}\label{Ms}
  \big(M(\underline s,w_M,\pi_M)(f_{\underline s})\big)(g)=\int_{N'}f_{\underline s}(w_M^{-1}n'g)\; dn'.
  \end{equation}
Cette intégrale converge absolument dans un certain cône de $\bbC^N$ et est analytique en $\underline s$
dans ce cône. Elle admet un prolongement méromorphe à $\bbC^N$.

 Remarquons que $M(\underline s,w_M,\pi_M)$ dépend du choix de l'élément  $w_M$ fait ci-dessus et pas seulement
 de sa classe dans $N_G(M)/M$.

 Posons $N_{M,\emptyset}=M\cap N_d$, $N_{M',\emptyset}=M'\cap N_d$.  
La donnée de Whittaker $(B_d,\chi)$ pour $G_N$ définit par restriction des données de Whittaker 
$(B_d\cap M, \chi_{N_{M,\emptyset}})$ et $(B_d\cap M', \chi_{N_{M',\emptyset}})$.
On a $w_M^{-1}\cdot N_{M',\emptyset}= N_{M,\emptyset}$ et $w_M\cdot \chi_{N_{M,\emptyset}}=\chi_{N_{M',\emptyset}}$.

Nous avons introduit en (\ref{epsilonprime}) les représentations standard $I_M^{\underline{\varepsilon}'}$ de $M$,  
$\varepsilon'\in \{\Triv,\sgn\}^{N'}$. 
Pour un tel $\underline{\varepsilon}'$, fixons  $\Omega_M^{\underline{\varepsilon}'}$
 une fonctionnelle de Whittaker non nulle sur 
$I_M^{\underline{\varepsilon}'}$. On vérifie facilement que  $\Omega_M^{\underline{\varepsilon}'}$
 est aussi une fonctionnelle de Whittaker pour $I_M^{\underline{\varepsilon}'}\otimes \chi_{\underline s}$  
 ainsi que  pour $w_M\cdot (I_M^{\underline{\varepsilon}'}\otimes \chi_{\underline s})$.

La formule intégrale (\ref{IntWh}),  ou plutôt son prolongement analytique, définit  ({\sl cf.} \ref{WhMG})  
des fonctionnelles de Whittaker non nulles 
sur $\Ind_P^{G_N}(I_M^{\underline{\varepsilon}'}\otimes \chi_{\underline s})$ et 
$\Ind_{P'}^{G_N} (w_M \cdot (I_M^{\underline{\varepsilon}'}\otimes \chi_{\underline s}) 
)$ que nous notons respectivement $\Omega^{\underline{\varepsilon}'}(\underline s)$ et 
${\Omega'}^{\underline{\varepsilon}'}(\underline s)$.

La proposition suivante utilise de manière cruciale les résultats de F. Shahidi \cite{Shahi}.
\begin{prop}\label{OpEntrBiz}
$(i)$   Il existe une fonction méromorphe $\underline s\mapsto r(\underline s)$ telle que 
 pour tout $\varepsilon'\in \{\Triv,\sgn\}^{N_r}$, l'opérateur 
  d'entrelacement  : 
\begin{equation}\label{NormOp} 
N(\underline s,w_M,I_M^{\underline{\varepsilon}'})= r(\underline s)^{-1}  
 M(\underline s, w_M,I_M^{\underline{\varepsilon}'})
\end{equation}  
  soit holomorphe et inversible au point $s=0$ et vérifie  
\begin{equation}\label{EntWh}
\Omega ^{\underline{\varepsilon}'}(\underline s)=  {\Omega'}^{\underline{\varepsilon}'}(\underline s) \circ N(\underline s,w_M,I_M^{\underline{\varepsilon}'}).
\end{equation}

$(ii)$ Pour toute représentation de longueur finie $\pi_M$ de $M$ dont tous les sous-quotients irréductibles
sont dans $ \scrE$, l'opérateur  
  \begin{equation}\label{NormOp2} 
N(\underline s,w_M,\pi_M)= r(\underline s)^{-1}   M(\underline s,w_M,\pi_M)
\end{equation}  
est alors   holomorphe et inversible au point $s=0$. 
(Remarquons que le facteur  $r$ ne dépend pas de $\pi_M$.)

$(iii)$ Si de plus $\pi_M$ est   standard  et 
si  $\Omega_M$ est une fonctionnelle de Whittaker non nulle sur $\pi_M$, alors 
$\Ind_P^{G_N}(\pi_M\otimes \chi_{\underline s})$ 
et    $\Ind_{P'}^{G_N} (w_M \cdot (\pi_M\otimes \chi_{\underline s}))$ admettent  
 des  fonctionnelles de Whittaker non nulles   $\Omega(\underline s)$
et $\Omega'(\underline s)$ fixées comme ci-dessus 
(prolongement analytique (\ref{WhMG}) de la  formule intégrale   (\ref{IntWh})) et 
  \begin{equation}\label{EntWh2}
\Omega(\underline s)=\Omega'(\underline s)\circ N(\underline s,w_M,\pi_M).
\end{equation}
  \end{prop}
 
  \medskip

\dem Commençons par montrer $(i)$. 
Posons pour tout $\underline{\varepsilon}'=(\varepsilon'_1,\ldots , \varepsilon'_{N'})\in \{\Triv,\sgn\}^{N'}$ et pour tout $i=1,\ldots ,a$, 
\[\delta_{1,k}^{\underline{\varepsilon}'}=  \delta_{1,k}=  \delta\left(\frac{-p-(a-1)}{2}+(k-1), 
\frac{p-(a-1)}{2}+(k-1)\right) . \]

Posons aussi  pour tout  $\underline{\varepsilon}'=(\varepsilon'_1,\ldots , \varepsilon'_{N'})\in \{\Triv,\sgn\}^{N'}$ et  pour tout $j=1,\ldots ,N'$, 
\[ \delta_{2,j}= \nu^{\lambda'_j}, \quad \delta_{2,j}^{\underline{\varepsilon}'}=  \varepsilon'_{j} \nu^{\lambda'_j}. \]

Soit $\bar w_M$ l'élément de groupe de Weyl $W_G$ dont $w_M$ est un représentant et écrivons
\[ \bar w_M=\bar \tau_{j_1},\ldots,\bar \tau_{j_\ell}\]
comme produit de longueur minimale de générateurs de $W_G$ (ce groupe est isomorphe à $\frS_N$
et  $\bar\tau_k$ correspond à la  transposition $(k,k+1)$ par cet isomorphisme. 
Soient $\tau_{j_1},\ldots,\tau_{j_\ell}$ des représentants des
$\bar \tau_i$ dans $G$, de sorte que $ w_M=\tau_{j_1}\ldots \tau_{j_\ell}$.
 On a alors une factorisation de l'opérateur d'entrelacement
\begin{equation}
M(\underline s,w_M,I_M^{\underline{\varepsilon}'}) : \; X(1)\nu^{s_1}\times  I(\underline{\varepsilon}',\lambda')\nu^{s'} 
\longrightarrow 
 I(\underline{\varepsilon}',\lambda')\nu^{s'} \times X(1) \nu^{s_1}
\end{equation} 
en produit d'opérateurs d'entrelacement correspondant aux $\tau_{j_k}$.
 En effet,  nous avons :
\[X(1)^{\underline{\varepsilon}'}:=X(1)=\times_{i=1,\ldots, a}  \, \delta_{1,i}^{\underline{\varepsilon}'}, \qquad 
  I(\underline{\varepsilon}',\lambda')= \times_{j=1,\ldots, N'} \,  \delta_{2,j}^{\underline{\varepsilon}'}. \]
Considérons un produit des $\delta_{i,i'}^{\underline{\varepsilon}'}\nu^{s_i}$ dans un ordre quelconque. Lorsque  deux couples $(j,j')$,  
$(k,k')$ dans ce produit sont  tels que   $\delta_{j,j'}^{\underline{\varepsilon}'}\nu^{s_j}$ apparaît immédiatement  à gauche de   
$\delta_{k,k'}^{\underline{\varepsilon}'}\nu^{s_k} $, 
disons en position $m$ et $m+1$, notons encore
 $M_{j,j',k,k'}^{\underline{\varepsilon}'}(\tau_m,\underline s)$ l'opérateur d'entrelacement échangeant ces deux facteurs par 
 \begin{equation} \label{Opj} \delta_{j,j'}^{\underline{\varepsilon}'}\nu^{s_j} \times  \delta_{k,k'}^{\underline{\varepsilon}'}\nu^{s_k} 
  \longrightarrow   \delta_{k,k'}^{\underline{\varepsilon}'}\nu^{s_k}  \times \delta_{j,j'}^{\underline{\varepsilon}'}\nu^{s_j}  \end{equation}
  et agissant   trivialement sur les autres. Ainsi $M(\underline s,w_M,I_M^{\underline{\varepsilon}'})$
peut s'écrire comme composition d'opérateurs $M_{j,j',k,k'}^{\underline{\varepsilon}'}(\tau_m,\underline s)$ avec $j=1$ et $k=2$.  
Soit $r_{j,j',k,k'}^{\underline{\varepsilon}'}(\underline s)$ la fonction méromorphe introduite 
par F. Shahidi dans \cite{Shahi}, \S3.1 
 comme facteur   les opérateurs   $M_{j,j',k,k'}^{\underline{\varepsilon}'}(\tau_m,\underline s)$, et posons 
\[N_{j,j',k,k'}^{\underline{\varepsilon}'}(\tau_m,\underline s)=r_{j,j',k,k'}^{\underline{\varepsilon}'}(\underline s)^{-1} 
M_{j,j',k,k'}^{\underline{\varepsilon}'}(\tau_m,\underline s).\]
Comme $\delta_{j,j'} ^{\underline{\varepsilon}'}\times  \delta_{k,k'} ^{\underline{\varepsilon}'}\simeq   
 \delta_{k,k'} ^{\underline{\varepsilon}'} \times \delta_{j,j'} ^{\underline{\varepsilon}'}$ est irréductible, 
$N_{j,j',k,k'}^{\underline{\varepsilon}'}(\tau_m,\underline s)$ est holomorphe et bijectif en $\underline s=0$.
D'autre part, et c'est crucial, remarquons que les fonctions $r_{j,j',k,k'}^{\underline{\varepsilon}'}(\underline s)$ 
ne dépendent pas 
de $\underline{\varepsilon}'$. C'est clair si $j=k=1$, et si $j=1$ et $k=2$ on a :
 \begin{align*} & \delta_{1,j'}^{\underline{\varepsilon}'}\nu^{s_1} \times  
 \delta_{2,k'}^{\underline{\varepsilon}'}\nu^{s_2} = 
  \delta\left(\frac{-p-(a+1)}{2}+j', \frac{p-(a+1)}{2}+j'\right)  \nu^{s_1}\times  
   \varepsilon'_{k'} \nu^{\lambda_{k'}} \nu^{s_2} \\
  &=\varepsilon'_{k'}\left(  \delta\left(\frac{-p-(a+1)}{2}+j', \frac{p-(a+1)}{2}+j'\right)  \nu^{s_1}\times  
   \nu^{\lambda_{k'}} \nu^{s_2} \right)\\
  &= \varepsilon'_{k'} \left(  
  \delta_{1,j'}\nu^{s_1} \times  
 \delta_{2,k'}\nu^{s_2}\right).
 \end{align*}
Ceci est immédiat car pour tout $\varepsilon\in \{\Triv,\sgn \}$ et toutes représentations $\rho_1$, $\rho_2$ de 
$G_{a}$ et $G_b$ respectivement, $\varepsilon_{a+b}(\rho_1\times \rho_2)
=(\varepsilon_a\rho_1)\times (\varepsilon_b\rho_2)$ et si $\delta$ est une
 série discrète de $G_2$, $\varepsilon \delta=\delta$.
L'opérateur d'entrelacement (\ref{Opj}) dans ce cas devient un opérateur 
\[ \varepsilon'_{k'} \left(    \delta_{1,j'}\nu^{s_1} \times   \delta_{2,k'}\nu^{s_2}\right) 
\longrightarrow   \varepsilon'_{k'} \left(     \delta_{2,k'}\nu^{s_2}  \times  \delta_{1,j'}\nu^{s_1} \right)  \]
 donné par prolongement analytique d'une expression intégrale qui ne dépend pas de $ \varepsilon'_{k'}$.
Le facteur  $r_{j,j',k,k'}^{\underline{\varepsilon}'}(\underline s)$ n'en dépend donc pas non plus,
 et nous le notons simplement  $r_{j,j',k,k'}(\underline s)$.
 Le facteur  $r(\underline s)$ est alors défini comme le produit correspondant des $r_{j,j',k,k'}(\underline s)$
 et $N(\underline s,w_M,I_M^{\underline{\varepsilon}'})= r(\underline s)^{-1}  
 M(\underline s, w_M,I_M^{\underline{\varepsilon}'})$ est le produit des 
 $N_{j,j',k,k'}^{\underline{\varepsilon}'}(\tau_m,\underline s)$.  L'opérateur d'entrelacement $N(\underline s,w_M,I_M^{\underline{\varepsilon}'})$
 est holomorphe et inversible au point $s=0$ car les $N_{j,j',k,k'}^{\underline{\varepsilon}'}(\tau_m,\underline s)$ 
 le sont, et la propriété (\ref{EntWh}) est immédiate puisque la normalisation choisie est celle de Shahidi ({\sl cf. \cite{Shahi}, \S 3.1 }) et que les opérateurs  de Shahidi vérifient la propriété voulue de factorisation
 (\S 3.2 de \cite{Shahi}). Ceci termine la démonstration du $(i)$.

\bigskip

  Montrons maintenant que  l'on a $(ii)$ pour une représentation $\pi_M  \in \scrE$ irréductible.
  Comme l'élément $w_M$ ne joue plus de rôle important dans ce qui suit, notons
simplement $N(\underline s,I_M^{\underline{\varepsilon}'})$ pour $N(\underline s,w_M,I_M^{\underline{\varepsilon}'})$.   
  D'après le  lemme  \ref{morsur}, il existe  des sous-modules $V_1$ et $V_2$ de $I_M^{\underline{\varepsilon}'}$ tels que 
  \[ V_1/V_2 \simeq \pi_M . \] 
  On note $I_{M \vert V_i}^{\underline{\varepsilon}'}$ la restriction de $I_{M}^{\underline{\varepsilon}'}$ à $V_i$, $i=1,2$.

On a un diagramme commutatif 
\[  \xymatrix{  \Ind_P^{G_N}(I_{M \vert V_1}^{\underline{\varepsilon}'}
\otimes \chi_{\underline s})  \ar^{ N(\underline s,I_{M\vert V_1}^{\underline{\varepsilon}'})}[rrr] 
\ar@{^{(}->}[d]
&&& \Ind_{P'}^{G_N} (w_M \cdot (I_{M\vert V_1}^{\underline{\varepsilon}'}
\otimes \chi_{\underline s})) \ar@{^{(}->}[d] \\
  \Ind_P^{G_N}(I_M^{\underline{\varepsilon}'}\otimes \chi_{\underline s})  \ar^{ N(\underline s,I_M^{\underline{\varepsilon}'})}[rrr] & &&
   \Ind_{P'}^{G_N} (w_M \cdot (I_M^{\underline{\varepsilon}'}\otimes \chi_{\underline s} )) } \]
où   les flèches verticales sont  les inclusions évidentes.
 Comme $N(\underline s,I_M^{\underline{\varepsilon}'})$ est holomorphe et inversible en $\underline s=0$,
  il en est de même de  $N(\underline s,I_{M\vert V_1}^{\underline{\varepsilon}'})$.

  De même, on a un diagramme commutatif
  \[  \xymatrix{0\ar[d] &&& 0 \ar[d]\\
   \Ind_P^{G_N}(I_{M \vert V_2}^{\underline{\varepsilon}'}\otimes \chi_{\underline s})
     \ar^{ N(\underline s,I_{M\vert V_2}^{\underline{\varepsilon}'})}[rrr] \ar@{^{(}->}[d]
&&& \Ind_{P'}^{G_N} (w_M \cdot (I_{M\vert V_2}^{\underline{\varepsilon}'}
\otimes \chi_{\underline s})) \ar@{^{(}->}[d] \\
 \Ind_P^{G_N}(I_{M \vert V_1}^{\underline{\varepsilon}'}\otimes \chi_{\underline s})  
 \ar^{ N(\underline s,I_{M\vert V_1}^{\underline{\varepsilon}'})}[rrr] \ar[d]
&&& \Ind_{P'}^{G_N} (w_M \cdot (I_{M\vert V_1}^{\underline{\varepsilon}'}
\otimes \chi_{\underline s})) \ar[d] \\
  \Ind_P^{G_N}(\pi_M\otimes \chi_{\underline s})\ar[d]    \ar^{ N(\underline s,\pi_M)}[rrr] & &&
   \Ind_{P'}^{G_N} (w_M \cdot (\pi_M\otimes \chi_{\underline s} )) \ar[d]  \\
   0 &&& 0  } ,\]
  les colonnes étant exactes. Il en découle immédiatement 
  que  $N(\underline s,I_{M\vert V_2}^{\underline{\varepsilon}'})$ et 
   $N(\underline s,\pi_M)$ sont holomorphes en $\underline s=0$, 
  le premier étant aussi  inversible et le second surjectif. Si 
  $v \in \ker( N(0,\pi_M))$, relevons $v$ en un élément $w\in  \Ind_P^{G_N}(I_{M \vert V_1}^{\underline{\varepsilon}'})$. 
  L'élément 
   $N(0,I_{M\vert V_1}^{\underline{\varepsilon}'})(w)$  se projette sur l'élément nul de  
   $\Ind_{P'}^{G_N} (w_M \cdot \pi_M )$
   et est donc dans l'image de $\Ind_{P'}^{G_N} (w_M \cdot (I_{M\vert V_2}^{\underline{\varepsilon}'}))$.
    L'élément $w$ est donc 
   dans l'image de  $\Ind_P^{G_N}(I_{M \vert V_2}^{\underline{\varepsilon}'})$, 
   et ainsi $v=0$, ce qui montre que  $N(0,\pi_M)$ est injectif.
La première assertion du $(ii)$ est donc démontrée pour $\pi_M$ irréductible.

Passons au cas général. Supposons que   $N(\underline s,\pi_M)$ ait un pôle en $\underline s=0$. Choisissons
une fonction holomorphe $\underline s\mapsto l(\underline s)$ telle que $l(0)=0$ et 
  $\underline s\mapsto l(\underline s) N(\underline s,\pi_M)$ soit holomorphe non nul 
 en  $\underline s=0$. Soit $V$ un sous-module irréductible non nul de 
 \[  \Ind_P^{G_N}(\pi_M)/ \ker\left[   \biggl( l(s) N(\underline s,\pi_M)\biggr)_{s=0} \right] . \]
  Comme tout sous-quotient de $\pi_M$ est dans $\scrE$, il existe un $\pi'_M \in \scrE  $ 
  tel que $V\simeq \Ind_P^{G_N}(\pi'_M)$. L'opérateur  $l(\underline s) N(\underline s,\pi'_M)$ est holomorphe et non nul 
  en $\underline s=0$ par définition de $l$ et de $V$. D'autre part, nous avons montré ci-dessus que  
   $N(\underline s,\pi'_M)$ est holomorphe   en $\underline s=0$. On obtient alors 
     $  \biggl( l(s) N(\underline s,\pi'_M) \biggr)_{s=0}=0$ ce qui constitue une contradiction. 
     L'opérateur  $N(\underline s,\pi_M)$ est donc holomorphe  en $\underline s=0$.
  
  Supposons $\ker[ N(\underline 0,\pi_M) ]$ non nul, et soit $V$ un sous-module 
  irréductible de $\ker[ N(\underline 0,\pi_M) ]$.
  Il existe  $\pi'_M \in \scrE$
  tel que $V\simeq \Ind_P^{G_N}(\pi'_M)$. Nous avons montré ci-dessus que l'opérateur  $N(\underline s,\pi'_M)$
  est holomorphe  et inversible en $s=0$.  Ceci est manifestement incompatible avec le fait que 
   $N(\underline 0,\pi_M)$ soit nul sur $V$. Ceci montre  que $N(\underline 0,\pi_M)$
  est injectif. On montre de même que $N(\underline 0,\pi_M)$
  est surjectif en considérant un quotient irréductible $V$ de $\coker[ N(\underline 0,\pi_M) ]$ .
  Ceci termine de montrer $(ii)$
  
Il reste à montrer le dernier point $(iii)$ concernant les fonctionnelles de Whittaker. 
Soit $\pi_M=\pi_1 \otimes  \pi'  \in \scrE$.
 Soit $X_M$ le module standard de $M$ dont $\pi_M$ est l'unique sous-module irréductible.
  C'est un produit tensoriel des modules standard $X_1$ de $\GL_{N_1}(\bbR)$ 
  et  $X'$ de  $\GL_{N'}(\bbR)$ dont les $\pi_1$ et $\pi'$ sont 
   les uniques sous-module irréductibles. 
Rappelons les morphismes surjectifs $h_1$ et $h'$  du lemme \ref{morsur}.
Notons $h_M(\underline s): \;  \left( X(1)\otimes  I(\underline\varepsilon',\lambda') \right)
 \otimes \chi_{\underline s}\longrightarrow X_M\otimes \chi_{\underline s}$
 le morphisme surjectif (de représentations de $M$)
obtenu par produit tensoriel de $h_1$ et $h'$ et  
\[ h(\underline s)= \Ind_{P=MN}^{G_N} (h_M(\underline s)):\,    \Ind_{P=MN}^{G_N}\left(  \left( X(1)\otimes  I(\underline\varepsilon',\lambda')   \right)
\otimes \chi_{\underline s}\right)  \longrightarrow \Ind_{P=MN}^{G_N} (X_M\otimes \chi_{\underline s})\]
le morphisme de représentations de $G_N$ obtenu  par le foncteur d'induction.
On a aussi un morphisme obtenu de manière analogue : 
\[ h_{w_M}(\underline s):\;    \Ind_{P'}^{G_N}(w_M\cdot (   \left( X(1)\otimes  I(\underline\varepsilon',\lambda')   \right)      
\otimes \chi_{\underline s})) \longrightarrow 
\Ind_{P'}^{G_N} (w_M \cdot (X_M\otimes \chi_{\underline s})).\]
Soient $\Omega_M$ une fonctionnelle de Whittaker non nulle sur $X_M$. Soient  $\Omega(\underline s)$ et
 $\Omega'(\underline s)$ les 
 fonctionnelles de Whittaker  sur  $\Ind_{P=MN}^{G_N} (X_M\otimes \chi_{\underline s})$ et $\Ind_{P'}^{G_N} 
 (w_M \cdot( X_M\otimes \chi_{\underline s}))$
 respectivement obtenues par (\ref{IntWh}) et (\ref{WhMG}).
En composant respectivement par $h(\underline s)$ et $h_{w_M}(\underline s)$, on obtient des 
 fonctionnelles de Whittaker non nulles sur
  $ \Ind_{P=MN}^{G_N}(    \left( X(1)\otimes  I(\underline\varepsilon',\lambda')   \right) \otimes \chi_{\underline s})$ et 
 $\Ind_{P'}^{G_N}(w_M\cdot ( \left( X(1)\otimes  I(\underline\varepsilon',\lambda')   \right)\otimes \chi_{\underline s}))$.
  Notons-les $\omega(\underline s)$ et $\omega'(\underline s)$.
Les fonctionnelles de Whittaker  $\omega(\underline s)$ et $\omega'(\underline s)$ sont obtenues 
  par (\ref{IntWh}) et (\ref{WhMG}) à partir de la  fonctionnelle
 de Whittaker 
 $f \mapsto \Omega_M( \left[ g\mapsto h_M(f(g)) \right])   $ sur     $  X(1)\otimes  I(\underline\varepsilon',\lambda')   $.

Considérons le diagramme suivant :
\[ \xymatrix{  
 \Ind_{P=MN}^{G_N}\left( \left( X(1)\otimes  I(\underline\varepsilon',\lambda')   \right)\otimes \chi_{\underline s}\right)   
\ar[rrr]^{N(\underline s,  \left( X(1)\otimes  I(\underline\varepsilon',\lambda')   \right) )}
\ar[d]^{h(\underline s)} &&&  
  \Ind_{P'}^{G_N}(w_M\cdot ( \left( X(1)\otimes  I(\underline\varepsilon',\lambda')   \right) \otimes \chi_{\underline s}))
  \ar[d]^{h_{w_M}(\underline s)}\\
 \Ind_{P=MN}^{G_N} (X_M\otimes \chi_{\underline s}) \ar[rrr]^{N(\underline s, X_M) } &&& \Ind_{P'}^{G_N}
  (w_M \cdot (X_M\otimes \chi_{\underline s}))  } \]
  Nous laissons au lecteur le soin de vérifier que ce diagramme est commutatif. 

On a vu  que 
$\omega(\underline s)=\omega'(\underline s)\circ N(\underline s,  \left( X(1)\otimes  I(\underline\varepsilon',\lambda')   \right)).$
On a donc comme $\omega=\Omega\circ h$ et $\omega'=\Omega'\circ h_{w_M}$,
\begin{align*}
\Omega(\underline s)\circ h(\underline s)&=
\Omega'(\underline s)\circ h_{w_M}(\underline s)\circ N(\underline s,  \left( X(1)\otimes  I(\underline\varepsilon',\lambda')   \right)) \\
&=\Omega'(\underline s)\circ N(\underline s,X_M)\circ h(\underline s)
\end{align*}
et comme $h(\underline s)$ est surjective, on a bien $\Omega(\underline s)
=\Omega'(\underline s)\circ N(\underline s,X_M)$.
Ceci termine la démonstration de la proposition.
\qed

\subsection{Induction}  \label{InductionTord}
Notre but est maintenant de construire un foncteur d'induction $\widetilde \Ind_{\widetilde M}^{\widetilde G_N}$.
 La catégorie de départ de ce foncteur est celle des  représentations $(\tilde \pi_M , \pi_M)$ de l'espace tordu 
 $\widetilde M=M\rtimes \theta_M$ telles que $\pi_M$ soit de longueur finie et ait tous ses sous-quotients irréductibles dans 
$\scrE$. La catégorie d'arrivée est celle des représentations  de l'espace tordu $\widetilde G_N$.
Remarquons que le fait que $\pi_M$  soit $\theta_M$-stable
impose des conditions sur le caractère infinitésimal $\lambda'$ de $\GL_{N'}(\bbR)$ :
celui-ci doit être stable par changement de signe de toutes ses coordonnées, ce que l'on suppose désormais.

Soit $(\tilde \pi_M , \pi_M)$ une  telle représentation de l'espace tordu $\widetilde M=M\ltimes \theta_M$.
Posons $A_{\theta_M}=\tilde \pi(\theta_M)$, que l'on peut voir comme un opérateur d'entrelacement
 \[ A_{\theta_M} :\,  \pi_M \rightarrow \pi_M^{\theta_M}.\] 
 Soit $\caA$ le morphisme obtenu en appliquant le foncteur d'induction parabolique
$\Ind_{P=MN}^{G_N}$ : 
\[ \caA=\Ind_{P=MN}^{G_N}( A_{\theta_M}):\;\Ind_{P=MN}^ {G_N}(\pi_M) \longrightarrow   \Ind_{P=MN}^{G_N}(\pi_M^{\theta_M}).  \] 

On a aussi un opérateur d'entrelacement 
\[ \vartheta : \; \Ind_{P'}^{G_N} (\pi_M^{\theta_N})\simeq  \Ind_{P'}^{G_N}(w_M\cdot(\pi_M^{\theta_M}))\longrightarrow 
\left(\Ind_{P}^{G_N}(\pi_M)\right)^{\theta_N}, \qquad f\mapsto f\circ \theta_N. \]

On obtient par composition un morphisme $A_{\theta_N}$ : 
\begin{equation}\label{ATheta} \xymatrix{
\Ind_{P}^{G_N}(\pi_M) \ar[r]^{\caA\quad } & \Ind_{P=MN}^{G_N}(\pi_M^{\theta_M}) \ar[rr]^{N(0,\pi_M^{\theta_M})} && 
 \Ind_{P'}^ {G_N}(w_M\cdot(\pi_M^{\theta_M}) ) \ar[r]^{\vartheta}& \left(\Ind_{P}^{G_N}(\pi_M)\right)^{\theta_N}
}.\end{equation}

Définissons alors $(\tilde \pi, \pi)=\widetilde \Ind_{\widetilde M}^{\widetilde G_N}(\tilde \pi_M , \pi_M)$
par 
\begin{equation}
\pi=\Ind_P^{G_N}(\pi_M), \qquad \tilde \pi(\theta_N)=A_{\theta_N}
\end{equation}

Si $\varphi$ est un morphisme entre $(\tilde \pi_M , \pi_M)$  et $(\tilde \pi'_M , \pi'_M)$, on a le morphisme
 \[\Ind_{P}^{G_N}(\varphi): \, \pi=\Ind_{P}^{G_N}(\pi_M)\longrightarrow \pi'=\Ind_{P}^{G_N}(\pi'_M),\] et il faut voir 
 qu'il entrelace aussi les actions de $A_{\theta_N}= \tilde \pi(\theta_N) $ 
et  $A'_{\theta_N}= \tilde \pi'(\theta_N) $, c'est-à-dire que le diagramme 
\begin{equation}\label{AThetaFunc} \xymatrix{
\Ind_{P}^{G_N}(\pi_M) \ar[r]^{\caA\quad } \ar[d]^{\Ind_P^{G_N}(\varphi)}
& \Ind_{P=MN}^{G_N}(\pi_M^{\theta_M}) \ar[rr]^{N(0,\pi_M^{\theta_M})}\ar[d] && 
 \Ind_{P'}^ {G_N}(w_M\cdot(\pi_M^{\theta_M}) ) \ar[r]^{\vartheta}\ar[d]& \left(\Ind_{P}^{G_N}
 (\pi_M)\right)^{\theta_N}\ar[d]^{\Ind_P^{G_N}(\varphi)}\\
 \Ind_{P}^{G_N}(\pi'_M) \ar[r]^{\caA\quad } & \Ind_{P=MN}^{G_N}({\pi'_M}^{\theta_M}) \ar[rr]^{N(0,{\pi'_M}^{\theta_M})} && 
 \Ind_{P'}^ {G_N}(w_M\cdot({\pi'_M}^{\theta_M}) ) \ar[r]^{\vartheta}& \left(\Ind_{P}^{G_N}(\pi'_M)\right)^{\theta_N}
 }\end{equation}
commute. 

Commençons par le carré de gauche. La flèche verticale  à droite de celui-ci est encore $\Ind_P^{G_N}(\varphi)$
car $\varphi$ est aussi un opérateur d'entrelacement entre $\pi_M^{\theta_M}$ et ${\pi'_M}^{\theta_M}$.
Les opérateurs $\caA$ sont obtenus par composition à gauche par $A_{\theta_M}$, les opérateurs
$\Ind_P^{G_N}(\varphi)$ par composition à gauche par $\varphi$ et la commutativité du diagramme
vient donc du fait que $A_{\theta_M}$ et $\varphi$ commutent.
Passons ensuite au carré central. La flèche verticale de droite est cette fois $\Ind_{P'}^{G_M}(\varphi)$.
Il est équivalent de vérifier la commutativité du diagramme où l'on remplace 
$\pi_M^{\theta_M}$ et ${\pi'_M}^{\theta_M}$ par $\pi_M$ et ${\pi'_M}$.
On vérifie la commutativité  par unicité du prolongement analytique et calcul direct avec les 
 formules intégrales dans le domaine de convergence :
 soit $f\in \Ind_{P=MN}^{G_N}(\pi_M)$.    On a alors 
\begin{equation}
\left(M(\underline s, \pi_M)(f_{\underline s})\right)(g)=\int_{N'}
 f_{\underline s}(w_M^{-1}n'g)\; dn' \end{equation}
et 
\begin{equation}\label{tres}
\Big( \Ind_{P'}^{G_N}(\varphi)\big( M(\underline s,\pi_M)(f_{\underline s})\big)\Big)(g)
=\varphi \left( \int_{N'}  f_{\underline s}(w_M^{-1}n'g)\; dn'\right) .\end{equation}
D'autre part 
\begin{equation}
\big( \Ind_{P}^{G_N}(\varphi) (f_{\underline s})\big)(g)
=\varphi( f_{\underline s}(g)) \end{equation}
et 
\begin{equation}\label{quatro}
\Big( M(\underline s, \pi'_M )\big( \Ind_{P}^{G_N}(\varphi) \left( f_{\underline s} \right) \big)\Big)(g) =  
 \int_{N'}  \varphi(f_{\underline s}(w_M^{-1}n'g))\; dn' \end{equation}
La commutativité du diagramme vient de  l'égalité des membres de droite de (\ref{tres}) et (\ref{quatro}), évidente par linéarité de l'intégrale, 
que l'on multiplie par le facteur $r(s)$ et que l'on évalue  en $\underline s=0$, 

Enfin, pour le carré de droite, on fait de même : 
soit $f\in \Ind_{P'}^{G_N}(w_M\cdot(\pi_M^ {\theta_M}))$.    On a alors 
\begin{equation}
\left(\vartheta(f)\right)(g)= f(\theta_N(g)) \end{equation}
et 
\begin{equation}\label{cinque}
\Big( \Ind_{P}^{G_N}(\varphi)\big( \vartheta (f)\big)\Big)(g)
=\varphi \left( f(\theta_N(g))   \right) .\end{equation}
D'autre part 
\begin{equation}
\big( \Ind_{P'}^{G_N}(\varphi)  (f)\big)(g)
=\varphi( f(g)) \end{equation}
et 
\begin{equation}\label{seis}
\Big( \vartheta \big( \Ind_{P'}^{G_N}(\varphi) \left( f \right) \big)\Big)(g)
= \varphi( f(\theta_N(g))).
   \end{equation}
Ceci termine la démonstration de la commutativité du diagramme. \qed

\begin{rmq}\label{exact}
Il découle de l'exactitude du foncteur $\Ind_P^ {G_N}$ que le foncteur $\widetilde \Ind_{\widetilde M}^{\widetilde G_N}$
est lui aussi exact.
\end{rmq}

 \begin{lemme} \label{bibi} Soient $X_M$ un module standard  $\theta_M$-stable de $M$
 dont tous les sous-quotients irréductibles sont dans $\scrE$ et $\Omega_M$ 
 une fonctionnelle de Whittaker non nulle sur   $X_M$. Soit 
 $A_{\theta_M}: X_M\rightarrow X_M^{\theta_M}$ un opérateur d'entrelacement 
 permettant de définir la représentation tordue $(\widetilde X_M,X_M)$ de $\widetilde M$. 
 Supposons que $A_{\theta_M}$ soit tel que $\Omega_M=\Omega_M \circ A_{\theta_M}$.
 Formons $(\tilde \pi,\pi)=\widetilde \Ind_{\widetilde M}^{\widetilde G_N}(\widetilde X_M,X_M)$ comme ci dessus, avec
  $\tilde \pi(\theta_N)=A_{\theta_N}$. Soit $\Omega$ la  fonctionnelle de Whittaker  sur  
 $\Ind_{P=MN}^{G_N} (X_M)$   obtenue par (\ref{IntWh}) et (\ref{WhMG}).
  Alors  $\Omega=\Omega \circ A_{\theta_N}$.
\end{lemme} 
 
Par unicité à un scalaire près des fonctionnelles de Whittaker sur le module standard
$X_M$, on voit que quitte à multiplier $A_{\theta_M}$ par un scalaire non nul, on peut toujours 
effectivement supposer que  $\Omega_M=\Omega_M \circ A_{\theta_M}$.
 
 \medskip
 
\dem La fonctionnelle de Whittaker $\Omega_M$ est aussi une fonctionnelle de Whittaker pour  $X_M^{\theta_M}$, 
à partir de laquelle on obtient des fonctionnelles de Whittaker $\Omega_{\theta_M}$ et  $\Omega'_{\theta_M}$
par  (\ref{IntWh}) et (\ref{WhMG}) sur $\Ind_{P}^ {G_N}(X_M^{\theta_M})$ et  
$\Ind_{P'}^ {G_N}(w_M\cdot(X_M^{\theta_M}))$ respectivement.
 D'après la proposition \ref{OpEntrBiz} 
 \[  \Omega_{\theta_M}=\Omega'_{\theta_M} \circ N(0,X_M^{\theta_M}). \]
On  montre facilement que $\Omega_{\theta_M}\circ \caA =\Omega$ et 
 $\Omega\circ \vartheta=\Omega'_{\theta_M}$, d'où   
\[\Omega=\Omega_{\theta_M}\circ \caA 
   = \Omega \circ \vartheta \circ N(0,X_M^{\theta_M}) \circ \caA = \Omega\circ A_{\theta_N}.\] 
   \qed

\begin{lemme}\label{colle}
 Si $\pi_M \in \scrE$ (en particulier $\pi_M$ est irréductible) est $\theta_M$-stable, et si 
$A_{\theta_M}$ est l'opérateur d'entrelacement permettant de construire
l'extension canonique $\pi_M^ +$ de $\pi_M$ à $M^ +$ comme dans la section \ref{NormExt}, alors $A_{\theta_N}$
 est l'opérateur d'entrelacement permettant de construire
l'extension canonique $\Ind_P^{G_N}(\pi_M)^+$ de $\Ind_P^{G_N}(\pi_M)$ à  $G_N^ +$. 
\end{lemme}

\dem Soit $X_M$ le module standard de $M$ admettant $\pi_M$ comme sous-module irréductible.
 C'est un produit tensoriel $X_1\otimes X'$ de modules standard respectivement de $\GL_{N_1}(\bbR)$ et $\GL_{N'}(\bbR)$.
Le corollaire \ref{passtdstd} nous dit que $\Ind_P^{G_N}(X_M)$ est isomorphe à une représentation standard
$X$ de $G_N$. Notons $L:\, \Ind_P^ {G_N}(X_M) \rightarrow X$ cet isomorphisme. Le module 
$X$ étant un module standard $\theta_N$-stable, il est muni de l'opérateur d'entrelacement canonique 
$A_{X}$ défini dans la section  \ref{NormExt}. L'opérateur  d'entrelacement $A_{\theta_N}: \Ind_P^ {G_N}(X_M) \rightarrow 
 \left(\Ind_P^ {G_N}(X_M) \right)^{\theta_N} $ 
est celui du lemme précédent. Il s'agit alors de montrer la  commutativité du diagramme :
\[ \xymatrix{
\Ind_{P}^{G_N}(\pi_M) \ar[rr]\ar[d]^{A_{\theta_N}} && \Ind_{P}^{G_N}(X_M) \ar[rr]^{L} \ar[d]^{A_{\theta_N}} && 
X \ar[d]^{A_{X}}\\
 \left( \Ind_{P}^{G_N}(\pi_M) \right)^{\theta_N} \ar[rr] && \left( \Ind_{P}^{G_N}(X_M)\right)^{\theta_N} \ar[rr]^{L} 
  && X^{\theta_N} 
 }\]  
 ce que nous laissons au lecteur.\qed

\section{Paquets d'Adams-Johnson}\label{ParAJ}

\subsection{Paramètres d'Adams-Johnson}

Nous redonnons dans cette section des éléments concernant les paramètres d'Arthur
considérés par Adams et Johnson \cite{AdJo}. 
Nous renvoyons aux exposés d'Arthur \cite{Art89}, 
Kottwitz \cite{Kott} et Taïbi \cite{Tai}.  Soit $\mathbf G$ un groupe algebrique connexe réductif défini sur $\bbR$
et que l'on suppose quasi-déployé. 
Les paramètres d'Adams-Johnson sont  les paramètres d'Arthur
\[\psi : \, W_\bbR \times \SL_2(\bbC) \longrightarrow {}^L G \]
vérifiant des conditions  que nous ne redonnons pas ici; nous nous contentons d'en donner
certaines conséquences, la plus importante étant que le caractère infinitésimal de $G$ qui leur est associé
 est celui d'une représentation irréductible de dimension finie $F$, en particulier, il est entier et régulier.

Fixons un épinglage $\mathbf{spl}_{\widehat G}=(\caB,\caT, \{\caX_\alpha\})$ de $\widehat G$ et supposons que 
le $L$-groupe ${}^LG$ est construit via cet épinglage, c'est-à-dire que l'on suppose que l'action
de $W_\bbR$ sur $\widehat G$ préserve  $\mathbf{spl}_{\widehat G}$.
Si $\psi$ est un paramètre d'Adams-Johnson pour $G$, alors il existe un sous-groupe de levi $\caL$ de $\widehat G$
contenant $\caT$ (c'est le centralisateur de $\psi(\bbC^\times)$), un sous-groupe algébrique connexe réductif
 $\mathbf L_*$ de $\mathbf G$ défini sur $\bbR$ dont le dual de Langlands
est $\caL$, et un plongement de $L$-groupes
\[\iota_{L,G} : {}^LL\longrightarrow {}^LG \]
tels que $\psi$ se factorise via ce plongement par un paramètre d'Arthur $\psi_L$ de $L_*$, c'est-à-dire
$\psi=\iota_{L,G}\circ \psi_L$. On suppose que le $L$-groupe ${}^LL$ est lui aussi construit via 
 l'épinglage $\mathbf{spl}_{\caL}=(\caB\cap \caL,\caT, \{\caX_\alpha\})$. 
D'autre part, ce $\psi_L$ est le paramètre d'Arthur d'un caractère unitaire, c'est-à-dire qu'il est obtenu 
à partir d'un morphisme $\SL_2(\bbC)\rightarrow \caL$ principal et d'un élément de $H^1(W_\bbR,Z(\caL))$.
Nous en dirons plus dans la section suivante sur ce sous-groupe $\mathbf L_*$, ainsi que sur ses formes intérieures
$\mathbf L$ lorsque nous décrirons les paquets d'Adams-Johnson.

Le paramètre d'Adams-Johnson $\psi$ détermine aussi un sous-groupe parabolique 
 $\caQ=\caL\caU$ de $\widehat G$, dont  $\caL$ est facteur de Levi de $\caQ$.

\subsection{Paquets d'Adams-Johnson}
\label{AJpackSec}
  Décrivons maintenant le paquet  $\Pi^{AJ}_\psi$ attaché par Adams et Johnson à un paramètre
  $\psi$ comme dans la section précédente.
  Fixons une paire de Borel $(\mathbf{B},\mathbf{T})$ de $\mathbf{G}$ où $\mathbf{T}$ est un tore maximal
  défini sur $\bbR$ et maximalement anisotrope.
Soient   $\Sigma_\mathbf{B}$ l'ensemble des racines simples du système de racines positives 
$R^+_{\mathbf B}=R(\mathbf T,\mathbf B)$
et  $\Sigma_\caB$ l'ensemble des racines simples du système de racines positives $R^+_{\caB}=R(\caT,\caB)$.
 On a un isomorphisme canonique entre les données radicielles basées
  \[ (X^*(\mathbf{T}), \Sigma_\mathbf{B}, X_*(\mathbf{T}), \Sigma_\mathbf{B} \check{}\, )\quad  \text{ et } \quad
 (X_*(\caT), \Sigma_\caB \check{}\, , X^*(\caT), \Sigma_\caB)   \]
  et l'on peut associer au sous-groupe parabolique $\caQ=\caL\caU$ un sous-groupe parabolique
   $\mathbf{Q}=\mathbf{L}\mathbf{U}$ de $\mathbf{G}$ contenant $\mathbf{T}$.

Notons $\mathbf{L}$ le sous-groupe de $\mathbf{G}$ ainsi défini. Il est défini sur $\bbR$ et l'on note  $L$ le groupe de ses points réels. Considérons l'ensemble $\caS_\caQ$ des classes de conjugaison   (sous $G$)
 de couples $(\mathbf{Q},\mathbf{L})$
  où $\mathbf{Q}$ est un sous-groupe parabolique de  $\mathbf{G}$ et  $\mathbf{L}$ est un sous-groupe de Levi
  défini sur $\bbR$ obtenus de cette façon. Adams et Johnson démontrent que parmi ces classes de conjugaison, 
  l'une au moins $(\mathbf{Q_*},\mathbf{L_*})$ est telle que $\mathbf L_*$ est quasi-déployé. C'est ce $\mathbf L_*$
  que nous avons introduit dans la section précédente.
  
 L'ensemble $\caS_\caB$ des classes de conjugaison de couples $(\mathbf{B},\mathbf{T})$, peut être identifié
 à $W(G,T)\backslash W(\mathbf{G},\mathbf{T})$ en choisissant un point de base. La surjection 
 naturelle de $\caS_\caB$ vers $\caS_\caQ$ identifie alors ce dernier avec 
 \[ W(G,T)\backslash W(\mathbf{G},\mathbf{T})^\tau/W (\mathbf{L},\mathbf{T}).\]
 (Voir la section 10 de \cite{AdJo} pour expliquer l'apparition du $\tau$ lorsqu'on ne suppose pas que 
 le rang de $G$ est égal au au rang de $K$ et qui correspond à la condition  que  $\mathbf L$ soit défini sur $\bbR$.)
Remarquons que cet ensemble s'identifie aussi naturellement avec 
$\ker[ H^1(\bbR,\mathbf{L})\rightarrow   H^1(\bbR,\mathbf{G})]$.

Pour une classe  $(\mathbf{Q},\mathbf{L})$ fixée dans $\caS_\caQ$, on a un isomorphisme canonique
$\caL\simeq \widehat L$.  Pour toute autre classe $(\mathbf{Q'},\mathbf{L'})$   dans $\caS_\caQ$, il existe un unique 
élément $g\in \mathbf{G}(\bbC)/ \mathbf{L}(\bbC)$ tel que $g\cdot(\mathbf{Q},\mathbf{L})=(\mathbf{Q'},\mathbf{L'})$
 (action par conjugaison),  qui donne un isomorphisme canonique entre L-groupes ${}^LL'\simeq {}^LL$.
 L'isomorphisme $\widehat L\simeq \caL$ s'étend en un $L$-plongement 
 \begin{equation}\label{iotaLG} \iota_{L,G}:\, {}^LL\rightarrow {}^LG \end{equation}
 déjà introduit dans la section précédente et qui factorise le paramètre $\psi$ en $\psi=\iota_{L,G}\circ \psi_L$. De plus  
  $  \psi_L $  détermine une représentation unitaire de dimension $1$  de $L$, que nous notons $\lambda_ {\psi,\mathbf{Q},\mathbf{L}}$.

 En appliquant le foncteur d'induction cohomologique de Vogan-Zuckerman (voir par exemple \cite{KV}),
  Adams et Johnson définissent une représentation 
 \begin{equation}\label{pipsiQL}\pi_ {\psi,\mathbf{Q},\mathbf{L}} =\caR^i_\frqqq(\lambda_ {\psi,\mathbf{Q},\mathbf{L}})
 \end{equation}
  de $G$, où 
 $\frqqq$ est l'agèbre de Lie de $\mathbf{Q}$ et $i$ 
est la dimension de la partie compacte du radical unipotent $\mathbf{U}$ de $\mathbf{Q}$, aussi donné par la formule  
  $$\frac{1}{2} (\dim G- \dim L)-(q(G)-q(L)).$$ 
 La définition de $q$ est donnée ci-dessous. On sait d'après Vogan  (voir \cite{KV}, Chapter 9) que cette représentation est unitaire.

\begin{rmq}\label{Lqd} Si   $(\mathbf{Q},\mathbf{L}) \in \caS_\caQ$ est tel que $\mathbf{L}$ est quasi-déployé, 
alors le paramètre de Langlands de  $\lambda_ {\psi,\mathbf{Q},\mathbf{L}}$ est $\phi_{\psi_{\mathbf{Q},\mathbf{L}}}$
({\sl cf. }  (\ref{ArtLang})) et celui de $\pi_ {\psi,\mathbf{Q},\mathbf{L}}$ est $\phi_\psi=\iota_{L,G}\circ \phi_{\psi_{\mathbf{Q},\mathbf{L}}}$. 
L'ensemble des  $\pi_ {\psi,\mathbf{Q},\mathbf{L}}$ avec 
$(\mathbf{Q},\mathbf{L}) \in \caS_\caQ$  tel que $\mathbf{L}$ est quasi-déployé est le paquet de Langlands de $G$ 
 attaché à  $\phi_\psi$.
\end{rmq}  

\medskip 

Nous pouvons maintenant définir le paquet d'Adams-Johnson de paramètre $\psi$.
\begin{defi} \label{AJpack}
 Le paquet d'adams-Johnson de paramètre $\psi$ est, avec les notations qui précèdent :
 \[  \Pi_\psi^{AJ}=\{  \pi_ {\psi,\mathbf{Q},\mathbf{L}}\, \vert \, (\mathbf{Q},\mathbf{L}) \in \caS_\caQ  \}.\]
\end{defi}

L'une des conditions devant être vérifiée par les paquets conjecturaux d'Arthur est que  ceux-ci doivent être le support
d'une distribution stable. Adams et Johnson montrent que tel est le cas pour les paquets 
qu'ils définissent. Enonçons ceci de manière précise. Nous avons besoin de la définition suivante.

\begin{defi}\label{rkKrkG}
Soit $\mathbf G$  un groupe algébrique connexe réductif défini sur $\bbR$. Soit $c_0^G$ la moitié de la dimension
de la partie déployée d'un sous-groupe de Cartan fondamental de $G$,  
et posons 
\begin{equation}\label{qGdef}
 q(G)=\frac{1}{2}  (\dim G-\dim K)-c_0^ G.
\end{equation}
C'est un entier, et  dans le cas où les rang de $G$ et de $K$ sont égaux (autrement dit, $G$ est forme intérieure
 d'une forme  compacte),  $q(G)=\frac{1}{2}  (\dim G-\dim K)$. Si $\mathbf G^*$ est une forme intérieure quasi-déployée
 de $G$, $e(G)= (-1)^{q(G^*)-q(G)}$ est le signe de Kottwitz de $G$.
\end{defi}

 Posons alors 
\begin{equation}\label{StablePiAJ}
 [\Pi_\psi^{AJ}]^{st} : =\sum_{ (\mathbf{Q},\mathbf{L}) \in \caS_\caQ  } e(L) \, 
[\pi_ {\psi,\mathbf{Q},\mathbf{L}}], \qquad 
 \Theta_{\Pi_\psi^{AJ}}^{st} : =\sum_{ (\mathbf{Q},\mathbf{L}) \in \caS_\caQ  }e(L) \, 
 \Theta_{\pi_ {\psi,\mathbf{Q},\mathbf{L}}}.
\end{equation} 
La première expression est une représentation virtuelle et la seconde est son caractère.

\begin{thm}(\cite{AdJo}, Thm. 2.13) 
 Soit $\psi$ un paramètre d'Adams-Johnson pour $G$, et  $\Pi_\psi^{AJ}$ le paquet associé. 
La distribution $\Theta_{\Pi_\psi^{AJ}}^{st} $ est stable.
\end{thm}

Adams et Johnson (\cite{AdJo}, Thm. 2.21)  montrent   aussi que 
les paquets $\Pi_\psi^{AJ}$ vérifient bien les identités endoscopiques (non tordues) attendues.

 \subsection{Formule des caractères stables pour les paquets d'Adams-Johnson}
 \label{ResJohPaq}
 Nous expliquons dans ce paragraphe comment Adams et Johnson expriment la représentation virtuelle stable (\ref{StablePiAJ})
 en terme de  représentations virtuelles stables associée aux pseudo-paquets de Langlands.
On reprend les notations de la section \ref{AJpackSec} : $\psi$ est un paramètre d'Adams-Johnson pour $G$ 
 et
 \[  \Pi_\psi^{AJ}=\{  \pi_ {\psi,\mathbf{Q},\mathbf{L}}\, \vert \, (\mathbf{Q},\mathbf{L}) \in \caS_\caQ  \}.\]
 est le paquet d'Adams-Johnson qui lui est associé. Le caractère infinitésimal de ces représentations  est celui 
 d'une représentation de dimension finie $F$ de $G$. 
 
 Dans \cite{Joh1},  Johnson montre qu'il existe des résolutions des 
 représentations   $\pi_ {\psi,\mathbf{Q},\mathbf{L}}$ constituant les paquets d'Adams-Johnson
par des représentations standards. Les résolutions des représentations de Speh en (\ref{ResSpeh0})
en sont un exemple.
Ces résolutions donnent des égalités de représentations virtuelles de la forme 
\[   [\pi_ {\psi,\mathbf{Q},\mathbf{L}}] =
  \sum_{\gamma \in \scrP_{ (\mathbf{Q},\mathbf{L}) }} (-1)^{\ell_V(\pi_ {\psi,\mathbf{Q},\mathbf{L}})-\ell_V(\gamma)} [X(\gamma)].\]
  Ici  $\scrP_{ (\mathbf{Q},\mathbf{L})} $ est l'ensemble des paramètres des modules standard
  apparaissant dans la résolution de Johnson de $\pi_ {\psi,\mathbf{Q},\mathbf{L}}$ ({\sl cf.} \cite{Joh1}). Ces modules standard sont donnés
  explicitement dans le Theorem 8.2 de   \cite{AdJo}.
Le signe est donné par la longueur de  Vogan  $\ell_V$ (\cite{VIC3}, \cite{VIC4}), et l'on  a alors 
   $l_V(\pi_ {\psi,\mathbf{Q},\mathbf{L}}) =q(L)$. 
  
On tire de  ceci l'égalité de représentations virtuelles 
\begin{align}\label{CES}
 \nonumber [\Pi_\psi^{AJ}]^{st}&=\sum_{ (\mathbf{Q},\mathbf{L}) \in \caS_\caQ  } e(L)
 [\pi_ {\psi,\mathbf{Q},\mathbf{L}}] 
 =\sum_{ (\mathbf{Q},\mathbf{L}) \in \caS_\caQ  } e(L)
  \sum_{\gamma \in \scrP_{ (\mathbf{Q},\mathbf{L}) }} (-1)^{\ell_V(\pi_ {\psi,\mathbf{Q},\mathbf{L}})-\ell_V(\gamma)} [X(\gamma)]\\
  &=(-1)^{q(L^*) }\sum_{ (\mathbf{Q},\mathbf{L}) \in \caS_\caQ  }  
  \sum_{\gamma \in \scrP_{ (\mathbf{Q},\mathbf{L}) }} (-1)^{\ell_V(\gamma)} [X(\gamma)].
   \end{align} 
  Il est démontré dans \cite{AdJo} que tous les $L$ intervenant dans le somme sont dans la même classe de forme intérieures.
  Notons $L^*$ une forme quasi-déployée. 
On a utilisé $e(L)=(-1)^{q(L)-q(L^*)}$ pour faire sortir la constante $(-1)^{q(L^*)}$.   
\medskip  
 
 Considérons un module standard $X$ intervenant dans la résolution de Johnson de l'un des 
  $ \pi_ {\psi,\mathbf{Q},\mathbf{L}}$, disons  de $ \pi_ {\psi,\mathbf{Q_1},\mathbf{L_1}}$, 
  $  (\mathbf{Q_1},\mathbf{L_1}) \in \caS_\caQ$.  
Soit $X'$ un module standard appartenant au m\^eme pseudo-paquet que $X$    
  ({\sl cf.} Définition \ref{pseudopaquets}). Alors d'après \cite{AdJo}, Lemma 8.8
  il existe un unique  $  (\mathbf{Q},\mathbf{L}) \in \caS_\caQ$
  tel que $X'$ apparaisse dans la résolution de Johnson de  $\pi_ {\psi,\mathbf{Q},\mathbf{L}}$.

  Une propriété bien connue des pseudo-paquets 
  (voir \cite{VIC4}, ou \cite{ADC}) est que la longueur de deux éléments du même pseudo-paquet
  est la même. De la discussion ci-dessus et de (\ref{CES}) on tire 
   \begin{align}\label{RepVStableAJ}
  [\Pi_\psi^{AJ}]^{st}&=(-1)^{q(L^*)}
  \sum_{\phi \in \Phi(\psi) } (-1)^{\ell_V(\phi)} [X_\phi]
  \end{align}
 où   $\Phi(\psi)$ est  l'ensemble des paramètres de Langlands des modules standard 
intervenant dans la résolution de l'un des  $\pi_ {\psi,\mathbf{Q},\mathbf{L}}$,   
$\ell_V(\phi)$ la longueur d'un de ses éléments et $[X_\phi]:= \sum_{X \in \widetilde \Pi_\phi} [X]$. 

Si l'on réécrit ceci comme une identité de distributions, on obtient
  \begin{equation}\label{IdStableAJ}
\Theta_{\Pi_\psi^{AJ}}^{st}= (-1)^{q(L^*)}  \sum_{\phi \in \Phi(\psi) } (-1)^{\ell_V(\phi)} \Theta_ {\widetilde \Pi_ \phi} , 
  \end{equation}
où toutes les distributions sur $G$ intervenant dans cette équation sont stables.

\begin{rmq}\label{cryptic}
Supposons $G$ quasi-déployé. On sait alors qu'il existe $(\mathbf{Q}_*, \mathbf{L}_*)\in \caS_\caQ$ tel que 
$ \mathbf{L}_*$ soit aussi quasi-déployé. Il découle alors de la discussion ci-dessus, de la remarque \ref{Lqd}
 et de la description
explicite des modules standard intervenant dans les résolutions des $ \pi_ {\psi,\mathbf{Q},\mathbf{L}}$ 
que la distribution stable $\Theta_{\Pi_\psi^{AJ}}^{st}$ est obtenue en stabilisant l'écriture  en modules
 standard de la somme des représentations dans le paquet de Langlands 
 paramétré par $\phi_\psi$.
\end{rmq}

\section{Groupes classiques et endoscopie tordue}\label{GCET}

\bigskip
\subsection{Les  groupes classiques}\label{grcla}
  
  Les \og groupes classiques \fg \; que nous considérons sont ceux qui apparaissent
 dans les données endoscopiques elliptiques simples des $\widetilde G_N$ selon la terminologie d'Arthur {\sl cf.}
 \cite{Art13}, \S I.2.   Le plus commode est encore d'en faire la liste, et de fixer quelques notations pour
  pouvoir s'y référer  facilement. Pour $n\in \bbN^\times$, on considère les groupes de rang $n$ suivants: 
 
 \begin{itemize}
 \item[\bf A.] Le groupe symplectique $\Sp(2n,\bbR)$.  
 
 \noindent C'est un groupe déployé. Son dual de Langlands est $\SO(2n+1,\bbC)$ et son $L$-groupe est le 
 produit direct $\SO(2n+1,\bbC)\times W_\bbR$. 
 
  \item[\bf B.] Le groupe spécial orthogonal impair $\SO(n,n+1 )$. 
  
  \noindent  C'est un groupe déployé. Son dual de Langlands 
  est $\Sp(2n,\bbC)$ et son $L$-groupe est le   produit direct $\Sp(2n,\bbC)\times W_\bbR$.
   
    \item[\bf C.]  Le groupe spécial orthogonal pair $\SO(n,n)$. 
   
   \noindent C'est un groupe déployé. Son dual de Langlands  est $\SO(2n,\bbC)$ et son $L$-groupe est le 
 produit direct $\SO(2n,\bbC)\times W_\bbR$.

    \item[\bf D.]  Le groupe spécial orthogonal pair $\SO(n-1,n+1)$. 
  
   \noindent  C'est un groupe quasi-déployé. Son dual de Langlands  est $\SO(2n,\bbC)$ et son $L$-groupe est le 
 produit semi-direct $\SO(2n,\bbC)\rtimes W_\bbR$.
    \end{itemize}

Pour chacun de ces groupes, on dispose d'une représentation naturelle du $L$-groupe dans un $\widehat \GL_N$ : 
\begin{equation}\label{DefStdG} \Std_G: \; {}^LG \longrightarrow  \GL_N(\bbC) \end{equation}
Dans le cas $\mathbf A$, elle est donnée l'inclusion de $\SO(2n+1,\bbC)$ dans $\GL_{2n+1}(\bbC)$, dans le cas 
$\mathbf B$, de l'inclusion de $\Sp(2n,\bbC)$ dans $\GL_{2n}(\bbC)$, dans le cas $\mathbf C$,
de l'inclusion de $\SO(2n,\bbC)$ dans $\GL_{2n}(\bbC)$. Le cas $\mathbf D$, est plus délicat
car le groupe étant non déployé, le $L$-groupe est un produit semi-direct non trivial
$\SO(2n,\bbC)\rtimes W_\bbR$. Mais l'action de $W_\bbR$ sur $\SO(2n,\bbC)$ est donnée par l'action d'un élément de 
$\Or(2n,\bbC)$, de sorte que l'on a un morphisme
\[{}^L\SO(n-1,n+1)=  \SO(2n,\bbC)\rtimes W_\bbR  \rightarrow  \Or(2n,\bbC), \]
et la composition avec  l'inclusion de $\Or(2n,\bbC)$ dans $\GL_{2n}(\bbC)$ nous donne la représentation voulue.
On a donc $N=2n$ dans les cas $\mathbf B$, $\mathbf C$, $\mathbf D$, et  $N=2n+1$ dans les cas $\mathbf A$. 

\medskip

Fixons $\mathbf G$ comme ci-dessus.
 La donnée de $(\mathbf G,\Std_G)$ est celle d'une {\sl  donnée endoscopique tordue}
elliptique pour $(G_N,\theta_N)$. Nous renvoyons le lecteur à \cite{KS} et \cite{Art13}
pour tout ce qui concerne la théorie de l'endoscopie tordue. Rappelons seulement
que dans une telle situation, Kottwitz et Shelstad définissent un  facteur de transfert, permettant de définir une 
application $\mathrm{Trans_{geo}}$ (\og {\sl transfert géométrique} \fg) entre l'espace des intégrales orbitales  
sur $\widetilde G_N$ et  l'espace des intégrales  intégrales orbitales stables sur $G$.
 Ceci est démontré par D. Shelstad dans \cite{Shel}.
Ce facteur de transfert n'est défini {\sl a priori } qu'à une constante multiplicative près, mais le choix de la donnée
 de Whittaker sur $G_N$ ({\sl cf.} section \ref{NormExt}) permet de fixer cette constante 
 (\cite{KS}, section 5.3), ce que nous   supposons fait dans la suite.

\medskip

Par dualité,  ce transfert d'intégrales orbitales définit une application entre espaces
de distribution invariantes : 
\begin{equation}\label{Trans}
\mathrm{Trans}_G^{\widetilde G_N}: \; \mathrm{Dist}(G)^{st}\longrightarrow \mathrm{Dist}(\widetilde G_N)^{G_N}  ,
\end{equation}
l'espace de départ étant celui  des distributions stables sur $G$ 
et l'espace d'arrivée  celui des distributions sur $\widetilde G_N$, invariantes sous l'action adjointe de $G_N$.

Soit $\psi_G \in \Psi(G)$ un paramètre d'Arthur pour le groupe $G$. Posons 
$\psi=\Std_G \circ \psi_G$. C'est un paramètre d'Arthur (auto-dual) pour le groupe 
$G_N$. Soit $\Pi_\psi$ la représentation irréductible autoduale de $G_N$ associée à ce paramètre
({\sl cf.} proposition \ref{refreq}). Rappelons qu'à un tel paramètre $\psi_G$, Arthur a  associé dans  \cite{Art13}
un paquet $\Pi_{\psi_G}$ de représentations unitaires de $G$ caractérisé par  les identités
 endoscopiques usuelles et l'identité endoscopique tordue que nous allons décrire, reliant le paquet 
  $\Pi_{\psi_G}$ à la représentation $\Pi_\psi$ de $G_N$. Tout d'abord, un tel paquet est le support d'une 
  distribution stable sur $G$, c'est-à-dire qu'il existe une combinaison linaire (à  coefficients  
dans $\bbZ$) des caractères des représentations dans 
$\Pi_{\psi_G}$, décrite explicitement par Arthur, disons
\[  \Theta_{\Pi_{\psi_G}}^{st}=\sum_{\pi \in \Pi_{\psi_G}} a_\pi \; \Theta_\pi   \] 
qui est une distribution stable. Cette distribution 
  se transfère donc par l'application (\ref{Trans})  en une distribution $\mathrm{Trans}( \Theta_{\Pi_{\psi_G}})$.
L'identité endoscopique est alors 
\begin{equation}\label{IdEnd}
\mathrm{Trans_G^{\widetilde G_N}}( \Theta_{\Pi_{\psi_G}}^{st})=\Tr_{\theta_N}(\Pi_\psi)
\end{equation}  
  le membre de droite étant la trace tordue  ({\sl cf.} Définition \ref{TrT}), normalisée par la donnée de Whittaker.

Lorsque le paramètre $\psi_G$ est trivial sur le facteur $\SL_2(\bbC)$, 
le paquet $\Pi_{\psi_G}$ est un paquet de Langlands tempéré, la distribution est simplement la somme  
des caractères $\Theta_\pi$, $\pi \in \Pi_{\psi_G}$, et l'identité (\ref{IdEnd}) est alors
démontrée par P. Mezo \cite{mezo}, à un facteur multiplicatif près. Nous démontrons dans l'annexe \ref{Moe}
qu'en fait  l'identité (\ref{IdEnd}) est valide, autrement dit que le facteur multiplicatif restant à déterminer
dans Mezo est $1$.
Soit $\phi_G$ un paramètre de Langlands pour $G$. Le résultat de Mezo, précisé dans l'annexe \ref{Moe},
 et le fait que  le transfert endoscopique commute à l'induction entraîne le résultat suivant
 pour les pseudo-paquets ({\sl cf.} Définition \ref{pseudopaquets}). 
\begin{prop}\label{TransPsPaq}
Posons $\phi=\Std_G \circ \phi_G$. On a alors 
 \[ \mathrm{Trans}(  \Theta_{ \widetilde \Pi_{\phi_G}})=\Tr_{\theta_N}(\widetilde \Pi_\phi). \]
  \end{prop}

    \subsection{Paramètres et paquets d'Adams-Johnson des groupes classiques} \label{AJcla}

  Soit $G$ l'un des groupes classiques de rang $n$ de la section \ref{grcla} et comme dans cette section, soient
   \[\psi_G :\; W_\bbR\times \SL_2(\bbC)\longrightarrow {}^L G\] 
     un paramètre d'Arthur  pour $G$,  $\psi=\Std_G \circ \psi_G$ et  $\Pi_{\psi}$ la représentation auto-duale 
   de $G_N$ attachée au paramètre d'Arthur $\psi$.

On suppose maintenant de plus que $\psi_G$ est un paramètre d'Adams-Johnson.
   On reprend les notations de la section \ref{AJpackSec} pour les objets  
   attachés au  paramètre $\psi_G$ et celle de la section \ref{PaqArtGN} pour les objets attachés à $\psi$.  
   On explicite ceux-ci pour  chaque famille de groupes classiques de \ref{grcla}.
   
Quitte à remplacer $\psi_G$ par un paramètre équivalent, on peut supposer que le    
    centralisateur $\caL$ de $\psi_G(\bbC^\times)$ dans $\widehat G$ est de la forme : 
   \begin{equation}\label{LeviG} \caL\simeq \caL_1\times \cdots \times \caL_r \end{equation}
où pour tout $i=1, \ldots, r-1$,  $\caL_i\simeq \GL(n_i,\bbC)$, $n_i \in \bbN^\times$
 et suivant les cas :
 \[\caL_r\simeq \begin{cases} \SO( 2n_r+1,\bbC) \quad \text{(cas {\bf A})},\\
  \Sp( 2n_r,\bbC) \quad \text{(cas {\bf B})},\\  
  \SO( 2n_r,\bbC)\quad \text{(cas {\bf C} et {\bf D})},\end{cases}\]
avec $n_r \in \bbN$ (le cas $n_r=0$ est une possibilité parfaitement loisible, auquel cas $\caL_r=\{1 \}$ 
   que l'on peut omettre dans (\ref{LeviG})). On a d'autre part  $\displaystyle \sum_{i=1}^r n_i=n$.

Soit $(\mathbf{Q},\mathbf{L})$ représentant une classe dans $\caS_\caQ$ ({\sl cf.} section \ref{AJpackSec}).   On laisse 
au lecteur vérifier que l'on a  alors un isomorphisme
   \[ L\simeq L_1\times \cdots \times L_r\]
   où  pour tout $i=1, \ldots, r-1$, $L_i=U(b_i,c_i)$,   avec $b_i+c_i=n_i$  et suivant les cas :
  \begin{itemize}
   \item[\bf A.]   $L_r=\Sp(2n_r,\bbR)$, 
         \item[\bf B.]  $L_r=\SO(b_r,c_r)$ avec $b_r+c_r=2n_r+1$, 
      $\displaystyle \sum_{i=1}^{r-1}2b_i+b_r=n$, $\displaystyle \sum_{i=1}^{r-1}2c_i+c_r=n+1$,
       \item[\bf C.] $L_r=\SO(b_r,c_r)$,   avec $b_r+c_r=2n_r$,  
         $\displaystyle \sum_{i=1}^{r-1}2b_i+b_r=n$, $\displaystyle \sum_{i=1}^{r-1}2c_i+c_r=n$,
   \item[\bf D.] $L_r=\SO(b_r,c_r)$,   avec $b_r+c_r=2n_r$,    $\displaystyle \sum_{i=1}^{r-1}2b_i+b_r=n-1$, 
   $\displaystyle \sum_{i=1}^{r-1}2c_i+c_r=n+1$.
\end{itemize}

 \medskip

Le paramètre $\psi$ se décompose comme en (\ref{opsii}) en somme de paramètres \og presque irréductibles\fg\, 
 $\psi=\oplus_{i=1, \ldots ,r}\psi_i$, où pour $i=1,\ldots ,r-1$,  
\[ \psi_i \simeq V\left(-\frac{p_i}{2},\frac{p_i}{2}\right)\otimes R_{n_i} :   W_\bbR\times \SL_2(\bbC) \longrightarrow  \widehat G_{N_i}, \]
 \[ \Pi_{\psi_i}=\Speh(\delta_i,n_i) \] 
avec $\delta_i=\delta\left(-\frac{p_i}{2}, \frac{p_i}{2} \right)$, 
 $p_i+n_i$ impair dans les cas {\bf A}, {\bf C}, {\bf D}, pair dans le cas {\bf B}.
\medskip 
Pour $i=r$, on a suivant les cas: 

   {\bf A.} $N_r=2n_r+1$,   $\psi_{r}=\epsilon \otimes R_{N_r}$, $\epsilon= \Triv$ ou $\mathbf{sgn}$,
 
   {\bf B.} $N_r=2n_r$,  $\psi_{r}=\epsilon \otimes R_{N_r}$,  $\epsilon= \Triv$ ou $\mathbf{sgn}$, 
  
  {\bf C} et {\bf D}. $N_r=2n_r$,  $\psi_{r}=\epsilon \otimes R_{N_r-1}\oplus \epsilon $,
         $\epsilon= \Triv$ ou $\mathbf{sgn}$ ou
        $\psi_{r}=\epsilon_{1} \otimes R_{N_r-1}\oplus \epsilon_{2} $,  
   $(\epsilon_1, \epsilon_2)=( \Triv, \mathbf{sgn})$ ou $(\sgn,\Triv)$. 

Comme on le voit, les $\psi_i$ sont irréductibles sauf $\psi_r$ dans les cas {\bf C} et {\bf D}.
Nous dirons que $\psi=\oplus_{i=1, \ldots ,r}\psi_i$ est une décomposition en paramètres élémentaires.

\begin{rmq}
      Dans le cas {\bf A},  les $\psi_i$, $i=1,\ldots,r$,  sont à valeurs dans $\SO(2n_i,\bbC)$ si $n_i$ est pair,
       mais seulement dans  $\Or(2n_i,\bbC)$ si $n_i$ est impair. On  a donc $\psi_r=\epsilon \otimes R_{N_r}$ avec 
      $\epsilon=\Triv$ si le nombre des $i$ dans $\{1, \ldots, r-1\}$ avec $n_i$  impair est pair,
       et $\epsilon=\sgn$ s'il est impair.
      
      Dans le cas {\bf B}, tous les $\psi_i$ se factorisent par le groupe $\Sp(2n_i,\bbC)$. Pour $i=r$, si 
      $\epsilon=\Triv$, $\psi_r$ se factorise par le paramètre de la représentation triviale de $\SO(b_r,c_r)$, si 
      $\epsilon=\sgn$, $\psi_r$ se factorise par le paramètre   du caractère valant $-1$ sur la composante 
        connexe non triviale de  $\SO(b_r,c_r)$.
        
        Dans les cas {\bf C} et  {\bf D},  les $\psi_i$ se factorisent par le groupe $\SO(2n_i,\bbC)$ si $n_i$ est pair
        et par le groupe $\Or(2n_i,\bbC)$ si $n_i$ est impair. 
       Pour $i=r$, dans le cas {\bf C},  si le nombre des $i$ dans $\{1, \ldots, r-1\}$ avec $n_i$  impair est pair,
        alors  $\psi_{r}=\epsilon \otimes R_{N_r-1}\oplus \epsilon $, $\epsilon= \{\Triv,\mathbf{sgn}\}$   et si ce nombre est impair, alors 
           $\psi_{r}=\epsilon_{1} \otimes R_{N_r-1}\oplus \epsilon_{2} $,  avec
   $(\epsilon_1, \epsilon_2)=( \Triv, \mathbf{sgn})$ ou $(\sgn,\Triv)$. 
 Pour $i=r$, dans le cas {\bf D}, c'est l'opposé du cas {\bf C}.
         
         Le paramètre  $\psi_r=\epsilon \otimes R_{N_r-1}\oplus \epsilon$ se factorise par le paramètre de la représentation triviale 
         de $\SO(b_r,c_r)$ si $\epsilon=\Triv$, et si 
      $\epsilon=\sgn$, $\psi_r$ se factorise par le paramètre   du caractère valant $-1$ sur la composante 
        connexe non triviale de  $\SO(b_r,c_r)$.

         Le paramètre   $\psi_{r}=\epsilon_{1} \otimes R_{N_r-1}\oplus \epsilon_{2} $,  
 avec $(\epsilon_1, \epsilon_2)=( \Triv, \mathbf{sgn})$, $\psi_r$ se factorise par le paramètre 
        de la représentation triviale de $\SO(b_r,c_r)$, et si $(\epsilon_1, \epsilon_2)=(\sgn,\Triv)$, 
        $\psi_r$ se factorise par le paramètre   du caractère valant $-1$ sur la composante 
        connexe non triviale de  $\SO(b_r,c_r)$.
\end{rmq}        

\medskip

La représentation irréductible  $\Pi_{\psi_r}$ de $G_{N_r}$
associée au paramètre $\psi_r$ est 
\begin{equation}
\Pi_ {\psi_r} =  \begin{cases} \varepsilon_{N_r} \text{  si  } \psi_{G,r}=\epsilon \otimes R_{N_r}, 
\quad \epsilon\in \{\Triv,\sgn\}, \quad \text{(cas {\bf A},  {\bf B})}\\
 \varepsilon_{1, N_r-1}\times \varepsilon_2  \text{  si  } \psi_{G,r}=\epsilon_1 \otimes R_{N_r-1}\oplus \epsilon_2
 \quad \epsilon_1, \epsilon_2\in \{\Triv,\sgn\}  \quad \text{(cas {\bf C},  {\bf D}) }.\end{cases}
 \end{equation}

  \section{Démonstration du résultat principal}\label{principal}
  
\subsection{Enoncé}  
  Soit $\mathbf{G}$ l'un des groupes classiques de la section \ref{grcla} et $\St_G :\; {}^L G \rightarrow \widehat G_N$
la représentation standard de son $L$-groupe. Soit $\psi_G : \; W_\bbR \times \SL_2(\bbC) \rightarrow {}^L G $ 
un paramètre  d'Adams-Johnson ({\sl cf.} section  \ref{AJcla}) pour $\mathbf{G}$ et $\Pi^{AJ}_{\psi_G}$ 
le paquet d'Adams-Johnson
qui lui est associé ({\sl cf.} Définition \ref{AJpack}). Notons $\psi=\St_G \circ \psi_G$ : c'est un  paramètre d'Arthur
de $G_N$, auquel est  associé une représentation $\Pi_\psi=\Pi_{\phi_\psi}$  de ce groupe.

Notre but est de démontrer l'égalité de distributions sur $\widetilde G_N$ : 
\begin{equation}\label{butfinal}
\Tr_{\theta_N} (\Pi_\psi) =\Trans_G^{G_N} (\Theta_{ \Pi_{\psi_G}^{AJ}}^{st}) 
\end{equation}
où $\Theta_{ \Pi_{\psi_G}^{AJ}}^{st}$ est la distribution stable sur $G$ définie en (\ref{StablePiAJ}).

D'après (\ref{IdStableAJ}), où les notations sont expliquées, on a une égalité de distributions sur $G$
\begin{equation}
\Theta_{ \Pi_{\psi_G}^{AJ}}^{st}=  (-1)^{q(L^*)} \sum_{\phi_G \in \Phi(\psi_G)} (-1)^{\ell_V(\phi_G)}
  \;  \Theta_ {\widetilde \Pi_ {\phi_G}}
\end{equation}
Tous les termes de cette formule sont des distributions stables, et la fonction longueur $\ell_V$ est celle de Vogan.
On en déduit que 
\begin{equation}\label{11}
\Trans_G^{G_N} ( \Theta_{ \Pi_{\psi_G}^{AJ}}^{st}  ) =  (-1)^{q(L^*)}
\sum_{\phi_G \in \Phi(\psi_G)} (-1)^{\ell_V(\phi_G)}  \;  \Trans_G^{G_N} ( \Theta_ {\widetilde \Pi_ {\phi_G}}  ).
\end{equation}
Or, d'après la proposition \ref{TransPsPaq}, on a pour tout $\phi_G \in \Phi(\psi_G)$
\begin{equation}\label{22}
\Trans_G^{G_N} ( \Theta_ {\widetilde \Pi_ {\phi_G}} ) =\Tr_{\theta_N} (\widetilde \Pi_{\St_G \circ \phi_G}) ,
\end{equation}
 où $\widetilde\Pi_{\St_G \circ \phi_G}$ est la représentation standard de $G_N$ associée par Langlands au paramètre 
  de Langlands $\phi=\St_G \circ \phi_G$.
Ainsi l'identité à établir est équivalente à :
\begin{equation}\label{trator}
\Tr_{\theta_N} (\Pi_\psi) =   (-1)^{q(L^*)}  \sum_{\phi_G \in \Phi(\psi_G)} (-1)^{\ell_V(\phi_G)}  \; \Tr_{\theta_N} 
(\widetilde \Pi_{\St_G \circ \phi_G}). 
 \end{equation}

Dans les deux prochains paragraphes, nous allons établir ce résultat lorsque  $\psi_G$ est un paramètre élémentaire.
  
  \subsection{Identités endoscopiques pour les paramètres élémentaires : représentations de dimension $1$}  
  \label{dim1}      On suppose dans ce paragraphe que dans la décomposition $\psi=\oplus_{i=1}^r\psi_{i}$, on a $r=1$,  $n=n_r>0$  et 
   $\psi=\psi_{r} $. Rappelons que nous sommes alors dans un des cas suivants :  
    
  {\bf A} :   $\psi=\Triv\otimes R_{N}$,   $N=2n+1$; 
   {\bf B} :    $\psi=\epsilon \otimes R_{N}$,  $\epsilon= \Triv$ ou $\mathbf{sgn}$, $N=2n$;  
   {\bf C} : $\psi=\epsilon \otimes R_{N-1}\oplus \epsilon $,
         $\epsilon= \Triv$ ou $\mathbf{sgn}$,   $N=2n$; 
  {\bf D} :  $\psi=\epsilon_{1} \otimes R_{N-1}\oplus \epsilon_{2} $,   
   $(\epsilon_1,\epsilon_2)= (\Triv,\mathbf{sgn}) $ ou $(\sgn,\Triv)$,    $N=2n$. 

\medskip
   
   On a alors  : 
   
   \begin{lemme}
  Le paquet d'Arthur $\Pi_{\psi_G}$ est un singleton, plus précisement $\Pi_{\psi_G}=\Pi_{\phi_{\psi_G}}$.
  \end{lemme}
  
\dem  On calcule explicitement le paramètre de Langlands $ \phi_{\psi_G}$. Dans le cas {\bf A}, 
c'est celui de la représentation triviale. Dans autres cas, c'est celui  de la   représentation triviale 
ou du  caractère du groupe spécial orthogonal   $G$  trival sur la composante neutre et valant $-1$ sur l'autre. 
Il suffit dans tous les cas de traiter le cas de la représentation triviale, 
les paquets étant préservés par tensorisation par 
un caractère.  Le paquet $\Pi_{\psi_G}$ contient donc la représentation triviale, et toutes les autres représentations de 
ce paquet ont même caractère infinitésimal que celle-ci. Soit $\pi \in  \Pi_{\psi_G}$. 
 On globalise la situation sur $\bbQ$ en considérant le paramètre
global donné par un $\SL_2$ régulier dans $\widehat G$. En toute place, la représentation triviale est dans le 
paquet local associé à ce paramètre. Considérons la représentation globale dont la composante locale
 à toutes les    places finies est la représentation triviale sauf à la place archimédienne où elle est égale à  $\pi$.
D'après le théorème 1.5.2 de \cite{Art13}, cette représentation globale apparaît dans le spectre discret.
Elle se réalise donc dans un ensemble de fonctions de carré intégrable, invariantes par translation à gauche sous 
$\mathbf G(\bbQ)$.  Mais une fonction qui est de plus invariante par translation à gauche par $\mathbf G(\bbQ_p)$, 
  pour tout nombre premier $p$,  est une fonction constante, et $\pi$ est nécessairement la représentation triviale. 
  \qed
  
\medskip  
  
Comme le paquet défini par Adams et Johnson est lui aussi dans ce cas réduit au singleton 
$\Pi_{\psi_G}^{AJ}=\Pi_{\phi_{\psi_G}}$, on en déduit que $\Pi_{\psi_G}^{AJ} = \Pi_{\psi_G}$. 
Par définition de $\Pi_{\psi_G}$ dans \cite{Art13}, on a donc 
  \begin{equation}\label{IdEndParEl1}
\mathrm{Trans}_G^{G_N}( \Theta_{\Pi_{\psi_G}^{AJ}}^{st})=\Tr_{\theta_N}(\Pi_\psi).
\end{equation}  
  L'identité (\ref{butfinal}) s'obtient donc directement dans ce cas. 
  Nous avons vu qu'elle est équivalente à (\ref{trator}), c'est-à-dire que nous avons
 \begin{equation}\label{iDe}
\Tr_{\theta_N}( \Pi_\psi)= (-1)^{q(G)}  \sum_{\phi_G \in \Phi(\psi_G) } (-1)^{\ell_V(\phi_G)}  \Tr_{\theta_N} (\widetilde\Pi_{\Std_G\circ \phi_G}) . 
  \end{equation}
 Dans ce cas, remarquons qu'avec les notations de la section \ref{ResJohPaq}, le sous-groupe de Levi $\caL$ attaché à 
$\psi_G$   est $\widehat G$, et ainsi $\caS_\caQ$ est bien  le  singleton $(\mathbf Q,\mathbf L)=(\mathbf G,\mathbf G)$.
Ainsi dans ce cas $L=L^*=G$ car $G$ est quasi-déployé, et 
l'ensemble des paramètres de Langlands $\Phi (\psi_G)$ 	est donc l'ensemble des paramètres des représentations
 standard intervenant dans la résolution de Johnson de  la représentation dimension $1$ de $G$ considérée.

  \subsection{Paramètres de Speh}\label{ParSpehE}
Nous nous replaçons dans le m\^eme contexte que la section précédente, mais l'on 
    suppose dans ce paragraphe que $r=2$ et $n_r=0$. Autrement dit $\psi$ est irréductible, 
    et plus précisément, $\psi$ est le paramètre
   d'Arthur d'une représentation de Speh $\Pi_\psi=\Speh\left( \delta\left( -\frac{p}{2}, \frac{p}{2} \right),n \right)$.
  On a par hypothèse $p>n-1$ ou $p\geq n-1$ dans les cas {\bf C} et {\bf D}. 
   Si $p+n$ est impair, on est dans un des cas  {\bf C}
  ou {\bf D} (selon que $n$ est pair ou impair respectivement)  et si  $p+n$ est pair, on est dans le cas {\bf B}.

 Le membre de gauche de l'identité  (\ref{trator}),  est d'après (\ref{FormTTSpeh}),   
   \begin{equation}\label{TrSpeh4}
 \Tr_{\theta_N}(\Pi_\psi) = \sum_{s \in \frI_n}
 (-1)^{\ell_\theta(s)} \; \Tr_{\theta_N}(X(s)).
\end{equation}
  Il s'agit donc d'établir que 
   \begin{equation}\label{jesaisplus} \sum_{s \in \frI_n}  (-1)^{\ell_\theta(s)} \; \Tr_{\theta_N}(X(s))= (-1)^{q(L^*)}
\sum_{\phi_G \in \Phi(\psi_G) } (-1)^{\ell_V(\phi_G)} \Tr_{\theta_N}(\widetilde \Pi_{\Std_G\circ \phi_G}).\end{equation}  
Il suffit pour cela d'établir qu'il existe une bijection $s \in \frI_n  \mapsto \phi_G \in \Phi(\psi_G)$
telle que 
\begin{equation}\label{longw0} X(s)=\widetilde \Pi_{\phi}, \quad \phi=\Std_G\circ \phi_G, \text{ et } q(L^*)+\ell_V(\phi_G)=\ell_\theta(s) \mod 2.\end{equation} 

Rappelons comment sont obtenus l'ensemble des paramètres $\Phi(\psi_G)$ en reprenant les notations de la section
\ref{ResJohPaq} pour le paramètre d'Adams-Johnson $\psi_G$.
L'ensemble  $\Phi(\psi_G)$ est alors l'ensemble des paramètres de Langlands 
des représentations standard apparaissant dans les résolutions des $\pi_{\psi_G, \mathbf Q,\mathbf L}$, 
$(\mathbf Q,\mathbf L)\in \caS_\caQ$ et $L\simeq \U(p,q)$, $p+q=n$. 
Celles-ci sont obtenus par induction cohomologique à partir des représentations standard apparaissant
 dans la résolution des caractères unitaires $\lambda$ de $L$ tel que 
 $\pi_{\psi_G, \mathbf Q,\mathbf L}=A_\frqqq(\lambda)$. D'après la section \ref{BBUpq}, 
ces représentations standard sont paramétrées par $\frI_n^{p,q,\pm}$, et leur paramètres de Langlands par 
$\frI_n$. Plus précisément, si $\bar s \in \frI^{p,q,\pm}$ est le paramètre
 d'une représentation standard 
apparaissant dans la résolution du caractère $\lambda=\lambda_{\psi_G,\mathbf{Q}, \mathbf{L}}$ de $\U(p,q)$
(un tel caractère est caractérisé par son caractère infinitésimal car $\U(p,q)$ est connexe, et donc facile à déterminer), 
et si $\phi_s:\, W_\bbR \rightarrow {}^LL$ est son paramètre de Langlands (qui ne dépend que de l'involution
$s\in \frI_n$ sous-jacente), alors $\Pi_{\Std_G\circ i_{L,G} \circ \phi_s}=X(w_0s)$. 
Notons l'apparition de  $w_0$, l'élément le plus long de $\frS_n$.   Ainsi 
$s\in \frI_n \mapsto  i_{L,G} \circ \phi_{w_0s}$ est la bijection voulue de $\frI_n$ dans $\Phi(\psi_G)$.

Vérifions maintenant que les longueurs coïncident. Nous donnons deux arguments, l'un utilisant
l'égalité (\ref{complinv}) laissée en exercice au lecteur,
 et l'autre le changement de base. Soit $X(\gamma)$ une représentation standard apparaissant dans la 
résolution de $\pi_{\psi_G, \mathbf Q,\mathbf L}$, de paramètre de Langlands $\phi_G$,  et $X_L(\gamma')$ 
la représentation standard  apparaissant  dans la résolution du caractère unitaire $\lambda$ de $L$
qui lui correspond (ici $\gamma$ et $\gamma'$ sont des paramètres pour les représentations
comme dans \cite{Vgreen}  ou \cite{Joh1}, nous n'explicitons pas leur
forme exacte). On a alors $\ell_V(\gamma')=\ell_V(\gamma)$, les longueurs étant les longueurs 
de Vogan dans les groupes $L$ et $G$ respectivement, ceci est une conséquence facile de la définition de la longueur de Vogan.
 Soit $Q_{\gamma'}$ l'orbite 
 attachée à $\gamma'$ dans la variété des  drapeaux $\scrB_L$ de $\mathbf L$ dans la paramétrisation
 de Beilinson-Bernstein des représentations des groupes unitaires ({\sl cf.} section \ref{BBUpq}).
  La longueur de Vogan et la dimension de l'orbite qui lui correspond sont reliées par (\cite{VIC3}, p. 396)
\begin{equation} \label{lVdimQ}\ell_V(\gamma')=\dim Q_{\gamma'}-\dim \caB_{L}+q(L) \end{equation}
Or nous avons vu dans la section \ref{BBUpq} que la dimension des orbites $Q_{\gamma'}$  dans $\caB_{L}$
est donnée par la fonction longueur définie sur $\frI_n^{p,q,\pm}$. Si $s\in \frI_n$ est
l'involution qui donne  le  paramètre de Langlands $\phi_G$, on a en utilisant (\ref{complinv})
{\small \begin{align*}
\ell_V(\phi_G)=\ell_V(\gamma)=\ell_V(\gamma')=\dim Q_{\gamma'}-\dim \caB_{L}+q(L)
 = \ell_\theta(s)+\frac{1}{2}(p(p-1)+q(q-1))-\dim \caB_{L}+q(L)
\end{align*}}
Or $\dim \caB_{L}=\frac{n(n-1)}{2}=\frac{(p+q)(p+q-1)}{2}$, et $q(L)=\frac{1}{2}( (p+q)^2-p^2-q^2 )=pq$, 
d'où finalement 
\begin{equation} \label{egalitelongueur}  \ell_V(\phi_G)=\ell_\theta(s)
\end{equation}
Il correspond à $\phi_G$ l'élément $w_0s$ par la bijection ci-dessus, et on a $\ell_\theta(w_0s)=
\ell_\theta(w_0)+\ell_\theta(s)\mod 2$. Pour montrer la relation sur la longueur dans (\ref{longw0}), il suffit donc de vérifier que 
$q(L^*)=\ell_\theta(w_0)$. Or, la relation (\ref{egalitelongueur}) nous dit que $\ell_\theta(w_0)$ est la longueur de Vogan maximale pour 
un paramètre de Langlands  des groupes unitaires $\U(p,q)$, et (\ref{lVdimQ}) nous dit que l'on a un tel paramètre lorsque l'orbite correspondante
 est ouverte et $q(L)$ est maximal. Or $q(L)$ est maximal pour la forme quasi-déployée. C'est une propriété générale, qui se voit bien ici pour 
$q(\U(p,q))=pq$.

Pour le second argument, on part de la résolution de Johnson de la représentation triviale 
du groupe quasi-déployé $\U(b,c)$, où $b+c=n$. Les représentations standard 
de $\U(b,c)$ de même caractère infinitésimal que la représentation triviale sont indexées par $\frI_n^{b,c, \pm}$ comme dans la section
\ref{BBUpq}. Appelons ici $Y(\bar s)$ le module standard indexé par  par $\bar s\in \frI_n^{b,c, \pm}$
 On sait que la résolution de Johnson donne une identité de représentations virtuelles de la forme
\[[\Triv_{\U(b,c)}]=  \sum_{\bar s\in \frI_n^{b,c, \pm} }  a(\bar s)\,  [Y(\bar s)], \]
où $a(\bar s)$ est un signe. Mais le changement de base vers $\GL_n(\bbC)$ envoie la représentation triviale
sur la réprésentations triviale, et le  module   standard  $Y(\bar s)$ sur  le modules standard $\bar \theta_n$-stable de $\GL_n(\bbC)$, 
noté $X(\rho,-s  \rho)$ dans la section \ref{cartorC}, où $s$ est l'involution dans $\frS_n$ donnée par $\bar s$ en oubliant les signes.
 Ceci découle de \cite{CloBC} (voir aussi la 
section \ref{GrUnit}) et il en résulte que les signes
$a(\bar s)$ sont donnés par $(-1)^{\ell_{\theta}(w_0s)}$. Comme dans \cite{Joh1}, 
on tensorise par le caractère unitaire idoïne de $\U(b,c)$, et on induit cohomologiquement de $\U(b,c)$
à $G$ pour obtenir une formule des caractères pour un élément du paquet d'Adams-Johnson $\Pi_{\psi_G}^{AJ}$.
La formule pour la distribution stable $\Theta_{\Pi_{\psi_G}^{AJ}}^{st}$ est obtenu en \og stabilisant \fg  \, 
cette formule du caractère, comme cela a été mentionné dans la remarque \ref{cryptic}.
Ceci se fait en y ajoutant des formules des caractères obtenues de la même manière pour les formes intérieures
$\U(b',c')$ de $\U(b,c)$, affectées d'un signe que nous n'avons pas besoin de calculer. Notons $Z(\bar s)$ la représentation standard de 
$G$ cohomologiquement induite de $Y(\bar s)$. La remarque \ref{cryptic} montre que les pseudo-paquets
contenant les représentations standard $Y(\bar s)$ sont indexés par $\frI_n$ et  notons  $\tilde Z(s)$ la somme des représentations standard
dans le pseudo-paquet contenant $Z(\bar s)$.
On obtient alors
\[ \Theta_{\Pi_{\psi_G}^{AJ}}^{st}= \sum_{s\in \frI_n} (-1)^{\ell_{\theta}(w_0s)} \Theta_{\widetilde Z(s)}
= (-1)^{q(L^*)}\sum_{s\in \frI_n} (-1)^{\ell_{\theta}(s)} \Theta_{\widetilde Z(s)}. \]
Après  transfert endoscopique $\Trans_G^{\widetilde G_N}$, ceci est égal au second membre de (\ref{jesaisplus}), et par identification,
 fournit une bijection entre $\frI_n$ et  $\Phi(\psi_G)$ avec les propriétés requises.
Ceci achève d'établir (\ref{jesaisplus}).\qed

\subsection{Réduction au cas élémentaire} \label{reduc}

Nous revenons au cas général. Soit   $\psi=\oplus_{i=1}^r \psi_i$ 
la décomposition de $\psi=\Std_G\circ \psi_G$  en paramètres élémentaires comme dans la section \ref{AJcla}.
Nous allons établir (\ref{trator}) par récurrence sur $r$. Le cas $r=1$ a été établi dans la section précédente.
Posons donc $\psi'=\oplus_{i=2}^r \psi_i$, de sorte que $\psi=\psi_1\oplus \psi'$. Posons 
aussi $N'=\sum_{i=2}^r N_i$.

On peut, quitte à permuter les indices, supposer que $p_1$ est le plus grand des paramètres
de série discrète apparaissant dans la décomposition, de sorte que l'inégalité (\ref{hyplemmegen})
est satisfaite. On peut donc utiliser les résultats de la section \ref{IndTor} et en particulier de \ref{InductionTord}.

Remarquons que dans le cas {\bf A}, le paramètre $\psi'$ n'est pas forcément obtenu 
à partir d'un paramètre d'Arthur pour un groupe symplectique plus petit;  c'est le cas si
$n_1$ est impair, car alors $\psi'$ n'est à valeur dans $\SO(N',\bbC)$ mais seulement dans 
$\Or(N',\bbC)$.  Ceci montre la nécessité de  tensoriser $\psi'$ par le caractère signe de $W_\bbR$, c'est-à-dire poser :
\[\psi''=\begin{cases} \psi' \text{ si } n_1 \text{ est pair }\\
\psi'\otimes  \sgn_{W_\bbR} \text{ si } n_1 \text{ est impair }\end{cases} .  \]
Cette torsion va aussi apparaître dans les autres  cas. Avec ces notations, $\psi''$ est de la forme 
 $\psi''=\Std_{G'}\circ \psi_{G'}$ pour un groupe classique $G'$ de même type que $G$ (sauf  les cas {\bf C} et {\bf D} qui sont échangés si 
si $n_1$ est impair).  Dans tous les cas, posons 
 \begin{equation}\label{etor} \psi''= \epsilon^1_{tor}\otimes \psi' . \end{equation}
 
 Notons $G_1$ le groupe classique tel que $\psi_1$ se factorise comme $\Std_{G_1}\circ \psi_{G_1}$.
 Nous avons établi dans le  paragraphe \ref{ParSpehE} que  l'on a :
\begin{equation}\label{tratori}
\Tr_{\theta_{N_1}} (\Pi_{\psi_1}) = (-1)^{q(L_1^*)} \sum_{\phi_{G_1} \in \Phi(\psi_{G_1})} (-1)^{ \ell_V(\phi_{G_1})}  \;
 \Tr_{\theta_{N_1}} (\widetilde \Pi_{\Std_{G_1}\circ \phi_{G_1}}) 
 \end{equation}
et  l'on suppose par hypothèse de récurrence que l'on a une formule du même type pour $\psi''$, c'est-à-dire : 
\begin{equation}\label{tratorprimeprime}
 \Tr_{\theta_{N'}} (\Pi_{\psi''}) =(-1)^{q({L'}^*)} \sum_{\phi_{G'} \in \Phi(\psi_{G'})} (-1)^{ \ell_V(\phi_{G'})}  \;
 \Tr_{\theta_{N'}} (\widetilde \Pi_{\Std_{G'}\circ \phi_{G'}}) .
  \end{equation}
   
  En tensorisant  par le caractère $\epsilon^1_{tor}$, on obtient 
  \begin{equation}\label{tratorprime}
 \Tr_{\theta_{N'}} (\Pi_{\psi'}) = (-1)^{q({L'}^*)} \sum_{\phi_{G'} \in \Phi(\psi_{G'})} (-1)^{ \ell_V(\phi_{G'})}  \;
 \Tr_{\theta_{N'}} (\widetilde \Pi_{(\Std_{G'}\circ \phi_{G'})\otimes \epsilon^1_{tor}}).
  \end{equation}
 
 Dans les équations ci-dessus, $L_1^*$ et ${L'}^*$ sont bien sûr les analogues de $L^*$ pour les groupes $G_1$ et $G'$ respectivement.
  
Comme dans la section \ref{InductionTord}, nous notons $M$ le sous-groupe de Levi standard de $G_N$ isomorphe
à $G_{N_1} \times G_{N'}$.
Le produit tensoriel  $\Pi_{\psi_1}\otimes \Pi_{\psi'}$ est alors une représentation de $M$ et l'on a avec les notations de cette section
\begin{align}
&\tr_{\theta_M} \left( \Pi_{\psi_1}  \otimes  \Pi_{\psi'} \right)= (-1)^{q({L_1}^*)} (-1)^{q({L'}^*)} \\
\nonumber &\times  \sum_{(\phi_{G_1},\phi_{G'}) \in   \Phi(\psi_{G_1})\times \Phi(\psi_{G'})  }
 (-1)^{ \ell_V(\phi_{G_1})+\ell_V(\phi_{G'})} \; 
 \tr_{\theta_M} \left( \widetilde  \Pi_{\Std_{G_1}\circ \phi_{G_1}}    \otimes \widetilde  \Pi_{(\Std_{G'}\circ \phi_{G'})\otimes  \epsilon^1_{tor}}   \right).\end{align}
 
On obtient alors par le lemme \ref{colle}, et le fait que $L_*=L_1^\times *{L'}^* $
  \begin{equation}
 \Tr_{\theta_N} (\Pi_\psi)=(-1)^{q(L^*)}
 \sum_{(\phi_{G_1},\phi_{G'}) \in   \Phi(\psi_{G_1})\times \Phi(\psi_{G'})  } (-1)^{ \ell_V(\phi_{G_1})+\ell_V(\phi_{G'})} \;   
  \Tr_{\theta_N}( \widetilde \Pi_{\Std_{G_1}\circ  \phi_{G_1}\oplus   (\Std_{G'}\circ \phi_{G'})\otimes  \epsilon^1_{tor} })
 \end{equation}
 
 L'égalité  à  démontrer (\ref{trator}), 
 est donc équivalente à 
  \begin{align}\label{trator5}&
 \sum_{\phi_G \in \Phi(\psi_G)} (-1)^{\ell_V(\phi_G)}  \; \Tr_{\theta_N} 
(\widetilde \Pi_{\St_G \circ \phi_G})\\
\nonumber &=\sum_{(\phi_{G_1},\phi_{G'}) \in   \Phi(\psi_{G_1})
\times \Phi(\psi_{G'})  } (-1)^{ \ell_V(\phi_{G_1})+\ell_V(\phi_{G'})} \;   
  \Tr_{\theta_N}( \widetilde \Pi_{ \Std_{G_1}\circ \phi_{G_1}\oplus  (\Std_{G'}\circ  \phi_{G'})\otimes  \epsilon^1_{tor} }).
 \end{align}

  Il reste donc à montrer  que l'on a  une bijection 
 \begin{equation}\label{bijparetlong}       \Phi(\psi_{G_1})\times \Phi(\psi_{G'})\longrightarrow 
     \Phi(\psi_G), \quad  (\phi_{G_1},\phi_{G'})   \mapsto \phi_{G}
 \end{equation} 
 telle que  $\Std_G\circ \phi_G=\phi_{G_1}\oplus   \phi_{G'}\otimes  \epsilon^1_{tor} $ et     $\ell_V(\phi_G)=\ell_V(\phi_{G_1})+\ell_V(\phi_{G'})$.

  Rappelons les définitions,  données dans la section \ref{ResJohPaq}. A $\psi_G$ est associé 
  un sous-groupe parabolique $\caQ=\caL\caU$ de $\widehat G$, puis un ensemble fini
  $\caS_\caQ$ de couples $(\mathbf Q,\mathbf L)$, où $\mathbf Q$ est un sous-groupe
  parabolique de $\mathbf G$ et $\mathbf L$ un facteur de Levi  de $\mathbf Q$ défini sur $\bbR$.
  A chaque $(\mathbf Q,\mathbf L)\in \caS_\caQ$ est attachée en (\ref{pipsiQL}) une représentation 
  $\pi_{\psi_G,\mathbf Q,\mathbf L}$, induite cohomologique d'un certain caractère 
  $\lambda_{\psi_G,\mathbf Q,\mathbf L}$ de $L$.
   Alors $\Phi(\psi_G)$ est l'ensemble des paramètres de 
  Langlands des modules standard intervenant dans la résolution de Johnson
  d'un des  modules $\pi_{\psi_G,\mathbf Q, \mathbf L}$. Mais on peut dire mieux, d'après la remarque \ref{cryptic}.
  En effet, $G$ étant quasi-déployé, il existe $(\mathbf{Q}_*,\mathbf{L}_*)\in \caS_\caQ$ tel que $\mathbf{L}_*$ 
  est quasi-déployé,   et   $\Phi(\psi_G)$ est l'ensemble des paramètres de 
  Langlands des modules standard intervenant dans la résolution de Johnson de  $\pi_{\psi_G,\mathbf{Q}_*, \mathbf{L}_*}$. 
  
 Cet ensemble est en bijection   avec les paramètres de    Langlands  des modules standard 
 intervenant dans la résolution de  Johnson du caractère  $ \lambda_*=\lambda_{\psi_G,\mathbf Q_*,\mathbf L_*}$   de $L_*$, 
 par  $\phi_L\mapsto \phi_G= \iota_{L,G}\circ \phi_L$, où $\iota_{L,G}$ est le plongement de $L$-groupe (\ref{iotaLG}), 
 la correspondance 
 entre modules standard étant obtenue par le foncteur d'induction cohomologique 
 $\caR_{\frqqq_*}^i$  ({\sl cf. }(\ref{pipsiQL})) qui  préserve les longueurs de Vogan des paramètres.
  Pour les groupes qui nous occupent, les ensembles $\caS_\caQ$ ont été décrit dans la section 
  \ref{AJcla}. Le sous-groupe de Levi $L_*$ est de la forme
  \begin{equation}\label{decLevi} L_*\simeq L_{1} \times \cdots \times L_{r}\end{equation}
où pour $i=1,\ldots r-1$,  $ L_{i}$ est un groupe unitaire quasi-déployé, et 
$ L_{r}$ est un groupe symplectique, ou un groupe spécial orthogonal quasi-déployé. 
On pose alors $L'=L_2\times \cdots \times L_r$, c'est un sous-groupe de Levi de $G'$, et $L_1$ est un sous-groupe de Levi 
de $G_1$.

Il est clair que $ \lambda_*$ est un produit tensoriel  
de deux  caractères $\lambda_1$ et $\lambda'$  de $ L_{1}$ et $L'$ respectivement, et  que  l'ensemble des modules standard
 intervenant dans la résolution  de $\lambda$ est en bijection avec le produit tensoriel de deux modules standard 
 intervenant respectivement  dans les résolutions  de $\lambda_1$, et $\lambda'$,  et que la longueur de Vogan d'un module standard
 est la somme des longueur des deux modules standard  qui lui correspondent.

La description de  $\Phi(\psi_{G_1})$ et $\Phi(\psi_{G'})$ est analogue. Un élément $\phi_{G_1}$ de 
$\Phi(\psi_{G_1})$ s'écrit  $\phi_{G_1}= \iota_{L_1,G_1}\circ \phi_{L_1}$ où $\phi_{L_1}$ est le paramètre de Langlands
d'une représentation standard intervenant dans la résolution de $\lambda_1$, et 
un élément $\phi_{G'}$ de 
$\Phi(\psi_{G'})$ s'écrit  $\phi_{G'}= \iota_{L',G'}\circ \phi_{L'}$ où $\phi_{L'}$ est le paramètre de Langlands
d'une représentation standard intervenant dans la résolution de $\lambda'$. En effet, 
nous affirmons que si 
$\phi_{G'}$ est la paramètre de Langlands d'une représentation standard de $G'$ apparaissant dans la résolution de l'induite cohomologique
de $L'$ à $G'$ du caractère  $\lambda'$,  alors 
$\phi=\Std_{G_1}\circ \phi_{G_1}\oplus (\Std_{G'}\circ \phi_{G'})\otimes \epsilon_{tor}^1$ 
se factorise en $\Std_G\circ \phi_G$, et que l'application $(\phi_{G_1},\phi_{G'} )\mapsto \phi_G$
est la bijection \ref{bijparetlong}  voulue.

Pour vérifier l'assertion ci-dessus, il faut savoir déterminer  le paramètre de Langlands d'une représentation standard
obtenue par   induction cohomologique   d'une représentation standard. Ceci se fait grâce aux théorèmes 
\og d'indépendance de polarisation \fg\, de \cite{KV}, Chapter 11, mais une forme commode de ceux-ci
pour le  calcul que nous avons à effectuer  est  écrite dans \cite{Matu}, Thm. 2.2.3.

 Ceci termine la démonstration du résultat principal. \qed

\begin{rmq}
L'identité à démontrer (\ref{trator}) est une identité de trace tordue pour $\widetilde G_N$, et l'introduction des groupes classiques
$G_1$ et $G'$ n'est en fait absolument pas nécessaire, leur rôle est simplement de permettre un raccourci 
 commode pour certaines notations et d'énoncer facilement l'hypothèse de récurrence.
  Les seuls groupes importants dans cette affaire sont les sous-groupes de Levi
$L_1\times \cdots \times L_r$ de $G$ associés au paramètre $\psi_G$. 
\end{rmq}

\section{Les groupes unitaires}\label{GrUnit}

Dans cette section, nous traitons le cas des groupes unitaires, en adaptant de manière 
relativement évidente ce qui a été fait pour les groupes classiques. Les résultats avaient dans ce cas 
déjà été obtenu par Johnson \cite{Joh2} à la suite des travaux de Clozel \cite{CloBC} dans le cas tempéré.

Nous passons donc en revue les adaptations et compléments
à apporter. Dans cette section, $\mathbf G_N$ désigne maintenant  le groupe $\Res_{\bbC/\bbR}(\GL_N)$ obtenu par restriction 
des scalaires de $\bbC$ à $\bbR$ du groupe $\GL_N$  et $G_N=\mathbf G_N(\bbR)$ est le groupe de ses points réels
de sorte que  que l'on identifie $G_N$ et $\GL_N(\bbC)$.
L'involution de Cartan sur $G_N=\GL_N(\bbC)$ est $\tau : \, g\mapsto {}^t \bar g^{-1}$.

Identifions $\mathbf G_N(\bbC)=\Res_{\bbC/\bbR}(\GL_N)(\bbC)$
à $\GL_N(\bbC)\times \GL_N(\bbC)$. En considérant la sous-algèbre de Cartan de $\caM_N(\bbC)$ constitué des matrices diagonales
que l'on identifie à $\bbC^N$, on obtient par l'isomorphisme  d'Harish-Chandra une identification entre les caractères infinitésimaux pour 
$G_N$ 
et les  orbites sous l'action de $\frS_N\times \frS_N$ dans l'espace $(\bbC\times \bbC)^N$, 
 le premier facteur $\frS_N$ agissant par permutation sur les premières composantes, et 
 le second sur les deuxièmes composantes. On peut voir une telle orbite comme un couple
de deux multi-ensembles à $N$ éléments dans  $\bbC$. Si le caractère infinitésimal  est donné
par un tel couple $(\{s_{i,1}\}_{ i=1,\ldots N}, \{s_{i,2} \}_{ i=1,\ldots N})$, alors celui-ci est entier si 
quels que soient $i,j\in \{1,\ldots, N\}$,  $s_{i,1}-s_{j,1}\in \bbZ$ et  $s_{i,2}-s_{j,2}\in \bbZ$.
 Il est régulier si quels que soient $i,j\in \{1,\ldots, N\}$,  $s_{i,1}-s_{j,1}\neq 0$ et  $s_{i,2}-s_{j,2}\neq 0$.
 C'est le caractère infinitésimal d'une représentation de dimension finie de $G_N$ si et seulement s'il
  est entier et régulier.
 
\medskip

\subsection{Classifications pour $\GL_N(\bbC)$}

La classification des représentations irréductibles de $\GL_N(\bbC)$ est similaire à celle de $\GL_N(\bbR)$,  le théorème \ref{Lgclass},  qui ramène
 la classification des représentations irréductibles à celle des séries discrètes  étant valable   tel quel. 
La différence entre les deux cas est que $\GL_N(\bbC)$ n'admet de séries discrètes que pour $N=1$. 
Celles-ci sont donc les caractères de $\bbC^\times$. Ces derniers sont paramétrés par les 
couples $s_1, s_2\in \bbC$ avec $s_1-s_2\in \bbZ$ : on pose 
\[  \eta(s_1,s_2) :\, \bbC^\times  \longrightarrow \bbC^\times , \qquad z\mapsto z^{s_1}\bar z^{s_2}= \vert z\vert^{s_1+s_2}
\left(\frac{z}{\vert z \vert}\right)^{s_1-s_2}.  \]

Un paramètre de Langlands pour $\GL_N(\bbC)$
est un morphisme algébrique d'image semi-simple
\[ \phi: \, W_\bbC=\bbC^\times \longrightarrow \GL_N(\bbC) .\]
Une telle représentation de $\bbC^\times$ est complètement réductible, les représentations
irréductibles  $\chi_{s_{i,1}, s_{i,2}}$  de $\bbC^\times$ ont été déterminées en (\ref{Chisn}).

Les paquets de Langlands pour $G_N$ sont des singletons, et la classification de Langlands
avec $L$-groupe constitue dans ce cas une classification complète de $\Pi(G_N)$.
Pour une classe de conjugaison de paramètres $\phi$, notons $\Pi_\phi$
la représentation correspondante. Ecrivons $\phi$ comme somme de représentations
iréductibles $\chi_{s_{i,1}, s_{i,2}}$ de $\bbC^\times$. Alors 
$\Pi_\phi=\times_i^\rightarrow \eta(s_{i,1},s_{i,2}   )$. 

\medskip

Un paramètre d'Arthur pour $G_N$ est un morphisme
\[  \psi: W_\bbC \times \SL_2(\bbC) \longrightarrow \GL_N(\bbC).\]
Une telle représentation de $W_\bbC \times \SL_2(\bbC)$ est complètement réductible, et les 
représentations irréductibles de $W_\bbC \times \SL_2(\bbC)$ sont de la forme 
$\chi_{s_1,s_2}\otimes R_a$, $s_1,s_2\in \bbC$, $s_1-s_2\in \bbZ$, $a\in \bbN^\times$.

Les paquets d'Arthur pour $G_N$ sont des singletons. On note $\Pi_\psi$ la représentation
attachée à la classe de conjugaison du paramètre d'Arthur $\psi$.
Si $\psi=\chi_{s_1,s_2}\otimes R_a$, $s_1,s_2\in \bbC$, $s_1-s_2\in \bbZ$, $a\in \bbN^\times$ est un
paramètre d'Arthur irréductible, on pose 
\[ I(s_1,s_2,a) =  \eta(s_1-\frac{a-1}{2}, s_2-\frac{a-1}{2}) \times \eta(s_1-\frac{a-3}{2}, s_2-\frac{a-3}{2})
\times \cdots \times \eta(s_1+\frac{a-1}{2}, s_2+\frac{a-1}{2}). \]
C'est une représentation standard dont l'unique sous-module irréductible de $ I(s_1,s_2,a)$ est alors le caractère 
\[ \eta(s_1,s_2)_a:\, g\mapsto  \eta(s_1, s_2)(\det(g))
 \]
de $G_a$. On a alors $\Pi_\psi=\eta(s_1,s_2)_a$.
Si $\psi$ se décompose en $\psi=\sum_{i=1,\ldots ,r}\psi_i$, alors $\Pi_\psi=\times_i \Pi_{\psi_i}$.

\begin{rmq}
Nous avons utilisé ci-dessus la forme \og  complexe \fg \, du L-groupe et des paramètres de Langlands de $G_N=\GL_N(\bbC)$.
Si nous considérons ce groupe comme le groupe $\Res_{\bbC/\bbR}(\GL_N)$, il nous faut aussi introduire 
$\widehat G_N=\GL_N(\bbC)\times \GL_N(\bbC)$ et ${}^LG_N=\widehat G_N\rtimes W_\bbR$, où l'élément
$j$ de $W_\bbR$ agit sur  $\GL_N(\bbC)\times \GL_N(\bbC)$ par $(g_1,g_2)\mapsto (g_2,g_1)$.
Les classes de conjugaison de paramètres de Langlands  \og complexes \fg\, 
$  \phi: \, W_\bbC \rightarrow \GL_N(\bbC) $ sont alors en bijection (que nous ne rappelons pas ici, voir 
\cite{Bor}) avec  classes de conjugaison de paramètres de Langlands  
$  \phi: \, W_\bbR  \rightarrow {}^ LG_N $, et de m\^eme pour les paramètres d'Arthur.
\end{rmq}

\subsection{Résolution de Johnson de $\eta(-\frac{p}{2},\frac{p}{2})_N$}
Soient $p,N\in \bbZ\times \bbN^\times$ tels que $\vert p \vert >N-1$ et soit  $\eta(-\frac{p}{2},\frac{p}{2})_N$
la caractère de $G_N$ défini ci-dessus. Pour $s\in \frS_N$, notons 
\[ X(s)=\times_{i=1}^N  \eta\left (  \frac{-p-(N-1)}{2}+(i-1),\frac{p-(N-1))}{2}+(s(i)-1)\right). \]
La résolution de Johnson s'écrit alors (remarquons que $q(G_N)=\frac{N(N-1)}{2}$) : 
\[0\rightarrow \eta(-\frac{p}{2},\frac{p}{2})_N \rightarrow X_0\rightarrow X_1\rightarrow \cdots \rightarrow X_{\frac{N(N-1)}{2}}
\rightarrow 0\]
où $X_i$ est la somme directe des $X(s)$, $s\in \frS_N$ de longueur $N$ dans le groupe de Coxeter $\frS_N$.

\subsection{L'espace tordu $\widetilde G_N$}
La seule chose qui change par rapport à la section \ref{torduGLN} est qu'il faut maintenant 
prendre pour définition de l'automorphisme $\theta_N$ de $G_N$ :
\[\theta_N: \, g\mapsto  J_N({}^t\bar g^{-1}) J_N .\]
 On a 
\[ \eta(s_1,s_2)^{\theta_N}= \eta(-s_2,-s_1).\]
Les paramètres d'Arthur irréductibles $\theta_N$-stables sont donc les $\psi=\chi_{-\frac{p}{2}, \frac{p}{2}}\otimes R_N$, 
$p\in \bbZ$. 

\subsection{Caractères tordus des $\eta\left(-\frac{p}{2},\frac{p}{2}\right)_N$} 
Les résultats de la section \ref{caractor} se transposent immédiatement au calcul des caractères tordus
de $\pi= \eta\left(-\frac{p}{2},\frac{p}{2}\right)_N$, on obtient :
\begin{equation}\label{TretaNbis}
\Tr_{\theta_N}(\pi)=\sum_{s \in \frI_N}  (-1)^{\ell_\theta(s)} \; \Tr_{\theta_N}(X(s)).
\end{equation}

\subsection{Groupes unitaires et changement de base}

  Les \og groupes unitaires \fg \; que nous considérons sont ceux qui apparaissent
 dans les données endoscopiques elliptiques simples des $\widetilde G_N$  {\sl cf.}
 \cite{mok}. Il s'agit donc, si $N$ est pair, du groupe $G=\U(\frac{N}{2},\frac{N}{2})$ et 
 si $N$ est impair du groupe  $G=\U(\frac{N-1}{2},\frac{N+1}{2})$.
 \medskip

Pour chacun de ces groupes, on dispose d'un  plongement  du $L$-groupe de $G$ dans  $ {}^LG_N$ : 
\begin{equation}\label{DefStdGBC} \Std_G: \; {}^LG \longrightarrow  {}^LG_N\end{equation}
appelé parfois dans la littérature \og  plongement standard \fg\,  ou  \og  principal\fg.
 La donnée de $(\mathbf G,\Std_G)$ est celle d'une {\sl  donnée endoscopique tordue}
elliptique pour $(G_N,\theta_N)$. Le reste de la section \ref{grcla} se transpose alors aisément.

Soit $\psi_G \in \psi(G)$ un paramètre d'Arthur pour le groupe $G$. Posons 
$\psi=\Std_G \circ \psi_G$. C'est un paramètre d'Arthur $\theta_N$-stable  pour le groupe 
$G_N$. Soit $\Pi_\psi$ la représentation irréductible $\theta_N$-stable  de $G_N$ associée à ce paramètre.
  
  On suppose maintenant de plus que $\psi_G$ est un paramètre d'Adams-Johnson.
   On reprend les notations de la section \ref{AJpackSec} pour les objets  
   attachés au  paramètre $\psi_G$ et celle de la section \ref{PaqArtGN} pour les objets attachés à $\psi$.  
   On explicite ceux-ci pour  les groupes unitaires.
   
Quitte à remplacer $\psi_G$ par un paramètre équivalent, on peut supposer que le    
    centralisateur $\caL$ de $\psi_G(\bbC^\times)$ dans $\widehat G$ est de la forme : 
   \begin{equation}\label{LeviGC} \caL\simeq \caL_1\times \cdots \times \caL_r \end{equation}
où pour tout $i=1, \ldots, r$,  $\caL_i=\GL_{N_i}(\bbC)$, $N_i \in \bbN^\times$.
  On a d'autre part  $\displaystyle \sum_{i=1}^r N_i=N$.

Soit $(\mathbf{Q},\mathbf{L})$ représentant une classe dans $\caS_\caQ$ ({\sl cf.} section \ref{AJpackSec}).   
On a alors un isomorphisme
   \[ L\simeq L_1\times \cdots \times L_r\]
   où  pour tout $i=1, \ldots, r$, $L_i=\U(b_i,c_i)$,   avec $b_i+c_i=N_i$.   Si $N$ est pair 
   $\sum_i b_i=\sum_i c_i=\frac{N}{2}$ et si $N$ est impair $\sum_i b_i=\frac{N-1}{2}$, $\sum_i c_i=\frac{N+1}{2}$ .

Le paramètre $\psi$ se décompose en somme de paramètres irréductibles $\psi=\oplus_{i=1, \ldots ,r}\psi_i$
comme en (\ref{opsii}), où $\psi_i$ se factorise en 
\[ \psi_i :   \xymatrix{ W_\bbR\times \SL_2(\bbC)\ar[rr]^{ \qquad \psi_{G,i} }&&  {}^L L_i \ar[rr]^{\Std_{N_i}  }&& {}^L G_{N_i}},\]
et l'on a $ \psi_{i, \vert\bbC^\times} \simeq \chi_{-\frac{p_i}{2}, \frac{p_i}{2}  }\otimes R_{N_i}$, 
$\Pi_{\psi_i}=\eta\left(  -\frac{p_i}{2},\frac{p_i}{2} \right)_{N_i}$ avec $p_i$ pair.
On note alors  $I(\psi_i)=I( -\frac{p_i}{2},\frac{p_i}{2},N_i)$
 la représentation standard correspondante.

\begin{rmq}\label{parBC}
La condition de parité sur les $p_i$ vient de ce que l'on a choisi $\Std_G$ comme étant le plongement
associé au changement de base principal. On aurait pu faire un autre choix de plongement, le changement de base 
\og tordu \fg\, et la condition de parité aurait été inversée.

\end{rmq}
  Le reste de la section \ref{AJcla} se transpose sans difficulté. Il en est de même des résultats de la section \ref{IndTor}.

  \subsection{Démonstration du résultat principal}
    On reprend les notations de la section précédente, où $\psi_G$ est  
un paramètre  d'Adams-Johnson pour une groupe unitaire quasi-déployé $G=\U(b,c)$ avec $b+c=N$ et $\psi=\St_G \circ \psi_G$.
Notre but est de démontrer l'égalité de distributions sur $\widetilde G_N$ : 
\begin{equation}\label{butfinalBC}
\Tr_{\theta_N} (\Pi_\psi) =\Trans_G^{G_N} (\Theta_{ \Pi_{\psi_G}^{AJ}}^{st}). 
\end{equation}
D'après (\ref{IdStableAJ}), on a une égalité de distributions sur $G$
\begin{equation}
\Theta_{ \Pi_{\psi_G}^{AJ}}^{st}=  (-1)^{q(L^*)} \sum_{\phi_G \in \Phi(\psi_G)} (-1)^{\ell_V(\phi_G)}
  \;  \Theta_ {\widetilde \Pi_ {\phi_G}},
\end{equation}
où tous les termes sont des distributions stables.
On en déduit que 
\begin{equation}\label{11BC}
\Trans_G^{G_N} ( \Theta_{ \Pi_{\psi_G}^{AJ}}^{st}  ) = (-1)^{q(L^*)} 
\sum_{\phi_G \in \Phi(\psi_G)} (-1)^{\ell_V(\phi_G)}  \;  \Trans_G^{G_N} ( \Theta_ {\widetilde \Pi_ {\phi_G}}).
\end{equation}
Or, le résultat de la proposition \ref{TransPsPaq} est encore valide dans le
cas du changement de base, où il est d\^u à Clozel \cite{CloBC}\footnote{La normalisation de Clozel 
ne se fait pas avec les modèles de Whittaker 
mais avec le $\U(n)$-type minimal; comme il n'est pas clair que cette 
normalisation coïncide avec celle qui est utilisée ici, on vérifie que 
le transfert se fait \og sans constante\fg \, comme pour les groupes 
orthogonaux ou symplectiques  dans l'appendice \ref{Moe}.} :  on a pour tout $\phi_G \in \Phi(\psi_G)$
\begin{equation}\label{22bis}
\Tr_{\theta_N} (\widetilde \Pi_\phi) =\Trans_G^{G_N} ( \Theta_ {\widetilde \Pi_ {\phi_G}}) 
\end{equation}
où $\phi=\St_G \circ \phi_G$ et $\Pi_\phi$ est la représentation de $G_N$ associée par Langlands au paramètre 
de Langlands $\phi$.
L'identité (\ref{butfinalBC}) est donc équivalente à 
\begin{equation}\label{tratorBis}
\Tr_{\theta_N} (\Pi_\psi) = (-1)^{q(L^*)}  \sum_{\phi_G \in \Phi(\psi_G)} (-1)^{\ell_V(\phi_G)}  \; \Tr_{\theta_N} 
(\widetilde \Pi_\phi). 
 \end{equation}
La réduction au cas d'un paramètre irréductible marche ici comme dans la section \ref{reduc}.
Il s'agit donc d'établir les résultats analogues à ceux de la section \ref{ParSpehE}, en remplaçant 
$\Speh\left( \left(-\frac{p}{2},\frac{p}{2}\right), n\right)$ par $\eta\left(-\frac{p}{2},\frac{p}{2}\right)_N$ 
(avec ici $p$ pair, {\sl cf.} remarque \ref{parBC}), et 
plus précisément de montrer que :
\[  \sum_{s \in \frI_N}  (-1)^{\ell_\theta(s)} \; \Tr_{\theta_N}(X(s))= (-1)^{q(L^*)} 
\sum_{\phi_G \in \Phi(\psi_G) } (-1)^{\ell_V(\phi_G)} \Tr_{\theta_N}(\widetilde \Pi_\phi).\]  
Exhibons à cette fin  une bijection $s \in \frI_N  \mapsto \phi_G \in \Phi(\psi_G)$
telle que $X(s)=\widetilde \Pi_{\phi}$, $\phi=\Std_G\circ \phi_G$  et $q(L^*)+\ell_V(\phi_G)=\ell_\theta(s)\mod 2$. 

Dans le cas qui nous occupe, le paquet d'Adams-Johnson du paramètre  $\psi_G$
est un singleton,  constitué du  caractère unitaire $\pi_{\psi_G}$ du groupe unitaire quasi-déployé $G=\U(b,c)$
donné par $g\mapsto (\det g)^{p/2}$. 
L'ensemble  $\Phi(\psi_G)$ est alors l'ensemble des paramètres de Langlands 
des représentations standard apparaissant dans les résolutions de $\pi_{\psi_G}$. 
 D'après la section \ref{BBUpq},  ces représentations standard sont paramétrées par $\frI_N^{b,c,\pm}$,
  et leur paramètres de Langlands par  $\frI_N$. Ceci donne la bijection voulue par multiplication par
  l'élément le plus long $w_0$  , de manière analogue à la section \ref{ParSpehE}. 
  
On vérifie maintenant que les longueurs coïncident, en utilisant 
  $\ell_\theta(w_0)=q(L^*)$ de la même manière que
dans la section  \ref{ParSpehE}, en reliant  longueur de Vogan, dimension des orbites dans la variétés des drapeaux
pour le groupe unitaire et $\theta$-longueur dans $\frI_N$.\qed

\appendix

\section{Calcul du facteur de transfert spectral pour les repr\'esentations temp\'er\'ees dans le cas archim\'edien}
\label{Moe}
 Nous avons  besoin de compl\'eter le calcul de  \cite{mezo} pour  avoir pr\'ecis\'ement 
 ce qui est appel\'e par Mezo et Shelstad le   facteur de transfert spectral ({\sl cf.} (\ref{IdEnd})). 
 Le calcul est   ici  beaucoup plus simple que ce qui est fait en \cite{mezo}, mais uniquement pour les groupes 
 classiques (alors que \cite{mezo} est plus g\'en\'eral). Il   devrait pouvoir 
 être fait de manière purement  locale et nous n'avons aucun doute sur le fait que Shelstad et Mezo 
 donneront sous peu une   meilleure preuve. 
L'argument est de type local/global: on montre que l'on peut mettre la situation \`a une place archim\'edienne
 dans une situation globale o\`u   aux autres places on conna\^{\i}t le r\'esultat: c'est la m\'ethode suivie par 
 \cite{Art13} pour d\'eduire le cas $p$-adique   du cas archim\'edien supposé  établi.
  L'int\'er\^et de cet annexe  est   aussi de compl\'eter la preuve des toutes 
 les hypoth\`eses faites en \cite{Art13} et de rendre inconditionnels tous les travaux qui en dépendent.
On s'inspire de \cite{Art13} 6.2.2, on y a juste ajout\'e une construction tir\'ee de \cite{cc} pour traiter le cas  
des caract\`eres centraux non triviaux pour les  repr\'esentations des groupes classiques consid\'er\'es.
On reprend les notations des sections \ref{torduGLN} à \ref{principal}.

\subsection{Rappel des r\'esultats de \cite{mezo}\label{rappel}}
On normalise les facteurs de transfert g\'eom\'etriques et les mesures de sorte que le transfert g\'eom\'etrique  soit 
compatible au produit scalaire elliptique,  \`a un facteur pr\`es, not\'e $c(\tilde{G},{\bf G}')$ dans 
\cite{stabilisationI} \S 4.17.
 Ceci ne fixe pas compl\`etement les choix pour les facteurs de transfert, on reviendra sur ce probl\`eme ci-dessous. 
 On dualise pour obtenir le produit scalaire elliptique sur l'espace des repr\'esentations elliptiques. Pour ce produit 
 scalaire elliptique  la norme  d'une repr\'esentation discr\`ete de $\widetilde G_N$ est $2^{x}$ o\`u $x$ 
 est le nombre de facteur du sous-groupe de Levi de $G_N$ \`a partir duquel on induit,  d'apr\`es \cite{ftlt} \S 7.3.
  (où la d\'efinition de $\iota(\tau)$ est en \cite{ftlt} \S 2.11).
Le facteur d'isom\'etrie vaut $2$ quand le groupe classique consid\'er\'e comme endoscopique est $\Sp(2n,\bbR)$
 ou $\SO(n,n+1)$ (cas {\bf A} et {\bf B} du texte) 
et $4$ pour les groupes $\SO(n,n)$ ou $\SO(n-1,n+1)$  (cas {\bf C} et {\bf D}) car il y a en plus 
un groupe d'automorphismes ext\'erieurs (voir  le calcul de 
$c(\tilde{G},{\bf G}')$ dans  \cite{stabilisationI} \S 4.17).

On fixe le groupe classique $G$ que l'on suppose quasi-d\'eploy\'e et ayant des s\'eries discr\`etes et soit $\psi$ 
un param\`etre pour une repr\'esentation discr\`ete de $\widetilde{G}_N$ se factorisant par un paramètre $\psi_G$ de $G$
(c'est-à-dire que $\psi=\Std_G\circ \psi_G$). 
On note $\pi^G(\psi):=\sum \pi$, o\`u la somme porte sur les s\'eries discr\`etes de $G$ dans le paquet d\'etermin\'e 
par $\psi_G$; rappelons que dans le cas des groupes orthogonaux pairs, on  regroupe deux paquets conjugu\'es
 sous le groupe orthogonal quand ces paquets sont distincts.
On note $\Pi_\psi$ la repr\'esentation irr\'eductible de $G_N$ ayant $\psi$ comme param\`etre de Langlands.
 On étend $\Pi_\psi$ en une représentation de $G_N^+$, normalisée  
 comme dans la section \ref{NormExt} grâce à la donnée de Whittaker.

D'apr\`es \cite{mezo}, il existe un nombre complexe $z(\psi)$ tel que l'on ait l'\'egalit\'e de transfert:
\begin{equation} \label{TrasMezo}
\forall \tilde{f} \in \widetilde{\caH}_N,\qquad  \tr_{\theta_N}( \Pi_\psi)(\tilde{f}))=z(\psi)\; \tr( \pi^G(\psi)(f^G)),
\end{equation}
o\`u $f^G$ est un transfert de $\tilde{f}$ \`a $G({\mathbb R})$.
Notre  but est de montrer que $z(\psi)=1$.

\subsection{Valeur absolue du facteur de transfert spectral} \label{valeurabsolue}

Une première étape est de remarquer comme corollaire des choix faits le résultat suivant.

\begin{lemme}
 Le facteur de transfert, $z(\psi)$ dans la formule  (\ref{TrasMezo}) est de valeur absolue $1$.
\end{lemme}
\dem 
 Il faut calculer la norme elliptique de $\pi^G(\psi)$.  Dans le cas de $\Sp(2n,\bbR)$,   chaque s\'erie discr\`ete 
 a pour norme $1$ et il y en a  $2^{n}$; avec le facteur d'isométrie $2$, on obtient $2^{n+1}$.  La norme de la repr\'esentation elliptique de $\widetilde G_{2n+1}$ est bien ici $2^{n+1}$. 

Pour les groupes orthogonaux, il faut revenir aux $K$-groupes de Kottwitz ({\sl cf.} \cite{Art03}).  La norme de la distribution stable se calcule \`a 
l'aide du transfert \`a partir des $K$-groupes. Si $G=\SO(n,n+1)$,le paramètre $\psi_G$ est celui de s\'eries discr\`etes 
de tous les groupes $\SO(p,q)$ avec $(-1)^p=(-1)^n$,  il y en a  $2^{n-1}$; la norme de la 
repr\'esentation elliptique de $\widetilde G_{2n}$ est effectivement $2^n$ et on a encore le r\'esultat.

Dans le cas de $\SO(n,n)$ ($n$ pair) ou $\SO(n-1,n+1)$ ($n$ impair) (de sorte qu'il y a bien  des s\'eries discr\`etes),
 on distingue suivant que dans la décomposition de $\psi=\oplus_{i=1,\ldots ,r}\psi_i$ de la section 
 \ref{AJcla}, la paramètre $\psi_r$ est trivial ou non ($n_r=0$ ou $n_r\neq 0$). 
Dans le premier cas, les s\'eries discr\`etes sont invariantes sous le groupe orthogonal et leur nombre  se calcule
 en utilisant l'endoscopie tordue pour le groupe orthogonal (ou plut\^ot le $K$-groupe associ\'e), on trouve alors $2^{n-1}$.
  La norme de la repr\'esentation de $\widetilde G_{2n}$ est $2^{n+1}$. Le facteur d'isom\'etrie étant $4$ dans ce cas,
   on obtient bien l'\'egalit\'e. 
Dans le second cas, dans la classe de conjugaison sous $\Or(2n,{\mathbb C})$ de $\psi_G$, il y  a deux 
classes de conjugaison   sous $\SO(2n,{\mathbb C})$. Les s\'eries discr\`etes ne sont pas invariantes  sous le groupe  
orthogonal, il  y a deux paquets différents de celles-ci échangés  par le groupe  orthogonal.
 La norme de la somme des deux paquets est encore le nombre de s\'eries discr\`etes associ\'ees 
 au param\`etre $\psi_G$ et se calcule en utilisant l'endoscopie ordinaire \`a automorphisme pr\`es. On trouve $2^{n-1}$ donn\'ees   endoscopiques possibles   \`a automorphisme pr\`es. Ainsi la norme d'un paquet est  $2^{n-2}$. 
La norme de la repr\'esentation de $\widetilde G_{2n}$ est ici $2^{n}$ et on retrouve encore l'assertion.\qed 

\subsection{La fonctorialit\'e de Langlands pour les repr\'esentations temp\'er\'ees}
Ici on suppose que $F$ est soit un corps p-adique soit un corps archim\'edien. On conna\^{\i}t d\'ej\`a un certain nombre de 
r\'esultats sur la fonctorialit\'e de Langlands en particulier sur les param\`etres. Cela repose sur  l'\'etude de l'espace 
$I_{cusp}(G)$: cet espace vectoriel est l'image des fonctions $\caC^{\infty}_{c}(G)$ dont toutes les termes constants 
(au sens d'Harish Chandra) pour les sous-groupes paraboliques propres de $G$ sont nulles modulo l'espace des 
fonctions $\caC^{\infty}_{c}(G)$ dont toutes les int\'egrales orbitales pour des \'el\'ements semi-simples r\'eguliers sont nulles. 

On sait d\'ecomposer cet espace en utilisant la th\'eorie de l'endoscopie ordinaire. Pour cela et pour un groupe 
quasi-d\'eploy\'e $H$, on d\'efinit $SI_{cusp}(H)$ qui est le quotient de $I_{cusp}(H)$ par le sous-espace engendr\'e
 par les fonctions dont toutes les int\'egrales orbitales stables sont nulles. 

Ce qui nous importe ici est de comprendre $SI_{cusp}(G)$ en supposant que $G$ est quasi-d\'eploy\'e.
 A une repr\'esentation elliptique irr\'eductible on associe un pseudo-coefficient qui est un \'el\'ement de
  $I_{cusp}(G)$; on sait que la projection orthogonale (pour le produit scalaire elliptique) de ce pseudo-coefficient 
  sur $SI_{cusp}(G)$ est nulle si $\pi$ n'est pas une s\'erie discr\`ete. 
  Si $\pi$ est une s\'erie discr\`ete on note $f_{\pi,st}$
   cette projection.  Si $F$ est archim\'edien, on sait que $f_{\pi,st}$ est la somme des pseudo-coefficients pour
    les s\'eries discr\`etes dans le m\^eme paquet de Langlands que $\pi$. Si $F$ est $p$-adique, tant que l'on 
    ne peut pas utiliser \cite{Art13}, c'est-\`a-dire tant que l'on n'a pas d\'emontr\'e ce qui est l'objet de ces lignes,
     on connait uniquement les propri\'et\'es suivantes:

\begin{itemize}
\item[$\bullet$] si $\pi\not\sim \pi'$ sont des s\'eries discr\`etes, alors $f_{\pi,st}$ est un multiple de $f_{\pi',st}$ ou ces deux
 fonctions sont orthogonales pour le produit scalaire elliptique, 

\item[$\bullet$] dans tous les cas $SI_{cusp}(G)$ est engendr\'e comme espace vectoriel par l'ensemble des fonctions
 $f_{\pi,st}$ quand $\pi$ parcourt l'ensemble des s\'eries discr\`etes.
\end{itemize}

\medskip

Dans certains cas, comme ceux qui nous occupent  ici, l'espace $SI_{cusp}(G)$ peut se comprendre grâce à 
l'endoscopie tordue; pour cela on utilise le fait que $G$ fait partie d'une donn\'ee endoscopique elliptique d'un 
$\widetilde G_N$ (cet espace tordu est défini pour $F$ $p$-adique de la même façon que pour $F=\bbR$. De manière 
plus générale les notations du texte concernant $\GL$ et les groupes classiques s'adaptent de manière évidente
à un corps $F$ quelconque). 
En fait ici, la situation est particuli\`erement simple: si $G\neq \Sp(2n)$, alors $N=2n$ est pair et $G=\SO(2n)$ 
quasi-d\'eploy\'e ou $G=\SO(2n+1)$ et $G$ d\'etermine uniquement une donn\'ee endoscopique elliptique de 
$\widetilde G_N$. Si $G=\Sp(2n)$, les donn\'ees endoscopiques elliptiques qui ont comme groupe sous-jacent 
$G$ sont en bijection avec les caract\`eres quadratiques de $F^*$;  dans ce qui suit on consid\`ere la donn\'ee
 endoscopique elliptique correspondant au caract\`ere quadratique trivial si $G=\Sp(2n)$. 

On revient \`a $\pi$ une s\'erie discr\`ete irr\'eductible de $G$ et \`a $f_{\pi,st}$; si $G=\SO(2n)$, on note 
$f_{\pi,st}^{O}$ soit $f_{\pi,st}$ si $\pi$ est stable sous l'action de $\Or(2n)$ soit la somme $f_{\pi,st}+ f_{\pi',st}$ 
o\`u $\pi'$ est l'image de $\pi$ sous l'action de $\Or(2n)$. On sait a priori que $f_{\pi,st}$ ou $f_{\pi,st}^O$ sont 
des transferts pour l'ensocopie tordue d'un \'el\'ement de $I_{cusp}(\widetilde G_N)$ et le probl\`eme est de 
d\'eterminer cet \'el\'ement.

On a d\'ej\`a des renseignements assez pr\'ecis pour r\'esoudre ce probl\`eme. A $\pi$ est associ\'e un morphisme $\phi$
de $W'_{F}$ dans $\GL_N({\mathbb C})$ o\`u $W'_{F}=W_{F}$ si $F$ est un corps archim\'edien et 
$W'_{F}=W_{F}\times \SL_2({\mathbb C})$ si $F$ est un corps p-adique. 
On note $\Pi_\phi$ la repr\'esentation irr\'eductible de $G_N$ qui correspond \`a $\phi$ par la 
correspondance  de Langlands. On sait que $\Pi_{\phi}$ définit une repr\'esentation elliptique de $\widetilde G_N$,  
c'est-\`a-dire qu'il existe  un sous-groupe de Levi  $M$  de $G_N$ et $\sigma$ une s\'erie discr\`ete de $M$ tels 
que $\pi$ soit l'induite de $\sigma$ et tels que  le normalisateur de $(M,\sigma)$ dans $\widetilde G_N$ est un 
espace homog\`ene sous l'action de $M$, ou autrement dit que   $\mathrm{Norm}_{\widetilde G_N}(M,\sigma)/M$ 
a exactement un \'el\'ement. 
 On étend cette repr\'esentation en une repr\'esentation
  de  $\widetilde G_N^+$ comme dans la section $\ref{NormExt}$ grâce à la donnée de Whittaker.

A une telle repr\'esentation elliptique, $\Pi_\phi$ est associ\'e un pseudo-coefficient 
$f_{\phi}$. Pour le moment, on sait 
qu'il existe $z'_{\pi}\in {\mathbb C}^*$ tel que $f_{\pi,st}$ soit un transfert de $z'_{\pi}f_{\phi}$. 
Comme $\Pi_\phi$ n'est pas  de norme $1$ pour le produit scalaire elliptique, il est pr\'ef\'erable de 
traduire ce transfert en termes de transfert de trace de  repr\'esentations: d'abord on \'ecrit la 
combinaison lin\'eaire stable de s\'eries discr\`etes associ\'ees \`a $f_{\pi,st}$ ou 
$f_{\pi,st}^O$ sous la forme:
$$
\pi_{st}:=\sum_{\pi'\in \Pi_{st}} a(\pi') \; \tr  \pi',
$$
o\`u $\Pi_{st}$ est un ensemble fini de s\'eries discr\`etes contenant $\pi$ et o\`u $a(\pi')\in {\mathbb C}^\times $
 avec $a(\pi)=1$. Si $F$ est archim\'edien, on a m\^eme $a(\pi')=1$ pour tout $\pi'\in \pi_{st}$. A posteriori cela sera 
 aussi vrai si $F$ est p-adique mais pour le moment on ne le sait pas. Donc au sujet du transfert, on sait qu'il existe 
 $z_{\pi}\in {\mathbb C}^\times $ tel que pour tout $\tilde f\in \widetilde \caH_N$, on ait:
\begin{equation}
z_{\pi} \; \tr_{\theta_N}\left( \Pi_\phi(\tilde f)\right)=
\sum_{\pi'\in \Pi_{st}}a(\pi') \; \tr\left( \pi'(f^{G}) \right), \end{equation}
o\`u $f^G$ est le transfert de $\tilde f$ \`a $G$.

Nous voulons d\'emontrer  que $z_{\pi}=1$ si $F$ est archim\'edien.
et pour cela,  nous allons utiliser  quelques cas simples de ce résultat dans le cas o\`u $F$ est p-adique.

\subsection{Les cas simples des groupes p-adiques et des repr\'esentations de Steinberg}
Comme dans le paragraphe précédent,  on fixe encore $G=\mathbf G(F)$ un groupe classique quasi-d\'eploy\'e,  et 
 ici $F$ est un corps $p$-adique. La repr\'esentation 
de Steinberg est bien d\'efinie pour $G$; c'est l'image de la repr\'esentation triviale par l'involution d'Iwahori-
Matsumoto.  C'est une s\'erie discr\`ete.
Comme la repr\'esentation triviale a un caract\`ere qui est une distribution stable,
 il en est de m\^eme de la repr\'esentation  de Steinberg. On note $\mathrm{St}(G)$ cette repr\'esentation.
On veut calculer son transfert pour l'endoscopie tordue vers $\widetilde G_N$. On d\'efinit a priori 
ce qui sera le r\'esultat.  Si $G$ est un groupe endoscopique principal de $\widetilde G_N$,  
c'est-\`a-dire si $G\neq \SO(2n)$, alors on consid\`ere   $\mathrm{St}(N)$ la repr\'esentation de Steinberg de $G_N$.
 Si $G=\SO(2n)$, on consid\`ere l'induite   $\mathrm{St}(N-1)\times \eta_{G}$ o\`u
  $\eta_{G}$ est le caract\`ere trivial de $F^\times$ si $G$ est d\'eploy\'e et est le
   caract\`ere d\'efini par l'extension quadratique qui 
  d\'eploie $G$ sinon. Pour unifier les notations, on note $\widetilde{\mathrm{St}}(G)$ cette repr\'esentation.

\begin{lemme} Pour toute fonction $\tilde{f}\in \widetilde \caH_N$, 
$$
\tr_{\theta_N}\left( \widetilde{\mathrm{St}}(G)(\tilde{f})\right)= \tr\left( \mathrm{St}(G)(f^G)\right),
$$
o\`u $f^G$ est un transfert de $\tilde{f}$ \`a $G$.
\end{lemme}
\dem 
En \cite{seriesdiscretes} il est montr\'e (comme nous l'avons rappel\'e ci-dessus) qu'il existe un scalaire 
$z$ tel que l'on ait pour tout $\tilde{f}\in  \widetilde \caH_N$, 
$$
\tr_{\theta_N}\left(  \widetilde{\mathrm{St}}(G)(\tilde{f})\right)=z \; \tr\left(  \mathrm{St}(G)(f^G)\right),
$$
et le point est de montrer que $z=1$. Comme le transfert commute au module de Jacquet, on se ram\`ene ais\'ement 
\`a  montrer cette \'egalit\'e pour $n=1$. Les cas de $\SO(3)$ et $\SL(2)$ sont simples: dans ces cas, on a:
\begin{equation}\label{difftr}
\mathrm{St}(G)= \Ind_{B}^{G}(\vert\, \vert ^x)-1_{G},
\end{equation}
o\`u $1_{G}$ est la repr\'esentation triviale de $G$ et $x$ vaut $1/2$ si $G=\SO(3)$ et $1$ si $G=\SL(2)$. Distinguons 
encore suivant que $G=\SL(2)$ ou $\SO(3)$. Dans le premier cas,  on obtient que le transfert de $\mathrm{St}(G)$ est la 
diff\'erence des traces tordue de l'induite $\vert\, \vert \times 1 \times \vert\, \vert^{-1}$ et de la trace tordue de la 
repr\'esentation triviale, l'action de $\theta$ est celle qui est donn\'ee par la normalisation de Whittaker. 
On remarque que l'induite $\vert\, \vert \times 1 \times \vert\, \vert^{-1}$ a exactement deux sous-quotients
 irr\'eductibles dont la   classe d'isomorphie est stable  sous l'action de $\theta$ qui sont la repr\'esentation 
 triviale de $\GL_3(F)$ et la repr\'esentation de Steinberg; chacune de ces  repr\'esentations a multiplicit\'e un 
 et h\'erite de la \og bonne \fg \,  action de $\theta$ quand on a normalis\'e l'action de $\theta$ sur
   l'induite de fa\c{c}on \`a commuter \`a la fonctionnelle de Whittaker. Donc la trace tordue de 
   l'induite est exactement la somme    des traces tordues. 
Ainsi la trace tordue de $\widetilde{\mathrm{St}}(G)$ se calcule en faisant la diff\'erence des traces tordues comme dans
 \ref{difftr}  ce qui donne  le r\'esultat cherch\'e. 
Dans le cas de $G=\SO(3)$, c'est encore plus simple car l'induite \`a consid\'erer est   
$\vert\, \vert^{1/2}\times \vert\, \vert^{-1/2}$ qui est de longueur exactement deux contenant comme
 sous-quotients irr\'eductibles la repr\'esentation de Steinberg et la repr\'esentation triviale.

Il reste le cas o\`u $G=\SO(2)$ d\'eploy\'e ou non. Ici $\mathrm{St}(G)$ est la repr\'esentation triviale du tore. 
On va d'abord traiter (avec l'aide de Waldspurger) le cas de $\SO(2)$ non d\'eploy\'e. Le cas archim\'edien 
(dont on n'a pas besoin ici) est partiellement trait\'e dans \cite{miniature} il manque juste l'interpr\'etation d'un signe
en  termes de facteurs de transfert.
On consid\`ere la repr\'esentation induite $1\times \eta$ de $\GL_2(F)$ o\`u $\eta$ est le caract\`ere quadratique 
d\'eterminant une extension $E$ de $F$ tel que $\SO(2)$ soit les \'el\'ements de normes un de $E^\times$.
 On note $\tilde{\pi}$ l'extension de  cette repr\'esentation \`a $\widetilde G_2$ en imposant que
  $\tilde{\pi}(\theta)$ pr\'eserve le mod\`ele de Whittaker. On peut
  d\'ecrire facilement $\tilde{\pi}(\theta)$ (c'est analogue \`a \cite{miniature}): notons $B$ le sous-groupe 
  de Borel de $\GL_2(F)$ 
  et $B_{0}$ son intersection avec $\SL_2(F)$. Alors la restriction de $\Ind_{B}^{\GL_2(F)}(1\otimes \eta)$
   \`a $\SL_2(F)$ est
   exactement $\Ind_{B_{0}}^{SL_2(F)}(1\otimes \eta)_{\vert B_{0}}$, c'est-\`a-dire l'induite du caract\`ere
    $\eta$ de $B_{0}$
    \`a $\SL_2(F)$. Or cette derni\`ere induite est la somme de deux sous-modules irr\'eductibles, $\pi^+\oplus \pi^-$ 
    o\`u par      d\'efinition $\pi^+$ est la sous-repr\'esentation ayant le mod\`ele de Whittaker. 
    On remarque qu'ici $\theta(g)=g/\det(g)$ et en 
    particulier $\theta$ agit trivialement sur $\SL_2(F)$. Ainsi $\tilde{\pi}(\theta)$ agit n\'ecessairement 
    par $+1$ sur $\pi^+$.
    Evidemment $\tilde{\pi}(\theta)$ ne peut \^etre identiquement $1$ sinon $\pi(g)=\pi(\theta(g))$
     pour tout $g\in \GL_2(F)$ ce qui
     force le caract\`ere central de $\pi$ \`a \^etre trivial alors qu'il ne l'est pas. Ainsi $\tilde{\pi}(\theta)$
      vaut n\'ecessairement $-1$
      sur $\pi^-$. Pour $g\in \SL_2(F)$, le caract\`ere de $\tilde{\pi}$ en $g\theta$ est exactement 
      $\pi^+(g)-\pi^-(g)$. Cette diff\'erence 
      de caract\`ere a \'et\'e calcul\'ee en terme d'endoscopie pour $\SL_2(F)$ par Labesse et Langlands: c'est le transfert 
      endoscopique du caract\`ere trivial du m\^eme $\SO(2,F)$ mais ce n'est pas la m\^eme correspondance endoscopique. 
      En reprenant les d\'efinitions, on v\'erifie que si $\delta\in \SO(2,F)$ correspond \`a $g\in \SL_2(F)$
       (ou plut\^ot \`a sa 
      classe de conjugaison stable) dans l'endoscopie pour $\SL_2(F)$ alors $\delta^2$ correspond \`a la classe de 
      conjugaison stable de $g$ pour l'endoscopie tordue pour $\GL_2$. Fort heureusement, pour un tel couple $(\delta,g)$ 
      on a l'\'egalit\'e des facteurs de transfert $$\Delta_{\SL_2(F)}(\delta,g)=\Delta_{\GL_2(F).\theta}(\delta^2,g).$$
Soit donc $g\in \SL_2(F)$ tel que les valeurs propres de $g$ soient
 $(\alpha(g), \alpha(g)^{-1})$ avec $\alpha(g)$ dans $E^\times$.
 N\'ecessairement $\alpha(g)$ est de norme un et d\'etermine donc un \'el\'ement $\delta\in \SO(2,F)$. Ainsi $\delta$ et 
 $\delta^{-1}$ sont les \'el\'ements de $\SO(2,F)$ tel que le facteur de transfert $\Delta_{\SL_2(F)}(\delta,g)\neq 0$ et 
 d'apr\`es \cite{labesselanglands}, on a l'\'egalit\'e de caract\`eres:
$$
\pi^+(g)-\pi^-(g)= \sum_{\delta\in \SO(2,F)} \Delta_{\SL_2(F)}(\delta,g).
$$
D'apr\`es ce que l'on a vu ci-dessus, pour ce point $g$, le terme de gauche est le caract\`ere tordu de $\tilde{\pi}$ 
et celui de droite est $\sum_{\delta\in \SO(2,F)} \Delta_{GL_2(F).\theta}(\delta^2,g)$ et ce sont les deux \'el\'ements 
$\delta^2$ qui correspondent \`a la classe de conjugaison tordue stable de $g$. On obtient donc bien $z=1$.

Il reste le cas o\`u $\SO(2,F)$ est d\'eploy\'e, c'est-\`a-dire $\SO(2,F)\simeq F^\times$. Dans ce cas, $\SO(2,F)$ n'est pas 
associ\'e \`a une donn\'ee endoscopique elliptique de $\widetilde G_2$ mais est \`a une donn\'ee endoscopique 
elliptique de l'espace de Levi $\tilde{M}:=(F^\times \times F^\times ).\theta$ 
o\`u $\theta(t,t')=(t^{'-1},t^{-1})$. Le transfert des 
fonctions est alors facile, $$\tilde{f} \mapsto \int_{F^{*}}\tilde{f}_{\tilde{M}}(t t^{'-1},t')dt'.$$
L'analogue de la repr\'esentation $\tilde{\pi}$ du paragraphe pr\'ec\'edent est tout simplement l'induite irr\'eductible du
 caract\`ere trivial du Borel de $\GL_2(F)$, l'action de $\theta$ \'etant l'action triviale et la formule de transfert des
  caract\`eres est alors imm\'ediate. \qed

\subsection{Le cas des repr\'esentations $\theta$-induites}
Le transfert commute \`a l'induction c'est une des propri\'et\'es cl\'es mais on va utiliser cette propri\'et\'e uniquement 
dans un cas tr\`es particulier. On suppose que $G$ est d\'eploy\'e, il contient donc un sous-groupe parabolique, $P$, 
de Levi isomorphe \`a $\GL_n(F)$; si $G\neq \SO(2n)$ \`a conjugaison pr\`es, il n'y a m\^eme qu'un seul tel sous-groupe
 parabolique, par contre si $G=\SO(2n)$, il y en a deux qui sont conjugu\'es sous $\Or(2n,F)$. On fixe $\omega$ une 
 repr\'esentation cuspidale irr\'eductible unitaire de $\GL_n(F)$. Soit $\mu$ un caract\`ere de $F^\times$; on consid\`ere 
 d'une part la repr\'esentation induite $\pi(\mu):=\Ind_{P}^G(\omega\mu)$ de $G$ (si $G=\SO(2n)$ on fixe l'un des choix) 
 et d'autre part la repr\'esentation induite de $\GL_N(F)$, $\tilde{\pi}(\mu)$ 
 qui est soit $\omega\mu \times \omega^*\mu^{-1}$
  si $N=2n$ soit $\omega\mu \times 1 \times \omega^*\mu^{-1}$ si $N=2n+1$. On munit $\tilde{\pi}(\mu)$ de l'action de 
  $\theta$ qui fixe la fonctionnelle de Whittaker; cette action co\"{\i}ncide avec l'action canonique de $\theta$ suivante: 
  on note $V_{\omega}$ un espace dans lequel on r\'ealise $\omega$; alors $\omega\mu$ se r\'ealise \'evidemment dans
   $V_{\omega}$ sans aucun choix suppl\'ementaire. Par contre, on fixe un automorphisme $A$ qui identifie l'espace 
   de $\omega^*$, $V_{\omega^*}$ avec $V_{\omega}$ et qui entrelace l'action de $\omega^*$ avec l'action
    $g\in \GL_n(F)\mapsto \omega(\, J_{n}^{-1}\, ^tg^{-1}J_{n})$. La repr\'esentation induite se r\'ealise dans l'espace 
    des fonctions sur $\GL_N(F)$ lisses \`a support compact modulo le bon parabolique et \`a valeurs dans 
    $V_{\omega}\otimes V_{\omega^*}$. Alors $\theta$ agit naturellement sur cet espace par $$
\forall g\in GL(N,F), \qquad \theta.f(g)=A^{-1}\otimes A f(\theta(g)).
$$
Cette action co\"{\i}ncide avec l'action de $\theta$ pr\'eservant la fonctionnelle de Whittaker. La trace tordue
 se calcule via le terme constant comme la trace tordue pour la repr\'esentation $\omega\mu \otimes \omega^*\mu^{-1}$ 
 de $(\GL_n(F)\times \GL_n(F)).\theta'$ o\`u $\theta'(g,g')= (\theta_{n}(g'),\theta_{n}(g)) $ (avec $\theta_{n}$ 
 d\'efini comme 
 $\theta$). L'endoscopie avec $\GL_n(F)$ est alors triviale. On  d\'eduit du fait que le transfert commute \`a l'induction,  
 l'\'egalit\'e  de transfert:
$$\forall \tilde{f}\in \widetilde \caH_N, 
\qquad
\tr_{\theta_N}\left( \tilde{\pi}(\mu)(\tilde{f})\right)=\tr\left( \pi(\mu)(f^G)\right).
$$
Il faut remarquer que $f^G$ est invariant sous l'action de $\Or(2n,F)$ si $G=\SO(2n,F)$ est d\'eploy\'e et que ce
 r\'esultat ne d\'epend donc pas du choix du parabolique pour d\'efinir $\pi(\mu)$.

\subsection{Globalisation d'apr\`es Arthur et Chenevier-Clozel}
Le paragraphe qui suit est tir\'e de \cite{Art13} 6.2.2; comme dans cette r\'ef\'erence il est suppos\'e que le corps 
de nombres sur lequel on se place a \og beaucoup \fg\,  de places archim\'ediennes et que 
nous voulons le r\'esultat sur ${\mathbb Q}$, on 
redonne les arguments. De plus on a ajout\'e une construction tir\'ee de \cite{cc} qui permet de ne pas avoir de 
probl\`eme avec les caract\`eres centraux des repr\'esentations des groupes classiques.  L\`a aussi on d\'evie du
 cadre de \cite{cc}, puisque dans cet article c'est le cas de $\SO(2n+1)$ qui \'etait trait\'e alors que l'id\'ee nous sert 
 justement \`a traiter le cas des groupes classiques ayant un centre. La situation est donc maintenant globale, notre 
 groupe classique $\mathbf G$ est défini sur un corps de nombres, noté $F$.

Le premier r\'esultat utilis\'e est un r\'esultat local d\'emontr\'e en \cite{lefschetz} 5.1
(voir aussi \cite{CR2}): on fixe $\pi$ une s\'erie discr\`ete
 de $G({\mathbb R})$ et on consid\`ere $f_{\pi,st}$. Ici $f_{\pi,st}$ est la somme des pseudo-coefficients pour
  les s\'eries discr\`etes de $G(\bbR)$ dans le m\^eme paquet stable que $\pi$. C'est une fonction cuspidale et stable. 
  Alors les int\'egrales orbitales de $f_{\pi,st}$ sont nulles en tout point non semi-simple.

On fixe un caract\`ere infinit\'esimal $\nu$ de $G({\mathbb R})$ que l'on suppose entier et loin des murs. 
Cela entra\^{\i}ne que pour toute repr\'esentation $\pi_{\infty}$ de $G({\mathbb R})$ unitaire, $\pi_{\infty}$ 
est une s\'erie discr\`ete: en effet d'apr\`es \cite{salamanca}, une telle repr\'esentation est l'une des 
repr\'esentations construites par Vogan et Zuckermann (\cite{VZ}). Or ces repr\'esentations n'ont un 
caract\`ere infinit\'esimal \og loin des murs \fg\,  que si elles sont des s\'eries discr\`etes. On reprend aussi la 
notation de \cite{Art13} 6.2.2 , $m\nu$ avec $m\in {\mathbb N}$ pour faire tendre $\nu$ vers l'infini dans 
une chambre de Weyl fix\'ee et \'evidemment $m\nu$ est toujours entier et loin des murs. On remarque
 alors aussi que quand $m\nu$ est fix\'e, les s\'eries discr\`etes de $\mathbf G({\mathbb R})$ ayant ce caract\`ere 
 infinit\'esimal sont dans un unique paquet de Langlands si $G\neq \SO(2n)$ et par contre sont dans deux 
 paquets de Langlands conjugu\'es sous $\Or(2n,{\mathbb R})$ si $G=\SO(2n)$. On note donc $f_{m\nu,st}$ 
 la somme des pseudo-coefficients comme ci-dessus. On remarque aussi pour la suite que si $\pi_{\nu}$ 
 est une s\'erie discr\`ete de $\mathbf G({\mathbb R})$ de caract\`ere infinit\'esimal $\nu$, son caract\`ere central 
 ne d\'epend que de $\nu$, on le note $\chi_{\nu}$ (cela n'a d'int\'er\^et que si $\mathbf G({\mathbb R})$ a un centre 
 non trivial) et $\chi_{m\nu}=\chi_{\nu}$: en effet pour calculer le caract\`ere central, on \'ecrit $\pi_{\nu}$ 
 comme sous-quotient d'une s\'erie principale;   dans cette s\'erie principale les exposants des caract\`eres 
 par leur multiple par $m$ sans changer leur restriction \`a $-1\in {\mathbb R}$ et la nouvelle s\'erie principale
  contient comme sous-quotient une s\'erie discr\`ete de caract\`ere infinit\'esimal $m\nu$. D'o\`u l'assertion.

On fixe $V$ un nombre fini de places de ${\mathbb Q}$ contenant la place archim\'edienne, le $V_{ram}$ 
de \cite{stabilisationX} (c'est-\`a-dire toutes les places de petites caract\'eristiques r\'esiduelles) et telles 
que le lemme fondamental tordu soit vrai hors de $V$. On suppose que $V$ a au moins deux places finies.
 On en fixe une, not\'ee $v_{0}$ et on suppose que $\mathbf G({\mathbb Q})$ est d\'eploy\'e en cette place 
 (la condition n'est restrictive que si $\mathbf G=\SO(2n)$ est quasi-d\'eploy\'e non d\'eploy\'e). En la place 
 $v_{0}$, on note $F_{0}$ le compl\'et\'e de ${\mathbb Q}$ et on fixe une repr\'esentation cuspidale
  irr\'eductible unitaire $\omega_{0}$ de $\GL_n(F_0)$ tel que pour tout caract\`ere non ramifi\'e $\nu$ 
  de $F_{0}^\times$, la repr\'esentation $\omega_{0}\nu$ n'est pas isomorphe \`a sa contragr\'ediente 
  $\omega_{0}^*\nu^{-1}$. Alors pour tout $\nu$ comme pr\'ec\'edemment et pour tout parabolique 
  $P$ de $\mathbf G(F_{0})$ de sous-groupe de Levi $\GL_{n}(F_0)$ l'induite  $\Ind_{P}^{G(F_{0})}\omega_{0}\nu$ 
  est irr\'eductible.
On reprend une id\'ee de \cite{cc}: on fixe une fonction $f_{0}$ sur $\mathbf G(F_{0})$ dont la composante de 
Paley-Wiener est nulle pour toute composante $(M,\sigma)$ o\`u $M$ est un sous-groupe de Levi de
 $\mathbf G(F_{0})$ diff\'erent de $\GL_n(F_0)$ ou \'egal \`a $\GL_n(F_0)$ mais alors $\sigma \not\sim \omega \nu$
 pour tout $\nu$ caract\`ere unitaire non ramifi\'e de $F_{0}^\times$. On impose en plus \`a $f_{0}$ de v\'erifier
  $f_{0}(1)\neq 0$ ce qui ne pose pas de probl\`eme en utilisant la formule de Plancherel. On note 
  $\chi_{\omega}$ le caract\`ere central de $\omega$ et on remarque que pour $z$ dans le centre de
   $\mathbf G(F_{0})$ (c'est-\`a-dire $z=1$ ou $z=-1$ ce qui ne peut se produire que si $G\neq \SO(2n+1)$), 
   on a $z.f_{0}=\chi_{\omega}(z)f_{0}$. 

On d\'efinit $\tilde{f}_{0}\in I(\widetilde G_N)$ de tel sorte que $f_{0}$ en soit un transfert: 
la composante de Paley-Wiener de cette fonction est nulle sauf sur la composante associ\'e \`a
 $\tilde{M}$, $\tilde{\omega}$ o\`u $\tilde{M}$ est le transfert du Levi $M=\GL_n(F_{0})$ de $\mathbf G(F_{0})$ et 
 $\tilde{\omega}$ est le transfert de $\omega$ c'est-\`a-dire essentiellement $\omega \otimes \omega$ 
 ou $\omega\otimes \omega \otimes 1$. Et les termes constants de $\tilde{f}_{0}$ ont  pour transfert le terme
  constant de  $f_{0}$ pour ces Levi mais l'op\'eration de transfert est ici assez \'evidente.

On fixe $\nu$ et $\omega$ et on demande \`a $\omega$ d'avoir comme caract\`ere central restreint au
 centre de $\mathbf G({\mathbb R})$ le caract\`ere $\chi_{\nu}$.
On note $V'$ le sous-ensemble de $V$ dont on a enlev\'e $v_{0}$ et l'infini.
Pour toute place $v$ de ${\mathbb Q}$ dans $V'$, on fixe $f_{v}$ un pseudo-coefficient de la repr\'esentation
 de Steinberg. On note $1_{K^V}$ la fonction caract\'eristique du compact $\mathbf G({\mathcal O}_{v})$. 
Alors la fonction $f^G(m\nu):=1_{K^V}f_{0}f_{m\nu}\times_{v\in V'}f_{v}$ est invariante sous l'action du centre 
de $\mathbf G({\mathbb Q})$. C'est en plus une fonction cuspidale en au moins deux places et dont les int\'egrales
 orbitales en les \'el\'ements non semi-simples sont nuls (car cela est vrai pour $f_{m\nu}$).
On fixe aussi $\tilde{f}_{m\nu}$ une fonction cuspidale de $I_{cusp}(\widetilde G_N)$ dont le transfert
 \`a $\mathbf G({\mathbb R})$ est la fonction $f_{m\nu}$ et pour tout $v\in V'$, une fonction 
 cuspidale dans $I_{cusp}(\widetilde{ \mathbf G}_N(F_{v'})))$ qui a $f_{v'}$ comme 
 transfert \`a $\mathbf G(F_{v'})$ (ici on utilise bien que ce ne soit pas utile le fait que 
 $f_{v'}$ est naturellement dans $SI_{cusp}^G$). 
 On pose $\tilde{f}(m\nu):= 1_{\tilde{K}^V}\tilde{f}_{0}\tilde{f}_{m\nu}\times_{v\in V'}\tilde{f}_{v}$ 
 o\`u $\tilde{K}^V$ est $\GL_N({\mathcal O}^V).\, \theta_N$. Comme la notation le sugg\`ere, $f^G$ 
 est un transfert de $\tilde{f}$ pour l'endoscopie tordue.

Soit $c^V$ un caract\`ere pour l'alg\`ebre de Hecke sph\'erique de $\mathbf G({\mathbb Q}^V)$ et $\tilde{c}^V$ son
 transfert \`a $\GL_N({\mathbb Q}^V)$.  Pour $m,\nu,\omega$ comme ci-dessus; on fixe aussi $\mu_{0}$ un
  caract\`ere non ramifi\'e de $F_{0}^\times$ et on note $\Pi^G(c^V,m\nu, \mu_{0})$ l'ensemble des repr\'esentations 
  automorphes de carr\'e int\'egrable de $\mathbf G({\mathbb A}_{{\mathbb Q}})$ qui ont $m\nu$ comme caract\`ere 
  infinit\'esimal \`a l'infini, qui sont non ramifi\'ees hors de $V$ et se transforment sous le caract\`ere $c^V$ pour 
  l'action de l'alg\`ebre de Hecke sph\'erique, qui en $V'$ sont les repr\'esentations de Steinberg et qui en 
  $v_{0}$ sont isomorphes \`a l'une des induites $\Ind_{P}(\mu_{0}\otimes \omega)$ pour $P$ l'un des paraboliques 
  de $G(F_{0})$ de Levi $\GL_n(F_0)$. 
On note $\tilde{\pi}(\tilde{c}^V,m\nu,\mu_{0})$ la repr\'esentation de $\GL_N({\mathbb A}_{{\mathbb Q}})$ qui hors 
de $V$ est la repr\'esentation non ramifi\'ee correspondant \`a $\tilde{c}^V$ qui en toute place de $V'$ est le
 transfert de la repr\'esentation de Steinberg de $\mathbf G(F_{v'})$ et qui en $v_{0}$ est l'induite  
 $\omega\mu_{0}\times \omega^*\mu_{0}^{-1}$ si $N$ est pair et de $\omega\mu_{0}\times \omega^*\mu_{0}^{-1}\times 1$
  si $N$ est impair.

Arthur a d\'emontr\'e le r\'esultat suivant (en des termes diff\'erents); on note ici $i(G)$ la
 constante $i(\widetilde G_N,G)$ qui intervient dans la stabilisation de la formule des traces, 
 la seule chose  importante pour nous est que c'est un nombre
 rationnel positif.

\begin{thm}  [\cite{Art13} 6.2.2] Pour $m$ suffisamment grand, il existe $c^V$ et $ \mu_{0}$ tel que l'ensemble 
$\Pi^G(c^V,m\nu,\mu_{0})$ soit non vide. Alors $\tilde{\pi}(\tilde{c}^V,m\nu,\mu_{0})$ est une repr\'esentation
 automorphe qui intervient dans la partie discr\`ete de la formule des traces pour $\widetilde {\mathbf G}_N$.

De plus pour toute fonction $\tilde{f}^{\infty}_{V}$ dans $\otimes_{v\in V-\{\infty\}}I(\widetilde{\mathbf G}_N(F_{v}))$ 
dont  on note $f^{G,\infty}_{V}$ un transfert \`a $\mathbf G({\mathbb Q}_{V-\{\infty\}})$ on a l'\'egalit\'e de transfert:
$$
\tr_{\theta_N}\left( \tilde{\pi}(\tilde{c}^V,m\nu,\mu_{0})(\tilde{f}^{\infty}_{V}f_{m\nu}1_{\tilde{K}}^V)\right)=
 i(G) \sum_{\pi\in \Pi^G(c^V,m\nu, \mu_{0})}m(\pi) \; \tr\left( \pi(f^{G,\infty}_{V} f_{m\nu}1_{K^V})\right),
$$o\`u $m(\pi)$ est un entier positif.
\end{thm}

\dem 
On reprend la d\'emonstration de {\sl loc.cit.}
On note $I^G_{geo}()$ et $I^G_{spec}()$ le c\^ot\'e g\'eom\'etrique et le c\^ot\'e spectral
 de la formule des traces pour $G$ et $SI^G_{geo}$ le c\^ot\'e g\'eom\'etrique de la formule 
 des traces stables pour $G$ et $I^{\widetilde G_N }_{geo}$ $I^{\widetilde G_N }_{spec}$ les c\^ot\'es 
 de la formule des traces tordues; on  a supprim\'e le $V$ de la notation mais c'est bien pour
  cet ensemble de places fix\'e que ces objets sont d\'efinis.  
La premi\`ere remarque est $\tilde{f}_{m\nu}$ n'a de transfert non nul que pour la donn\'ee endoscopique
 associ\'ee \`a $\mathbf G({\mathbb R})$. Ainsi la stabilisation de la formule des traces pour 
 $\widetilde{\mathbf G}_N(\bbQ)=\GL_N({\mathbb Q}).\, \theta$
  se r\'eduit \`a:
$$
I^{\widetilde G_N }_{geo}(\tilde{f}(m\nu))=i(G)\;  SI^G_{geo}(f^G(m\nu)).
$$
Mais de m\^eme $f_{m\nu}$ est dans $SI_{cusp}(\mathbf G({\mathbb R}))$ et ses transferts aux groupes 
endoscopiques elliptiques propres de $\mathbf G({\mathbb R})$ sont nuls d'o\`u encore:
$$
SI^G_{geo}(f^G(m\nu))= I^G_{geo}(f^G(m\nu)).
$$
La fonction $f^G(m\nu)$ est cuspidale en au moins deux places; le c\^ot\'e g\'eom\'etrique de la formule 
des traces se r\'eduit donc aux int\'egrales orbitales. En plus $f_{m\nu}$ a ses int\'egrales orbitales 
non semi-simple nulles et ce c\^ot\'e g\'eom\'etrique se r\'eduit donc \`a une somme avec des coefficients
 explicites (des volumes) pour les int\'egrales semi-simples. La d\'ependance de ses int\'egrales quand $m$
  grandit est explicit\'e par Harish-Chandra comme le montre \cite{Art13}. Cela assure que $I^G_{geo}(f^G(m\nu))$
   est certainement non nul pour $m$ grand si la somme des int\'egrales orbitales en les \'el\'ements centraux de
    $\mathbf G({\mathbb Q})$ est non nulle. Or la fonction $f^G(m\nu)$ est invariante sous 
    l'action de ce centre par construction 
    et il suffit donc que $f^G(m\nu)(1)\neq 0$. Or il est clair que $1_{K^V}(1)=1$, $f_{m\nu}(1)$ n'est pas nul par le 
    calcul d'Harish-Chandra, $f_{0}(1)\neq 0$ car on l'a construit comme cela et pour tout $v\in V$ tel que $f_{v}$ est
     un pseudo coefficient d'une repr\'esentation de Steinberg, on a aussi $f_{v}(1)\neq 0$ par la formule de Plancherel. 
     D'o\`u la non nullit\'e cherch\'ee. 

Ainsi on a aussi une non nullit\'e du c\^ot\'e spectral:
$$
I^{\widetilde G_N }_{spec}(\tilde{f}(m\nu))= I^G_{spec}(f^G(m\nu))\neq 0,
$$
 pour $m$ suffisamment grand. Comme ${f}^G(m\nu)$ est non ramifi\'e hors de $V$, il existe un caract\`ere $c^V$ 
 comme ci-dessus tel que $\mathbf G({\mathbb A}_{{\mathbb Q}})$ ait des repr\'esentations non ramifi\'ees hors de $V$ se
  transformant suivant le caract\`ere $c^V$ pour l'alg\`ebre de Hecke sph\'erique et intervenant dans la partie discr\`ete 
  de la formule des traces pour $G$. En transf\'erant \`a $\widetilde{\mathbf G}_N$, on trouve 
  donc une repr\'esentation automorphe 
  $\tilde{\pi}(c^V)$ de $\GL_N({\mathbb A}_{{\mathbb Q}})$ invariante sous $\theta$ et dont la $\theta$-trace n'annule pas
   la fonction $\tilde{f}(m\nu)$. Par les th\'eor\`emes de multiplicit\'es un fort cette repr\'esentation est uniquement 
   d\'etermin\'e par le choix de $c^V$ ou plut\^ot son transfert $\tilde{c}^V$.
 
 Montrons que cette repr\'esentation a les propri\'et\'es annonc\'ees dans l'\'enonc\'e: d'abord on consid\`ere la
  place archim\'edienne. On note $\tilde{\pi}_{\infty}$ la composante de $\tilde{\pi}(c^V)$. On rappelle que 
  $\tilde{f}_{m\nu}$ est une fonction cuspidale. On \'ecrit la $\theta$-trace de $\tilde{\pi}_{\infty}$ dans
   le groupe de Grothendieck des $\theta$-repr\'esentations avec comme base les induites de
    repr\'esentations $\theta$-elliptiques. Et  $\tilde{f}_{m\nu}$ est de trace nulle sur les induites propres, il ne 
    reste que la contribution des $\theta$-elliptiques mais comme $\tilde{f}_{m\nu}$ est un pseudo coefficient de
     la repr\'esentation $\theta$-elliptique de caract\`ere infinit\'esimal $\tilde{m\nu}$  (le transfert de $m\nu$)le 
     caract\`ere infinit\'esimal de $\tilde{\pi}_{\infty}$ est aussi $\tilde{m\nu}$. Quand on \'ecrit $\tilde{\pi}(c^V)$ 
     comme induite de repr\'esentation de carr\'e int\'egrable, on voit que le fait que $m\nu$ soit loin 
     des murs force       $\tilde{\pi}_{\infty}$ \`a \^etre cette repr\'esentation $\theta$-elliptique. 
     Cela force en m\^eme temps $\tilde{\pi}(c^V)$ 
     \`a \^etre une induite de repr\'esentations cuspidales (et pas seulement de repr\'esentations de carr\'e int\'egrable).

 On consid\`ere maintenant une place $v\in V$ qui n'est ni l'infini ni $v_{0}$. Avec la m\^eme d\'emonstration que ci-dessus 
 on v\'erifie que $\tilde{\pi}(c^V)_{v}$ est un sous-quotient de la s\'erie principal qui admet aussi le transfert
  de la repr\'esentation   de Steinberg comme sous-quotient. En utilisant le fait que $\tilde{\pi}(c^V)_{v}$
   est n\'ecessairement une induite de composantes
   locales de repr\'esentations cuspidales (globales) cela force $\tilde{\pi}(c^V)_{v}$ a \^etre en plus temp\'er\'ee
    (bien que l'on ne     connaisse pas la conjecture de Ramanujan); il n'y a plus de choix si $\mathbf G\neq 
    \SO(2n)$, $\tilde{\pi}(c^V)_{v}$ est n\'ecessairement     la repr\'esentation de Steinberg. Et si $\mathbf G
    =\SO(2n)$ c'est n\'ecessairement l'induite $\mathrm{St}(2n)\times \eta_{v}$ o\`u $\eta_{v}$ est 
    le caract\`ere quadratique d\'eterminant la forme orthogonale.
 
 Il reste la place $v_{0}$; traitons le cas o\`u $\mathbf G=\Sp(2n)$ sinon il faut enlever 
 $\mu_{3}$ dans ce qui suit. En utilisant le fait que   $\tr \tilde{\pi}(c^V)_{v_{0}}(\tilde{f}_{0})\neq 0$, 
 on voit qu'il existe des caract\`eres non ramifi\'es $\mu_{1}, \mu_{2},\mu_{3}$ 
  tel que $\tilde{\pi}(c^V)_{v_{0}}=\omega\mu_{1}\times \omega^*\mu_{2}\times \mu_{3}$. Comme cette repr\'esentation est
   n\'ecessairement autoduale (en tenant aussi compte des propri\'et\'es de $\omega$), on a $\mu_{2}=\mu_{1}^{-1}$ et 
   $\mu_{3}^2=1$. Mais le caract\`ere central de cette repr\'esentation est n\'ecessairement trivial car il fait partie de la
    donn\'ee endoscopique li\'ee \`a $G$ et donc $\mu_{3}=1$. Ainsi en posant $\mu_{0}:=\mu_{1}$, on a montr\'e la
     totalit\'e des propri\'et\'es voulues.

 On \'ecrit la stabilisation de la partie discr\`ete de la formule des traces pour $ \widetilde{\mathbf G}_N$
  conform\'ement \`a \cite{stabilisationX}.   En fixant le caract\`ere $\tilde{c}^V$, dans ce c\^ot\'e spectral
   il n'y a que la trace tordue de $\tilde{\pi}(c^V)$. Dans la stabilisation, 
 a priori d'autres donn\'ees endoscopiques elliptiques que $G$ peuvent appara\^{\i}tre mais on ne calcule les
  distributions que sur   les fonctions $\tilde{f}_{V}^{\infty}1_{\tilde{K}^V}\tilde{f}_{m\nu}$ (avec les notations
   de l'\'enonc\'e). Le fait que les transferts de   $\tilde{f}_{m\nu}$ sont nuls sauf pour la donn\'ee 
   endoscopique elliptique li\'ee \`a $G$, simplifie la stabilisation en, pour tout 
 $\tilde{f}_{V}^{\infty}$
 \begin{equation}\label{C4}
 \tr_{\theta_N} \left( \tilde{\pi}(\tilde{c}^V,m\nu,\mu_{0})(\tilde{f}^{\infty}_{V}f_{m\nu}1_{\tilde{K}}^V)\right)
 = i(G) \sum_{\pi\in \Pi^G(c^V, m\nu)}m(\pi)\; \tr\left( \pi(f^{G,\infty}_{V} f_{m\nu}1_{K^V})\right), 
\end{equation}
o\`u $\Pi^G(c^V,m\nu)$ est l'ensemble des repr\'esentations automorphes de $G({\mathbb A})$ non
 ramifi\'ees hors de $V$ et y correspondant au caract\`ere $c^V$ de l'alg\`ebre de Hecke sph\'erique et
  ayant $m\nu$ comme caract\`ere infinit\'esimal; $m(\pi)$ est la multiplicit\'e dans la partie discr\`ete stable 
  du c\^ot\'e spectral. On a d\'ej\`a vu que l'on pouvait enlever le mot stable et $m(\pi)$ est donc la multiplicit\'e 
  dans la partie discr\`ete de la formule des traces pour $G$ et cette multiplicit\'e est positive 
  car la r\'egularit\'e du caract\`ere infinit\'esimal fait que l'on ne voit pas les termes venant 
  des paraboliques propres. Il reste juste \`a montrer 
   que l'on peut encore remplacer $\Pi^G(c^V,m\nu)$ par $\Pi^G(c^V,m\nu,\mu_{0})$: pour cela on remarque que 
   l'\'egalit\'e (\ref{C4}) est une \'egalit\'e de transfert en toute place de $V$ 
   sauf la place archim\'edienne. Mais on sait que 
   le terme de gauche est un transfert de la repr\'esentation irr\'eductible et stable 
   $\otimes_{v\in V-\{\infty,v_{0}\}}\mathrm{St}_{v}
    \otimes \Ind_{P}^{G}(\omega\mu_{0})$ (il n'y a pas \`a distinguer suivant les deux paraboliques possibles 
    si $\mathbf G=\SO(2n)$
     car le transfert des fonctions donnent une fonction invariante sous $\Or(2n)$). Avec l'ind\'ependance lin\'eaire des
      caract\`eres on obtient l'assertion cherch\'ee. Et ceci finit la d\'emonstration du th\'eor\`eme.\qed

\begin{cor} Le facteur de transfert $z({\psi})$ est égal à  $1$ d\`es que $m$ est suffisamment grand (o\`u $\psi$ est le
 param\`etre de Langlands du paquet de s\'eries discr\`etes de caract\`ere infinit\'esimal $m\nu$)
\end{cor}
\dem 
On reprend les d\'efinitions et notations du th\'ero\`eme pr\'ec\'edent. On remarque que pour tout $\tilde{f}_{V}^{\infty}$ 
on a l'\'egalit\'e de transfert:
$$
\tr_{\theta_N}\left( \otimes_{v\in V-\{\infty\}}\tilde{\pi}(c^V,m\nu,\mu_{0})(\tilde{f}_{V}^{\infty})\right)= \tr\left( 
\otimes_{v\in  V-\{\infty,v_{0}\}}
\mathrm{St}_{v}\otimes \Ind_{P}^{G}\omega\mu_{0}(f^{G,\infty}_{V})\right).
$$
On reporte dans le th\'eor\`eme en utilisant aussi le lemme fondamental hors de $V$ et on trouve finalement
$$
\tr_{\theta_N} \left( \tilde{\pi}(c^V,m\nu,\mu_{0})_{\infty}(\tilde{f}_{m\nu})\right)
=i(G)\sum_{\pi\in \Pi^{G}(c^V,m\nu,\mu_{0})}  m(\pi)\; \tr\left( \pi_{\infty}(f_{m\nu})\right).
$$
Et comme $f_{m\nu}$ est la somme des pseudo-coefficients pour toutes les s\'eries discr\`etes de $G({\mathbb R})$ 
ayant $m\nu$ comme caract\`ere infinit\'esimal, le terme de droit vaut exactement $i(G)\sum_{\pi\in 
\Pi^{G}(c^V,m\nu,\mu_{0})}m(\pi)$. 
Le terme de gauche vaut avec (\ref{rappel})  $z(\psi)^{-1}$. Comme $z(\psi)$ est de valeur absolue 1,
 cela force $z(\psi)=1$

\subsection{Le cas g\'en\'eral}

On reprend la m\^eme construction que pr\'ec\'edemment en suivant encore \cite{Art13} 6.2.2. Mais maintenant on a 
le r\'esultat aux places archim\'ediennes d\`es que le caract\`ere infinit\'esimal est suffisamment grand. 
On a toujours le r\'esultat aux places $p$-adiques si l'on ne consid\`ere que les repr\'esentations de Steinberg. 
On fixe donc une extension de ${\mathbb Q}$ ayant au moins deux places archim\'ediennes. On fixe $u$ l'une de
 ces places archim\'ediennes et on fixe un paquet de s\'eries discr\`etes en cette place index\'e par un param\`etre $\psi$; 
 on consid\`ere encore un pseudo coefficient dans $SI_{cusp}(\mathbf G({\mathbb R}))$ relatif \`a ce paquet et 
  on note $\psi$ 
 le param\`etre de Langlands du paquet. On fixe aussi $V$ un ensemble fini de places contenant $V_{ram}$
  (cf. \cite{stabilisationX}) et en ces places, quand elles sont $p$-adiques, on s'int\'eresse comme ci-dessus 
  \`a des repr\'esentations de Steinberg. Aux places archim\'ediennes, on a soit la place $u$ o\`u on a fix\'e le paquet
   de s\'eries discr\`etes soit des places o\`u encore une fois on s'autorise un caract\`ere infinit\'esimal 
   tendant vers l'infini     en \'etant loin des murs. 

\begin{thm} Le facteur de transfert $z(\psi)$ est égal à $1$.
\end{thm}
\dem C'est exactement la m\^eme d\'emonstration que ci-dessus. On n'a pas besoin de la place $v_{0}$ introduite 
ci-dessus et de la repr\'esentation induite en cette place puisqu'elle n'\'etait l\`a que pour ne pas devoir imposer 
aux paquets de s\'eries discr\`etes consid\'er\'es d'avoir un caract\`ere central trivial. Mais maintenant on a r\'egl\'e 
le cas de tous les paquets de s\'eries discr\`etes aux places archim\'ediennes d\`es que le caract\`ere infinit\'esimal
 est grand. \qed

\bigskip
\bibliographystyle{alpha}
\bibliography{AMR7}

\end{document}